\documentclass[a4paper,12pt]{article}
\usepackage[frenchb]{babel}
\usepackage[latin1]{inputenc}
\usepackage[T1]{fontenc}
\usepackage{amsmath}
\usepackage{amssymb}
\usepackage{pdfsync}
\usepackage{url}
\newtheorem{defi}{Définition}[section]
     \newtheorem{rem}[defi]{Remarque}
      \newtheorem{prop}[defi]{Proposition}
      \newtheorem{thm}[defi]{Théorème}
   
      \newtheorem{lem}[defi]{Lemme}

\newcommand {\tu}[1]{{}^{\varphi^{#1}}\!}

\def\ta {{}^\varphi}
\def\cqfd{{}\hfill $\square$ \vskip 1mm}

\def\Fil{{\mathrm {Fil}}}
\def\tc{{\widehat \otimes }}

\def\Hom{{\mathrm{Hom}}}
\def\R{{\mathbb R}}

\def\Z{{\mathbb Z}}
\def\N{{\mathbb N}}

\def\CC{\mathcal{C}}
\def\NN{\mathcal{N}}
\def\QQ{{\mathbb Q}}
\def\M{\mathcal{M}}

\def\O{\mathcal{O}}

\def\A{\mathcal{A}}

\def\F{\mathbb{F}_q}
\def\Fp{\mathbb{F}_p}
\def\V{\mathcal{V}}

\def\gS{\mathfrak{S}}

\def\Id{\mathrm{Id}}
\def\mod{\mathrm{mod}}
\def\det{\mathrm{det}}
\def\W{\mathcal{W}}

\def\vi{\varphi}

\def\naif{\text{naïf}}

\newcommand{\s}[1]{\langle #1 \rangle}

\def\wt{\widetilde }

\def\dotss{...}

\def\gM{\mathfrak{M}}
\def\gN{\mathfrak{N}}
\def\wS{\widehat \gS}
\def\Fil{\mathrm{Fil}}
\def\ModphiS{\mathrm{Mod}^{\vi}_{/\gS}}
\def\PIModphiS{\mathrm{PIMod}^{\vi}_{/\gS}}
\RequirePackage{graphics}

\begin{document}

\title{Structures de Hodge-Pink  pour les $\vi / \gS$-modules de Breuil et Kisin}
\author{Alain Genestier et Vincent Lafforgue}

\maketitle

Le but de cet article est d'appliquer les méthodes de~\cite{fontaine29}  aux $\vi / \gS$-modules de Breuil et Kisin (voir~\cite{breuil-ENS97, breuil-griffiths, breuil-nonpub, breuil-compositio, breuil-invent99, breuil,kisin}). 
On  démontre ainsi un léger renforcement  du corollaire 1.3.15 de~\cite{kisin}, sans utiliser les  résultats de Kedlaya~\cite{kedlaya,slopefil}. Le corollaire 1.3.15 de~\cite{kisin}  est un étape essentielle de la nouvelle démonstration 
du théorème ``faiblement admissible implique admissible'' de Colmez et Fontaine~\cite{colmez-fontaine} que Kisin a donnée dans~\cite{kisin}. 
En insérant dans~\cite{kisin} notre preuve du corollaire  1.3.15 
on obtient donc une démonstration relativement élémentaire du  théorème ``faiblement admissible implique admissible''. 

En revanche   les méthodes employées ici ne s'appliquent pas aux $(\varphi,\Gamma)$-modules 
(voir~\cite{fontaine-groth,berger-invent02,berger,berger-intro,kisin-ren}), pour lesquels
  le relèvement du Frobenius n'est plus $u\mapsto u^{p}$. 
De fa\c con plus précise il semble que l'on ne puisse pas retrouver avec notre méthode les résultats de~\cite{kisin-ren} (qui sont pourtant assez parallèles à ceux de~\cite{kisin} et utilisent les résultats de Kedlaya de la même fa\c con). 
  
 La théorie des $\vi / \gS$-modules de Breuil et Kisin est $\Z_{p}$-linéaire. 
 Dans le dernier paragraphe nous indiquerons un cadre général où $\Z_{p}$ est remplacé par l'anneau des entiers d'un corps local non archimédien $\O_{F}$. Les $\vi / \gS$-modules en ce sens généralisé sont exactement les chtoucas locaux~\cite{fontaine29} lorsque $\O_{F}$ est d'égales caractéristiques. Ce cadre plus général unifie une partie de~\cite{fontaine29} et le présent article (à l'exception du premier paragraphe), à des variantes de notations près.  
 
 Soit $k$ un corps parfait contenant $\Fp$. On note $W=W(k)$.  On possède un morphisme de Frobenius $\Z_p$-linéaire $\vi:W\to W$ qui relève l'endomorphisme $x\mapsto x^{p}$ de $k$.  On note $K_{0}=\mathrm{Frac} W=W[\frac{1}{p}]$. 
On note 
$\gS=W[[u]]$. 
On note $\vi:\gS\to \gS$ le morphisme égal à $\vi$ sur $W$ et envoyant $u$ sur $u^{p}$. Si $M$ est un module sur un anneau $A$ muni de $\vi:A\to A$
 (qui sera égal à $W$ ou $\gS$ ou à un localisé d'un de ces anneaux) on note $\ta M=M\otimes_{A,\vi}A$. Ce module 
 $\ta M$ est noté $\vi^{*}(M)$ dans~\cite{kisin} et ${}^{\tau}M$ dans~\cite{fontaine29}. Si $f:M\to M'$ est un morphisme de $A$-modules on note 
 $\ta f=f \otimes 1$ et si $x$ est un élément de $M$ on note $\ta x=x \otimes 1$.

Soit $K$ une extension totalement ramifiée de $K_{0}=W[\frac{1}{p}]$, $e=[K:K_{0}]$ l'indice de ramification de $K$ sur $K_{0}$, 
 $\pi_{K}$  une uniformisante de $K$,  
et $E$ le polynôme minimal de $\pi_{K}$ sur $K_{0}$, qui est un polynôme d'Eisenstein.   On a donc $E=u^{e}+c_{e-1}u^{e-1}+...+c_{0}$, avec $c_{e-1},...,c_{1}\in pW$ et $c_{0}\in pW^{\times}$. On notera $E(0)=c_{0}$. 

Soit $\O_{K}$ l'anneau des entiers de $K$. Le morphisme $\gS/E\gS\to \O_{K}$ qui envoie $u$ sur $\pi_{K}$ est un isomorphisme.  

\begin{defi} \label{def-ModphiS}
On appelle $\vi / \gS$-module un couple $(\gM, \vi_{\gM})$ où $\gM$ est un $\gS$-module libre de type fini et $\vi_{\gM}: \ta \gM[\frac{1}{E}]\to \gM[\frac{1}{E}]$ est un isomorphisme de $\gS[\frac{1}{E}]$-modules. 
 Pour $s,t\in \Z$ vérifiant $s\leq t$, on dit que $(\gM, \vi_{\gM})$ est d'amplitude $\subset [s,t]$ si on a 
 $E^{t}\ \gM \subset \vi_{\gM}(\ta \gM)\subset E^{s}\ \gM$.
 
 On note $\ModphiS$  
 la catégorie des $\vi / \gS$-modules et $\mathrm{Mod}^{\vi}_{/\gS,[s,t]} $  
 la catégorie des $\vi / \gS$-modules d'amplitude $\subset [s,t]$.

 On note $\ModphiS \otimes_{\Z_p}\QQ_p$ la catégorie à isogénie près et on appelle iso-$\vi / \gS$-module un objet de cette catégorie. 
 \end{defi}

\begin{rem} La catégorie $\ModphiS$ de~\cite{kisin} n'est pas exactement la nôtre mais en est la sous-catégorie pleine formée des objets d'amplitude $\subset [0,+\infty[$. 
\end{rem}

Soit $\wS$ le complété de $\gS[\frac{1}{p}]$ par rapport à l'idéal engendré par $E$. C'est un anneau local d'uniformisante $E$ dont le corps résiduel est $K=\gS[\frac{1}{p}]/E\gS[\frac{1}{p}]=\widehat \gS/E\widehat \gS$. 

\begin{defi}
On appelle $\vi$-module  un couple  
$D=(D,\vi_{D})$ où 
\begin{itemize}
\item $D$ est un $K_{0}$-espace vectoriel de dimension finie
\item $\vi_{D}:\ta D\to D$ est un isomorphisme de $K_{0}$-espaces vectoriels. 
\end{itemize}

Si $D$ est un $K_{0}$-espace vectoriel de dimension finie on appelle structure de Hodge-Pink sur $D$ un $\widehat \gS$-module libre $V$ qui est un réseau 
 dans $\ta D\otimes_{K_{0}}\widehat \gS[\frac{1}{E}]$. 

On appelle $\vi$-module de Hodge-Pink un triplet $D=(D,\vi_{D},V_{D})$ où 
\begin{itemize}
\item $(D,\vi_{D})$ est un $\vi$-module, 
\item $V_{D}$ est une structure de Hodge-Pink sur $D$. 
  \end{itemize}
  On notera toujours 
$U_{D}=\ta D\otimes_{K_{0}}\widehat \gS$ qui est un  réseau 
dans $\ta D\otimes_{K_{0}}\widehat \gS[\frac{1}{E}]$. 

On note $MHP(\vi)$ la catégorie des $\vi$-modules de Hodge-Pink.
Un morphisme de $\vi$-modules de Hodge-Pink  $D'\to D$ est un morphisme de $\vi$-modules $f$ tel que $f(V_{D'})\subset V_{D}$. 
 \end{defi}

La terminologie ``Hodge-Pink'' provient du cas d'égales caractéristiques (voir~\cite{pink,pink-expose,fontaine29}). De tels objets ont été introduits par Breuil~\cite{breuil-ENS97, breuil-griffiths,   breuil-invent99, breuil} mais toujours sous la condition de transversalité de Griffiths  (que nous rappelerons dans \eqref{cond-gr} ci-dessous). 

Si $D$ est de dimension $1$, la matrice de $\vi_D$ dans une base  est le produit de $p^k$ et d'une
unité. De plus on a $V_D=E^{l}U_D$. On note $t_N(D)=k$ et
$t_H(D)=-l$. 

Si $D$ est de dimension $r$, on note $\det(D)=\Lambda^rD$,
$\vi_{\det(D)}=\det(\vi_D)$ et $V_{\det(D)}=\det(V_D)$. Alors $\det(D)$
est un  $\vi$-module de Hodge-Pink  de dimension $1$. 
On pose alors $t_N(D)=t_N(\det(D))$ et
$t_H(D)=t_H(\det(D))$. 

On notera $t_N(D)=t_N(D,\vi_{D})$ et $t_H(D)=t_H(D,\vi_{D},V_{D})$ en cas d'ambiguïté. 

Si $(D,\vi_{D})$ est un $\vi$-module on appelle sous-$\vi$-module un
sous-$K_{0}$-espace vectoriel $D'$ tel que $\vi_D(\ta D')=D'$, et on prend pour $\vi_{D'}$ la restriction de $\vi_D$ à $\ta D'$.
 Si $D$ est
un $\vi$-module de Hodge-Pink  on munit $D'$ d'une structure de
 Hodge-Pink en posant $V_{D'}=V_D\cap (\ta D'\otimes _{K_{0}}\wS[\frac{1}{E}])$. 

\begin{defi}\label{def-fa-HP}
Un $\vi$-module de Hodge-Pink  $D$ est dit faiblement admissible si 
$t_H(D)=t_N(D)$ et pour tout sous-$\vi$-module $D'$ on a  
$t_H(D')\leq t_N(D')$. 
\end{defi}
\begin{rem} Soit $D$ faiblement admissible. 
Pour tout $\vi$-module de Hodge-Pink  $(D',\vi_{D'},V_{D'})$ et tout  morphisme injectif $f:D'\to D$ de $\vi$-modules de Hodge-Pink, on a $V_{D'}\subset V_D\cap (\ta D'\otimes _{K_{0}}\wS[\frac{1}{E}])$ et donc  $t_H(D')\leq t_N(D')$. \end{rem}

On va construire un foncteur de Dieudonné $\mathbb D_{\mathrm{iso}}: \ModphiS \otimes_{\Z_p}\QQ_p\to MHP(\vi) $ et notre résultat principal sera le théorème suivant. 

\begin{thm}\label{fa--a}
Le foncteur $\mathbb D_{\mathrm{iso}}$ est pleinement fidèle et son image essentielle est constituée exactement des $\vi$-modules de Hodge-Pink faiblement admissibles. 
\end{thm}

On appellera admissibles les objets de l'image essentielle  de ce foncteur. 
 La partie difficile de ce théorème est bien sûr l'implication ``faiblement admissible implique admissible''. Le théorème~\ref{fa--a} sera démontré dans le paragraphe~\ref{faibl-ad=ad} (en utilisant les résultats  des paragraphes~\ref{pseudo-iso-phi-gS} et~\ref{construc-Dieudonne}). 
 
  Le théorème~\ref{fa--a} est un léger renforcement du corollaire 1.3.15 de~\cite{kisin}  car on verra dans le premier paragraphe que 
    l'énoncé de Kisin équivaut au cas particulier du théorème~\ref{fa--a} où  les structures de Hodge-Pink vérifient la condition de transversalité de  Griffiths (rappelée  dans \eqref{cond-gr} ci-dessous). 
  Il est très probable que les arguments du premier paragraphe de~\cite{kisin} (qui reposent de fa\c con essentielle sur les résultats de Kedlaya) permettent aussi de montrer notre théorème~\ref{fa--a}. L'avantage de la démonstration que nous proposons ici est d'être plus élémentaire et de ne pas nécessiter l'introduction d'une clôture algébrique de $K$. 
  
  \vskip 1mm
  Voici le contenu de cet article. 
  
  Dans le premier paragraphe  nous donnons un dictionnaire entre 
 les $(\vi,N)$-modules filtrés au sens de Fontaine et 
  les $\vi$-modules de Hodge-Pink vérifiant la condition de transversalité de Griffiths  (rappelée dans \eqref{cond-gr} ci-dessous).  Les propriétés de faible admissibilité se correspondent de part et d'autre. Ceci fait du corollaire 1.3.15 de~\cite{kisin} une conséquence du 
  théorème~\ref{fa--a}.  

Dans le paragraphe~\ref{pseudo-iso-phi-gS} nous introduisons la notion de pseudo-iso-$\vi / \gS$-module et nous montrons qu'elle est équivalente à la notion  
  plus restrictive d'iso-$\vi / \gS$-module. 
  
   Dans le paragraphe~\ref{construc-Dieudonne} nous construisons un foncteur de Dieudonné de la catégorie des  pseudo-iso-$\vi / \gS$-modules vers celle des $\vi$-modules de Hodge-Pink. Nous montrons de plus qu'il est pleinement fidèle, que son image essentielle est formée d'objets faiblement admissibles et qu'elle est stable par extension.

  Dans le paragraphe~\ref{faibl-ad=ad} nous achevons la démonstration du théorème~\ref{fa--a}, qui est notre résultat principal. 
  
  Dans le paragraphe~\ref{fontainelaffaille1} nous développons une théorie entière : par des arguments très proches de ceux de la démonstration du théorème~\ref{fa--a}, nous redémontrons un résultat de Caruso et Liu
 (le  théorème 2.2.1 de~\cite{carliu}).  
  
   Enfin le dernier paragraphe introduit un cadre général où $\Z_{p}$ est remplacé par l'anneau des entiers d'un corps local non archimédien.

Cet article est essentiellement une transcription de~\cite{fontaine29}, aux changements de notations près. Au prix de certaines répétitions, nous avons fait en sorte que le lecteur n'ait pas besoin de se reférer à~\cite{fontaine29}.

Nous remercions Laurent Fargues pour des observations qui ont été déterminantes pour cet article. 
Nous remercions également Laurent Berger, Christophe Breuil, Olivier Brinon, Xavier Caruso, Pierre Colmez, Jean-Marc Fontaine, Mark Kisin et Ariane Mézard pour des discussions très utiles et leurs encouragements à écrire cet article. 

Dans cet article on notera souvent $f$ au lieu de $f\otimes 1$. Si $A$ est anneau, $M_{r}(A)$ désignera l'anneau des matrices de taille $r\times r$ à coefficients dans $A$. 

 \section{Transversalité de Griffiths}\label{sec-griffiths}
  
Le but de ce paragraphe est de montrer que  le théorème~\ref{fa--a}  implique le corollaire 1.3.15 de~\cite{kisin}. 

La définition ci-dessous est due à Fontaine~\cite{fontaine3}. 

\begin{defi}
On appelle $(\vi,N)$-module filtré  un quadruplet $$D=(D,\vi_{D},N_{D},(\Fil^{i}D_{K}))\text{ \ \ où  }$$ 
\begin{itemize}
\item $(D,\vi_{D})$ est un $\vi$-module, 
\item $N_{D}:D\to D$ est une application $K_{0}$-linéaire vérifiant 
$$N_{D}\vi_{D}=p\vi_{D}\ta N_{D}, $$
\item $ (\Fil^{i}D_{K})$ est une filtration décroissante de $D_{K}=D\otimes_{K_{0}}K$ par des sous-$K$-espaces vectoriels tels que 
$$\Fil^{i}D_{K}=D_{K}\text{ pour }i<<0\text{ et }\Fil^{i}D_{K}=0\text{ pour }i>>0.$$ 
\end{itemize}
On note $MF(\vi,N)$ la catégorie des $(\vi,N)$-modules filtrés.
 \end{defi}

Soit $D=(D,\vi_{D},N_{D},(\Fil^{i}D_{K}))$ un $(\vi,N)$-module filtré. On définit 
$t_{N}(D)$ comme précédemment. Si $D$ est de dimension $1$, on note 
$t_{H}(D)$ l'unique entier $i$ tel que $\Fil^{i}D_{K}=D_{K}$ et $\Fil^{i+1}D_{K}=0$. En général on note $t_{H}(D)=t_{H}(\det D)$. 

Si $(D,\vi_{D})$ est un $\vi$-module on appelle sous-$(\vi,N)$-module un
sous-$K_{0}$-espace vectoriel $D'$ tel que $\vi_D(\ta D')=D'$ et 
$N_{D}(D')\subset D'$. On définit alors 
$\vi_{D'}$  et $N_{D'}$ comme la restriction de $\vi_D$ à $\ta D'$ et de $N_{D}$ à $D'$ et   on munit $D'$ de la filtration  $\Fil^{i}D'_{K}=\Fil^{i}D_{K}\cap ( D'\otimes _{K_{0}}K)$. 

D'après~\cite{fontaine3} on a la définition suivante. 

\begin{defi}\label{def-f-a-fontaine}
Un $(\vi,N)$-module filtré   $(D,\vi_{D},N_{D},(\Fil^{i}D_{K}))$ est faiblement admissible si $t_{H}(D)=t_{N}(D)$ et si pour tout sous-$(\vi,N)$-module $D'$ on a $t_{H}(D')\leq t_{N}(D')$. 
\end{defi}

A  toute structure de Hodge-Pink $V$ (c'est-à-dire un $\wS$-réseau dans 
$\ta D\otimes _{K_0}\wS[\frac{1}{E}]$) on peut associer la structure de Hodge $(\Fil ^{i}D_{K})_{i\in \Z}$ (c'est-à-dire  une filtration décroissante, séparée et exhaustive, de $D_{K}=D\otimes _{K_0}K$ par des sous-$K$-espaces vectoriels) définie par 
\begin{gather}\label{h-hp}
\vi_{D}^{-1}(\Fil^{i}D_{K})=\Big(E^{i}V\cap U_{D}\Big)/\Big(E^{i}V\cap EU_{D}\Big)\end{gather}
dans $\ta D\otimes_{K_{0}} K=U_{D}/ EU_{D}$, pour tout $i\in \Z$.  

 Cette application (qui à une structure de Hodge-Pink associe une structure de Hodge) est surjective mais pas injective (sauf dans le cas minuscule). 
 Le lemme suivant, qui est dû à Breuil (voir~\cite{breuil-ENS97,breuil-griffiths,breuil-invent99,breuil}),  fournit une section de cette application (dépendant de  la donnée supplémentaire de $N_{D}$) 
  
  \begin{lem}\label{unicite-1-et-2}
  Soit  $D=(D,\vi_{D},N_{D},(\Fil^{i}D_{K}))$ un 
 $(\vi,N)$-module filtré.
  Il existe une unique structure de Hodge-Pink $V$ sur $D$ qui vérifie les deux conditions suivantes
\begin{itemize}
\item la filtration $(\Fil^{i}D_{K})$ est associée à $V$ comme dans \eqref{h-hp},   
\item $V$ vérifie la condition de transversalité de Griffiths relativement à la connexion $p\ta N_{D}\otimes 1+1\otimes u\frac{d}{du}$ sur $\ta D\otimes_{K_{0}}\widehat \gS[\frac{1}{E}]$, c'est-à-dire que 
\begin{gather}\label{cond-gr}\big(p\ta N_{D}\otimes 1+1\otimes u\frac{d}{du}\big)(V)\subset E^{-1}V.\end{gather} 
\end{itemize}
    \end{lem}
   \noindent {\bf Démonstration (d'après Breuil~\cite{breuil-griffiths,breuil}).}
     La formule suivante détermine $V$ par récurrence. Pour $i\in \Z$ assez petit on a $ E^{i}V\cap U_{D}=U_{D}$ et pour tout $i\in \Z$ on a 
     \begin{gather}\nonumber  E^{i+1}V\cap U_{D}=\{x\in U_{D},  \big(p\ta N_{D}\otimes 1+1\otimes u\frac{d}{du}\big)(x)\in E^{i}V\cap U_{D} \\\label{form-rec-Fil-i-V} \text{ et } x\text{ mod } EU_{D}\in 
    \vi_{D}^{-1}(\Fil^{i}D_{K}) \text{ dans } U_{D}/EU_{D}=\ta D\otimes _{K_0}K\}. \end{gather} De plus $V=E^{-i}(E^{i}V\cap U_{D})$ pour $i$ assez grand. 
           \cqfd

On note  
$$HP:  MF(\vi,N)\to MHP(\vi) $$  le foncteur qui à $(D,\vi_{D},N_{D},(\Fil^{i}D_{K}))\in MF(\vi,N)$
associe  $(D,\vi_{D},V_{D})$ où $V_{D}$ est l'unique structure de 
Hodge-Pink  qui satisfait les deux conditions \eqref{h-hp} et \eqref{cond-gr}.

On a alors le lemme suivant. 

\begin{lem}\label{lem-faHP-faFil}
L'image d'un objet  $(D,\vi_{D},N_{D},(\Fil^{i}D_{K}))$ de $ MF(\vi,N)$ par $HP$ est faiblement admissible  au sens de la définition~\ref{def-fa-HP}  
si et seulement si $(D,\vi_{D},N_{D},(\Fil^{i}D_{K}))$ est faiblement admissible au sens de la définition~\ref{def-f-a-fontaine}. 
\end{lem}

\noindent{\bf Démonstration.} On note $(D,\vi_{D},V_{D})$ l'image de  
$(D,\vi_{D},N_{D},(\Fil^{i}D_{K}))$ par $HP$. Il est évident que la faible admissibilité de $(D,\vi_{D},V_{D})$ implique la faible admissibilité de 
$(D,\vi_{D},N_{D},(\Fil^{i}D_{K}))$. En effet on a 
\begin{gather}\label{eg-tH-tH}
t_{H}(D,\vi_{D},V_{D})=
t_{H}(D,\vi_{D},N_{D},(\Fil^{i}D_{K}))
\end{gather} 
et cette relation reste vraie pour tout sous-$(\vi,N)$-module $D'$ de $D$. 

On suppose maintenant que $(D,\vi_{D},N_{D},(\Fil^{i}D_{K}))$ est faiblement admissible et on suppose par l'absurde que 
$(D,\vi_{D},V_{D})$ ne l'est pas. Il existe donc un sous-$\vi$-module $D'$ dans $D$ tel qu'en posant $V_{D'}=V_D\cap (\ta D'\otimes _{K_{0}}\wS[\frac{1}{E}])$ on ait 
$t_{H}(D',\vi_{D'},V_{D'})>t_{N}(D',\vi_{D'})$. On suppose $D'$ minimal pour cette propriété, ce qui implique pour tout $\vi$-module $\overline D'$ quotient non trivial de $D'$, en notant $V_{\overline D'}$ l'image de $V_{D'}$ dans 
$\ta{\overline D}'\otimes _{K_{0}}\wS[\frac{1}{E}]$ on a
$t_{H}({\overline D}',\vi_{{\overline D}'},V_{{\overline D}'})>t_{N}({\overline D}',\vi_{{\overline D}'})$. Pour tout $i\in \N$ on note $$D'_{i}=D'+N_{D}(D')+...+N_{D}^{i}(D').$$ On a $D'_{0}=D'$ et pour tout $i\in \N$, $D'_{i}$ est un sous-$\vi$-module de $D$. La suite $(D'_{i})_{i\in \N}$ est croissante et stationnaire et on note $D'_{\infty}$ sa limite, qui est un sous-$(\vi,N)$-module de $D$. On pose $D'_{-1}=0$.  Pour tout $i\in \N$
 on note $\overline D'_{(i)}$ le  quotient de $D'$ par le noyau du morphisme surjectif $D'\xrightarrow{N_{D}^{i}}D'_{i}/D'_{i-1}$, autrement dit 
$\overline D'_{(i)}$ est le  $\vi$-module quotient de $D'$ tel que $$ {\overline D}'_{(i)}\xrightarrow{N_{D}^{i}}D'_{i}/D'_{i-1}$$ soit un isomorphisme de $K_{0}$-espaces vectoriels. Comme $N_{D}\vi_{D}=p\vi_{D}\ta N_{D}, $ on a 
\begin{gather}\label{diff-tN-Ni}
t_{N}(D'_{i}/D'_{i-1})=t_{N}({\overline D}'_{(i)})-i\dim({\overline D}'_{(i)}).\end{gather}
Pour tout $i\in \N$ on note $V_{D'_{i}}=V_D\cap (\ta D'_{i}\otimes _{K_{0}}\wS[\frac{1}{E}])$ et $V_{D'_{i}/D'_{i-1}}$ l'image de $V_{D'_{i}}$ dans $
\ta (D'_{i}/D'_{i-1})\otimes _{K_{0}}\wS[\frac{1}{E}]$. 
On remarque que l'on a une suite exacte 
$
0\to V_{D'_{i-1}}\to V_{D'_{i}}\to V_{D'_{i}/D'_{i-1}}\to 0 $
d'où \begin{gather}\label{eg-suite-exacte-VD-i-i-1}
t_{H}(D'_{i},\vi_{D'_{i}},V_{D'_{i}}) -t_{H}(D'_{i-1},\vi_{D'_{i-1}},V_{D'_{i-1}})=t_{H}(D'_{i}/D'_{i-1},\vi_{D'_{i}/D'_{i-1}},V_{D'_{i}/D'_{i-1}}). 
\end{gather}
D'autre part il est évident que 
\begin{gather}\label{eg-tN-i-i-1}
t_{N}(D'_{i},\vi_{D'_{i}}) -t_{N}(D'_{i-1},\vi_{D'_{i-1}})=t_{N}(D'_{i}/D'_{i-1},\vi_{D'_{i}/D'_{i-1}}). 
\end{gather}
La condition \eqref{cond-gr} de transversalité de Griffiths implique que pour tout $i\in \N$, on a 
$V_{D'_{i}}\supset \Big(E\big( p\ta N_{D}\otimes 1+1\otimes u\frac{d}{du} \big)  \Big)^{i}(V_{D'_{0}})$, d'où $V_{D'_{i}/D'_{i-1}}\supset E^{i}\ \ta N_{D}^{i}(V_{D'_{0}})$, ce qui équivaut à $V_{D'_{i}/D'_{i-1}}\supset E^{i}\ \ta N_{D}^{i}(V_{{\overline D}'_{(i)}})$. On a donc \begin{gather}\label{diff-tH-Ni}t_{H}(D'_{i}/D'_{i-1},\vi_{D'_{i}/D'_{i-1}},V_{D'_{i}/D'_{i-1}}) \geq t_{H}({\overline D}'_{(i)},\vi_{{\overline D}'_{(i)}},V_{{\overline D}'_{(i)}})-i\dim({\overline D}'_{(i)}).\end{gather}
Par hypothèse on a 
$t_{H}({\overline D}'_{(i)},\vi_{{\overline D}'_{(i)}},V_{{\overline D}'_{(i)}})>t_{N}({\overline D}'_{(i)},\vi_{{\overline D}'_{(i)}})$ dès que ${\overline D}'_{(i)}$  est non trivial. Il résulte alors de \eqref{diff-tN-Ni} et \eqref{diff-tH-Ni} que pour tout $i\in \N$ tel que $D'_{i}/D'_{i-1}$ est non trivial on a 
\begin{gather}\label{ineg-i/i-1}
t_{H}(D'_{i}/D'_{i-1},\vi_{D'_{i}/D'_{i-1}},V_{D'_{i}/D'_{i-1}})>
t_{N}(D'_{i}/D'_{i-1},\vi_{D'_{i}/D'_{i-1}}). 
\end{gather}
Alors \eqref{ineg-i/i-1}, \eqref{eg-suite-exacte-VD-i-i-1} et \eqref{eg-tN-i-i-1} impliquent que 
\begin{gather}\label{ineg-infty}
t_{H}(D'_{\infty},\vi_{D'_{\infty}},V_{D'_{\infty}})>
t_{N}(D'_{\infty},\vi_{D'_{\infty}}). 
\end{gather}
Comme $D'_{\infty}$ est un sous-$(\vi,N)$-module de $D$ (et non simplement un sous-$\vi$-module de $D$), et grâce à l'égalité analogue à \eqref{eg-tH-tH} pour $D'_{\infty}$, 
ceci contredit la faible admissibilité de $(D,\vi_{D},N_{D},(\Fil^{i}D_{K}))$. 
\cqfd

Grâce au lemme~\ref{lem-faHP-faFil}, le théorème~\ref{fa--a} implique le corollaire 1.3.15 de~\cite{kisin}. 

Voici un exemple pour illustrer les relations entre les conditions de faible admissibilité pour les 
$(\vi,N)$-modules filtrés et les $\vi$-modules de Hodge-Pink. 
On considère le $\vi$-module $(D,\vi_{D})$ donné par $D=K_{0}^{2}$ et 
$\vi_{D}=\Id$. On notera que l'on a forcément $N_{D}=0$. On note $(e_{1},e_{2})$ la base canonique de $K_{0}^{2}$. Soit $(\Fil^{i}D_{K})$ la filtration définie par 
$$\Fil^{-1}D_{K}=D_{K},\ \  \Fil^{0}D_{K}=\Fil^{1}D_{K}=Ke_{1}, \ \ \Fil^{2}D_{K}=0.$$
Le $(\vi,N)$-module filtré  $(D,\vi_{D},N_{D},(\Fil^{i}D_{K}))$ n'est pas faiblement admissible au sens de Fontaine à cause du sous-$(\vi,N)$-module $D'=K_{0}e_{1}$. Soit maintenant $\alpha \in K$ et soit $V_{D}^{\alpha}$ la structure de Hodge-Pink définie par $$V_{D}^{\alpha}= \wS E^{-1}(e_{1}+\alpha E e_{2})+\wS E e_{2}$$
de sorte que pour tout $\alpha \in K$, $(\Fil^{i}D_{K})$ est la filtration de Hodge associée à $V_{D}^{\alpha}$ par \eqref{h-hp}. 
Alors le $\vi$-module de Hodge-Pink $(D,\vi_{D}, V_{D}^{\alpha})$ est faiblement admissible si et seulement si $\alpha \neq 0$. On remarque aussi que $V_{D}^{\alpha}$ vérifie la condition de transversalité de Griffiths si et seulement si $\alpha=0$.

\section{Pseudo-iso-$\vi / \gS$-modules}\label{pseudo-iso-phi-gS}

\begin{defi}\label{defiisochtouca}
Pour $s,t\in \Z$ vérifiant $s\leq t$,
 on appelle pseudo-iso-$\vi / \gS$-module d'amplitude $\subset [s,t]$ un couple 
$(\gN,\vi_{\gN})$  où $\gN$ est un $\gS[\frac{1}{p}]$-module libre de type fini et $\vi_{\gN}: \ta \gN[\frac{1}{E}]\to \gN[\frac{1}{E}]$ est un isomorphisme vérifiant la condition

\noindent {\bf (PIM)} : 
pour un (ou pour tout) $\gS$-module libre 
$\gM$ muni d'un isomorphisme 
$\gM[\frac{1}{p}]=\gN$, il existe une constante $C$, telle que pour tout $n\in \N^{*}$, 
\begin{gather}\nonumber \vi_{\gN} \ta\vi_{\gN} \dotss  \tu{n-1}\vi_{\gN}\in p^{-C}E^{s}\dotss  \vi^{n-1}(E)^{s}
 \mathrm{Hom}_{\gS}
(\tu{n}\gM,\gM)\\ 
\nonumber\text{ et \ \ \ \  }(\vi_{\gN} \ta\vi_{\gN} \dotss  \tu{n-1}\vi_{\gN} )^{-1}\in p^{-C}E^{-t}\dotss  \vi^{n-1}(E)^{-t}
 \mathrm{Hom}_{\gS}
(\gM,\tu{n} \gM).\end{gather}
 On note $\PIModphiS$  
 la catégorie des pseudo-iso-$\vi / \gS$-modules. 
\end{defi}

\begin{prop}\label{equivalence}
Le foncteur évident 
$(\gM,\vi_{\gM})\mapsto (\gM[\frac{1}{p}],\vi_{\gM}\otimes 1)$ de 
 la catégorie $\ModphiS \otimes_{\Z_p}\QQ_p$ des iso-$\vi / \gS$-modules dans la catégorie $\PIModphiS$ des pseudo-iso-$\vi / \gS$-modules est une équivalence de catégories. 
\end{prop}
\noindent {\bf Démonstration.}
Seule l'essentielle surjectivité doit être démontrée. Soient $s,t\in \Z$. 
Il s'agit d'établir qu'un pseudo-iso-$\vi / \gS$-module $\gN$ d'amplitude $\subset [s,t]$  est associé à un $\vi / \gS$-module $\gM$ d'amplitude $\subset [s,t]$ (évidemment non unique, mais unique à 
isogénie près). On  note $\O_{\mathcal E}$ le complété de $\gS[\frac{1}{u}]$ pour la topologie $p$-adique. C'est un anneau de valuation discrète complet  d'uniformisante $p$ dont le corps résiduel s'identifie à $k((u))$. On note $\mathcal E=\O_{\mathcal E}[\frac{1}{p}]$ son corps des fractions. 
Par le lemme~\ref{langtonOL} ci-dessous  il suffit de trouver un $\O_{\mathcal E}$-réseau $V$ dans 
le $\mathcal E$-espace vectoriel $\gN\otimes _{\gS[\frac{1}{p}]}\mathcal E$, qui soit préservé par $\vi _{\gN}
\otimes _{\gS[\frac{1}{E},\frac{1}{p}]}\Id_{\mathcal E}$
(c'est-à-dire $\vi _{\gN}
\otimes _{\gS[\frac{1}{E},\frac{1}{p}]}\Id_{\mathcal E}(\ta V)=V$). 
En effet si $\gM$ est le $\gS$-module libre associé à $\gN$ et $V$ comme dans le lemme~\ref{langtonOL},  on a 
$\vi_{\gM}(\ta \gM)\subset E^{s}\gM$ et $\vi_{\gM}^{-1}(\gM)\subset E^{-t}\  \ta \gM$ car $E^{s}\gS[\frac{1}{p}]\cap \O_{\mathcal E} =E^{s} \gS$ (et de même avec $-t$).

Montrons maintenant l'existence de $V$. 
Soit $V_{1}$ n'importe quel $\O_{\mathcal E}$-réseau de $\gN\otimes _{\gS[\frac{1}{p}]}\mathcal E$. D'après la condition (PIM) de la définition~\ref{defiisochtouca} il existe $C$ tel que  pour tout $n\in \N$ l'image de $\tu{n}V_{1}$ par $ \vi _{\gN}\ta\vi _{\gN}\dotss  \tu{n-1}\vi_{\gN}$ est comprise entre $p^{C}
V_{1}$ et $p^{-C}V_{1}$. On pose alors 
$$V=\bigcup_{k\in \N}\bigcap _{n\in \N, n\geq k}
\vi_{\gN} \ta\vi_{\gN} \dotss  \tu{n-1}\vi_{\gN} (\tu{n}V_{1}).$$
Il est clair que $V$ est un $\O_{\mathcal E}$-réseau de $\gN\otimes _{\gS[\frac{1}{p}]}\mathcal E$, qui est préservé par $\vi _{\gN}
\otimes _{\gS[\frac{1}{E},\frac{1}{p}]}\Id_{\mathcal E}$. \cqfd

\begin{lem}\label{langtonOL}
a) 
Soit $\gN$ un $\gS[\frac{1}{p}]$-module libre de rang $r$, $V$ un $\O_{\mathcal E}$-module libre de rang $r$ et $\alpha:\gN\otimes _{\gS[\frac{1}{p}]}\mathcal E\to V\otimes_{\O_{\mathcal E}}\mathcal E$ un isomorphisme de $\mathcal E$-espaces vectoriels. Alors il existe un unique triplet $(\gM,\beta, \gamma)$, où 
$\gM$ est un $\gS$-module libre de rang $r$, 
$\beta:\gM\otimes_{\gS}\gS[\frac{1}{p}]\to \gN$ est un isomorphisme de $\gS[\frac{1}{p}]$-modules libres, 
$\gamma:\gM\otimes_{\gS} \O_{\mathcal E}\to V$ est un isomorphisme de $\O_{\mathcal E}$-modules libres,  et $\alpha\circ (\beta\otimes 1)=\gamma\otimes 1$. 

b) Soient $\gM$ et $\gM'$ des $\gS$-modules libres de type fini, et \begin{gather*}g\in \Hom_{\gS[\frac{1}{p}]} (\gM\otimes_{\gS}\gS[\frac{1}{p}], \gM'\otimes_{\gS}\gS[\frac{1}{p}]) \\ \text{ et \ \  }h\in \Hom_{\O_{\mathcal E}}(\gM\otimes_{\gS} \O_{\mathcal E},
\gM'\otimes_{\gS} \O_{\mathcal E})\end{gather*} tels que $g\otimes 1$ et $h\otimes 1$ coïncident dans $$\Hom_{\mathcal E}(\gM\otimes_{\gS} \mathcal E,
\gM'\otimes_{\gS} \mathcal E).$$ Alors il existe un unique morphisme $f: \gM\to \gM'$ tel que $g=f\otimes 1$ et $h=f\otimes 1$.
\end{lem}

\noindent {\bf Démonstration.}
Partant d'un $\gS$-module libre $\gM_{0}$ de rang $r$, muni d'un isomorphisme $\gM_{0}\otimes_{\gS}\gS[\frac{1}{p}]\simeq \gN$, il s'agit de le modifier au point générique du diviseur $p=0$. On peut se limiter à traiter le cas d'une modification élémentaire supérieure, c'est-à-dire que l'on suppose 
$$\gM_{0}\otimes _{\gS} \O_{\mathcal E}\subset V \subset p^{-1}\gM_{0}\otimes_{\gS} \O_{\mathcal E}$$  avec $V/\gM_{0}\otimes _{\gS} \O_{\mathcal E}$ un $k((u))$-espace vectoriel de dimension $1$  et on doit montrer que $p^{-1}\gM_{0}\cap V$ est un $\gS$-module libre de rang $r$. D'abord  $V/ \gM_{0}\otimes_{\gS} \O_{\mathcal E}$ est une droite  du $k((u))$-espace vectoriel 
$(p^{-1}\gM_{0}/\gM_{0})\otimes _{k[[u]]}k((u))$. Il est évident que l'intersection de $k[[u]]^{r}$ avec une droite de $k((u))^{r}$ est un sous-module libre  de rang $1$ facteur direct de $k[[u]]^{r}$ : il est engendré par un générateur de cette droite dont le minimum des valuations des coordonnées est nul. Donc il existe une base $e_{1},\dotss  ,e_{r}$ du $k[[u]]$-module libre $p^{-1}\gM_{0}/\gM_{0}$ telle que $k((u))e_{1}=V/ \gM_{0}\otimes_{\gS} \O_{\mathcal E}$. Soit $f_{1},..,f_{r}$ une base du $\gS$-module libre $\gM_{0}$ telle que $p^{-1}f_{i} \text{ mod } \gM_{0}=e_{i}$ pour $i=1,\dotss  ,r$. Alors $(p^{-1}f_{1},f_{2},\dotss  ,f_{r})$ est une base du $\gS$-module $p^{-1}\gM_{0}\cap V$ et celui-ci est donc libre de rang $r$. \cqfd

\begin{rem}
Le lemme~\ref{langtonOL} résulte aussi du fait que tout module de type fini réflexif sur un anneau local régulier de dimension $2$ est libre (lemme 6 de~\cite{serre-iwasawa}).  
\end{rem}

\section{Construction du foncteur de Dieudonné}\label{construc-Dieudonne}
\label{HPsanstorsion}

Soit $\A=K_{0}[[u]]$ le complété de $\gS[\frac{1}{p}]$ pour la topologie $u$-adique. Cet anneau fourre-tout contient plusieurs anneaux que nous introduirons ensuite et les inclusions et  intersections de ces anneaux auront lieu dans $\A$. 
 
On note $\O$ le sous-anneau de $\A$ formé des élements qui, considérés comme fonctions de la variable $u$, convergent sur le disque ouvert de rayon $1$.   
Pour décrire $\O$ de fa\c con plus explicite on a besoin d'introduire une norme ultramétrique $|.|$ sur $\QQ_p$, que l'on étend de fa\c con évidente à $K_{0}$ puis à $K$. Le choix de $|p|\in  ]0,1[$ est arbitraire mais on a bien sûr $|\pi_{K}|=|p|^{\frac{1}{e}}$. 

Pour $r\in ]0,1[$ et  $x=\sum_{n\in \N}x_{n}u^{n}\in \A=K_{0}[[u]]$ on pose 
$$|x|_{r}^{\naif}=\sup_{n\in \N}|x_{n}|r^{n}\in [-\infty, +\infty].$$
Pour $x,y\in \A\setminus\{0\}$ on a $|xy|_{r}^{\naif}=|x|_{r}^{\naif}|y|_{r}^{\naif}$. 

 On a alors 
 $$\O=\{x\in \A, \forall r\in ]0,1[, |x|_{r}^{\naif}<+\infty\}$$
 et pour $ r\in ]0,1[$ on note $|.|_{r}$ la restriction de $|.|^{\naif}_{r}$ à $\O$. 
 
 On munit $\O$ de la topologie  telle que les parties  $$\{x\in \O, |x|_{r}\leq \varepsilon\}$$ pour  $r\in ]0,1[$ et  $\varepsilon\in \R_{+}^{*}$  forment une base de voisinage de $0$.

 On a les inclusions $\gS[\frac{1}{p}]\subset \O\subset \A$.

On introduit $$\lambda=\prod_{n=0}^{\infty}
\vi^{n}\Big(\frac{E}{E(0)}   \Big)  \in \O.$$ 
Comme $\lambda=1$ modulo $u$, $\lambda$ est une unité dans $\A$, donc $\lambda$ n'est pas diviseur de $0$ dans  $\O$. 
On note d'ailleurs que $\lambda$ est inversible dans $\gS[[ \frac{u^{e}}{p}]]\subset \A$. 
On étend 
$|.|_{r}$ à $\O[\frac{1}{\lambda}]$ en posant $|x\lambda^{-n}|_{r}=|x|_{r}|\lambda|_{r}^{-n}$. On note que $|.|_{r}$ n'est pas la restriction de 
 $|.|^{\naif}_{r}$ à $\O[\frac{1}{\lambda}]\subset \A$, d'où l'adjectif ``naïf''. 
  Pour la suite on note que \begin{gather}\nonumber   |E(0)|_{r}= |p|_{r}=|p|=|\pi_{K}|^{e}, \ 
  |\vi^{m}(E)|_{r}=\max(|\pi_{K}|^{e},r^{p^{m}e}), \\ \label{norme-vi-n-E-r}
  |\lambda|_{r}=\prod_{m=0}^{\infty}\max\Big(1,\Big(\frac{r^{p^{m}}}{|\pi_{K}|}\Big)^{e}\Big)\end{gather} 
 et qu'en particulier 
 \begin{gather}\nonumber
  |\vi^{m}(E)|_{|\pi_{K}|^{p^{-n}}}=|p|^{p^{m-n}} \text{ pour }m<n, 
  |\vi^{m}(E)|_{|\pi_{K}|^{p^{-n}}}=|p| \text{ pour }m\geq n, \\
   |p|_{|\pi_{K}|^{p^{-n}}}=|p|
   \text{ \ \ 
 et \ \ }
  \label{vi-m-n-r}
  |\lambda|_{|\pi_{K}|^{p^{-n}}}=|p|^{(p^{-1}+...+p^{-n})-n}
 .\end{gather}

\begin{lem}\label{conv-suites-CCB} L'anneau $\O$ est complet vis-à-vis des normes $|.|_{r}$. En d'autres termes, 
si $(y_{n})_{n\in \N}$ est une suite d'éléments de $\O$ telle que pour tout $r\in ]0,1[$ la suite $|y_{n}-y_{n+1}|_{r}$ tend vers $0$, alors la suite $(y_{n})$ converge dans $\O$ vers une limite $y$ et pour tout $r\in ]0,1[$, $|y|_{r}=\lim_{n\to \infty} |y_{n}|_{r}$. \cqfd
\end{lem}

Pour tout $n\in \N$ on note  $\wS_{n}$ le complété de $\gS[\frac{1}{p}]$ par rapport à l'idéal engendré par $\vi^{n}(E)$. C'est un anneau local d'uniformisante $\vi^{n}(E)$. 
En effet $\vi^{n}(E)$ est un élément irréductible de $\gS$ car c'est un polynôme d'Eisenstein. 
On note que $\wS_{0}=\wS$ et que $\vi:\gS\to \gS$ s'étend en des morphismes continus $\vi:\wS_{n}\to \wS_{n+1}$ et 
 $\vi:\wS_{n}[\frac{1}{\vi^{n}(E)}]\to \wS_{n+1}[\frac{1}{\vi^{n+1}(E)}]$.

\begin{lem}\label{cont-epi-k}
Pour tout $n\in \N$, $\mathrm{Id}_{\gS}$ s'étend de fa\c con unique en  un morphisme 
$\O\to \wS_{n} $ tel que 
pour tout $k\in \N^{*}$ le morphisme 
$$\O\to \wS_{n} /\vi^{n}(E)^{k}\wS_{n}
=(\gS /\vi^{n}(E)^{k}\gS)[\frac{1}{u}]
$$ qui s'en déduit  est continu lorsque l'on munit  l'espace d'arrivée  de la topologie ``$u$-adique'' pour laquelle les voisinages de $0$ sont les 
$u^{m}(\gS /\vi^{n}(E)^{k}\gS)$. \cqfd \end{lem}

\begin{lem}\label{enlever-poles}
Soit $x\in \O[\frac{1}{\lambda}]$, tel que pour tout $n\in \N$, l'image de $x$ dans $\wS_n[\frac{1}{\vi^n(E)}]$
appartienne à $\wS_n$. 
Alors $x$ appartient à $\O$. \cqfd
\end{lem}

\begin{lem}\label{regul}
Soit $x\in \O$. On suppose qu'il existe $C\in \R$ tel que 
$|x|_{r}\leq C$ pour tout $r\in ]0,1[$. Alors $x$ appartient à $\gS[\frac{1}{p}]$. 
\end{lem}
\noindent {\bf Démonstration.} On écrit $x=\sum_{n\in \N}x_{n}u^{n}$ avec $x_{n}\in K_{0}$. Soit $n\in \N$. On a  $|x_{n}|\leq Cr^{-n}$ pour tout $r\in ]0,1[$, d'où $|x_{n}|\leq C$. \cqfd

Soit $(\gN,\vi_{\gN})$ un pseudo-iso-$\vi / \gS$-module. 
Nous allons construire son image par $\mathbb D_{\mathrm{iso}}$, qui est un  $\vi$-module de Hodge-Pink $D=(D,\vi_D,V_{D})$,   c'est-à-dire un $K_0$-espace
vectoriel $D$ muni d'un 
isomorphisme de $K_0$-espaces vectoriels $\vi_D:{\ta D}\to D$, et
d'un $\wS$-réseau $V_{D}$ 
dans ${\ta D}\otimes _{K_0}\wS[\frac{1}{E}]$. 

On remarque d'abord que $$\gS[\frac{1}{p}]/ u\gS[\frac{1}{p}]
=\gS[\frac{1}{p},\frac{1}{E}]/u\gS[\frac{1}{p},\frac{1}{E}]=K_0. $$
Soit $(D,\vi_{D})$ le $\vi$-module défini par  
\begin{gather}\label{debut-def-Diso}
D=\gN \otimes _{\gS[\frac{1}{p}]} K_0=
\gN/u\gN\text{\ \  et\ }\vi_D=\vi_{\gN} \ \text{mod}\ 
u. \end{gather}

\begin{lem}\label{rigidite2}
Il existe un unique morphisme $\xi$ de $\O[\frac{1}{\lambda}]$-modules
 congru à $1$ modulo $u$ de $D\otimes _{K_0} \O[\frac{1}{\lambda}]$ vers  
$\gN\otimes _{\gS[\frac{1}{p}]}\O[\frac{1}{\lambda}]$ 
 vérifiant 
\begin{gather}
\label{eq-fonct}
\xi(\vi_D \otimes 1)=(\vi_\gN\otimes 1)\ta \xi.
\end{gather} 
De plus c'est un
 isomorphisme de $\O[\frac{1}{\lambda}]$-modules et il dépend fonctoriellement de $ (\gN, \vi_{\gN})$. 
\end{lem}

\noindent {\bf Démonstration.}
L'unicité de $\xi$ est évidente car \eqref{eq-fonct} détermine $\xi$ modulo $u^{p}, u^{p^{2}}$,... et $\O[\frac{1}{\lambda}]\subset \A$ est une inclusion. 

Montrons l'existence. Soit $s_{0}\in \mathrm{Hom}_{\gS[\frac{1}{p}]}(D\otimes _{K_0}\gS[\frac{1}{p}],\gN)$ égal à $1$ modulo $u$. 
On choisit aussi un $W$-réseau $\Delta$ de $D$ et un $\gS$-module libre $\gM$ avec un  isomorphisme $\gM[\frac{1}{p}]=\gN$.

 On pose alors \begin{gather*}\xi=\lim _{n\mapsto
+\infty}s_n \text{\ \  dans \ } \mathrm{Hom}_{\gS[\frac{1}{p}]}(D\otimes _{K_0}\gS[\frac{1}{p}],\gN) 
\otimes _{\gS[\frac{1}{p}]}\O[\frac{1}{\lambda}]\\ \text{\ \ où 
\ }s_n=\vi_{\gN} \dotss   \tu{n-1} \vi_{\gN}
\tu{n} s_{0}\tu{n-1}\vi_{D}^{-1}\dotss   \vi_{D}^{-1} \\ \in 
\mathrm{Hom}_{\gS[\frac{1}{p}]}(D\otimes _{K_0}\gS[\frac{1}{p}],\gN)[\frac{1}{E},\dotss  ,\frac{1}{\vi^{n-1}(E)}].
\end{gather*} 
Cette limite existe pour la raison suivante. 
Soit $C\in \N^{*}$ tel  que \begin{gather}\label{borne-viN-viD}\vi_{\gN}\in p^{-C}E^{-C}
\mathrm{Hom}_{\gS}(\ta\gM, \gM), \ \   \vi_{D}^{-1}\in p^{-C}\mathrm{Hom}_{W}(\Delta ,\ta  \Delta)\\ \label{borne-r}\text{ et \ \ \ } 
s_{0}\in p^{-C}\mathrm{Hom}_{\gS}(\Delta\otimes _{W}\gS,\gM).\end{gather} Alors 
$$(s_{n+1}-s_{n})\in p^{-C(2n+3)}E^{-C}\dotss  \vi^{n}(E)^{-C}u^{p^{n}}\mathrm{Hom}_{\gS}(\Delta \otimes_{W}\gS,\gM)$$   
$$\text{ donc }\lambda^{C}(s_{n+1}-s_{n})\in p^{-C(3n+4)}u^{p^{n}}\big(\tu{n+1} \lambda\big)^{C}\mathrm{Hom}_{\gS}(\Delta \otimes_{W}\gS,\gM)$$ et donc 
$$\lambda^{C}(s_{n+1}-s_{n})\in p^{-C(3n+4)}u^{p^{n}}\mathrm{Hom}_{\gS}(\Delta \otimes_{W}\gS,\gM)
\otimes_{\gS}\gS[[\frac{u^{ep^{n+1}}}{p}]].$$ 
Le lemme~\ref{conv-suites-CCB} montre que la suite $\lambda^{C}s_{n}$
 converge dans $\mathrm{Hom}_{\gS}(\Delta \otimes_{W}\gS,\gM)
\otimes_{\gS}\O$. Donc $\xi$ appartient à  $\mathrm{Hom}_{\gS[\frac{1}{p}]}(D\otimes_{K_0}\gS[\frac{1}{p}],\gN)
\otimes _{\gS[\frac{1}{p}]}\lambda^{-C}\O$. 
Ceci termine la démonstration de l'existence de $\xi$. 

De la même fa\c con, on choisit  $s'_{0}\in \mathrm{Hom}_{\gS[\frac{1}{p}]}(\gN,D\otimes _{K_0}\gS[\frac{1}{p}])$ égal à $1$ modulo $u$,  on pose $\xi'=\lim _{n\mapsto
+\infty}\vi_{D} \dotss   \tu{n-1} \vi_{D}
\tu{n} s'_{0}\tu{n-1}\vi_{\gN}^{-1}\dotss   \vi_{\gN}^{-1} $ et on montre 
que 
$\xi'$ appartient à  $\mathrm{Hom}_{\gS[\frac{1}{p}]}(\gN,D\otimes_{K_0}\gS[\frac{1}{p}])
\otimes _{\gS[\frac{1}{p}]}\O[\frac{1}{\lambda}]$. Une adaptation de l'argument d'unicité de $\xi$ montre que $\xi'\xi=1$ et $\xi\xi'=1$, donc $\xi$ est un isomorphisme. 
\cqfd

\begin{rem} C'est le lemme 1.2.6 de~\cite{kisin}. La démonstration précédente est exactement la même que celle donnée dans~\cite{kisin}. Elle est aussi parallèle à celle du 
 lemme 7.4    
  de~\cite{fontaine29} (la transcription est détaillée dans le dernier paragraphe).

\end{rem}

\begin{rem} 
Dans le cadre des $(\vi,\Gamma)$-modules apparaissent  d'autres relèvements de $\vi$ de $W$ à $\gS$, qui sont donnés par  \begin{gather}\label{rel-frob}u\mapsto u^{p}+h_{p-1}u^{p-1}+...+h_{1}u\text{ \ avec \ }h_{i}\in pW, 
\end{gather} au lieu de $u\mapsto u^{p}$. On renvoie à~\cite{kisin-ren} pour une situation proche de la nôtre. 
Pour un tel relèvement de $\vi$ la suite $s_{n}$ de la démonstration précédente   ne converge pas nécessairement (sauf évidemment en rang $1$).   D'ailleurs la démonstration du lemme 2.2.2 de~\cite{kisin-ren} est de nature différente, car elle utilise la connexion induite par l'action de $\Gamma$.  En fait
l'analogue brutal du lemme~\ref{rigidite2} est faux dans le cadre de~\cite{kisin-ren}, comme le montre l'exemple suivant. 
 On prend $E(u)=\frac{(1+u)^{p}-1}{u}$ et le relèvement de $\vi$ 
de $W$ à $\gS$ donné par  $u\mapsto (1+u)^{p}-1$. On prend $\gN=\gS[\frac{1}{p}]^{2}$ et $\vi_{\gN}=
\begin{pmatrix}
1  & x \\ 0 & E
\end{pmatrix}$ avec $x\in u+u^{2}\gS$. 
Par fonctorialité $\xi$ devrait être de la forme 
$\begin{pmatrix}
1  & y \\ 0 & \lambda
\end{pmatrix}$ avec $y\in u\O[\frac{1}{\lambda}]$ et 
d'après 
\eqref{eq-fonct},  $y$ devrait vérifier \begin{gather}\label{eq-fonct-y}y-p^{-1}\vi(y)=p^{-1}\ta \lambda x.\end{gather} On écrit $y=\sum_{n=1}^{\infty }y_{n}u^{n}$. Alors   $\vi(y)=
\sum_{n=1}^{\infty }y_{n}\Big((1+u)^{p}-1 \Big)^{n}=py_{1}u \text{ mod } u^{2}$. Donc $y-p^{-1}\vi(y)$ est divisible par $u^{2}$ alors que $p^{-1}\ta \lambda x$ est congru à $p^{-1}u$ modulo $u^{2}$, donc \eqref{eq-fonct-y} ne peut être satisfaite. En fait en revenant au cadre g\'en\'eral du lemme précédent, mais avec un
rel\`evement de $\vi$ comme dans \eqref{rel-frob},
l'existence de $\xi$ vérifiant   \eqref{eq-fonct} \'equivaut \`a celle 
d'une solution modulo une puissance finie de $u$. En effet on peut montrer que si $C$ est comme dans \eqref{borne-viN-viD}, et si $s_{0}$ vérifie 
 $s_{0}(\vi_D\otimes 1)=(\vi_\gN \otimes 1)\ta s_{0}$ modulo $u^{3C+1}$, alors la suite $s_{n}$ converge dans $\lambda^{-C}\mathrm{Hom}_{\gS[\frac{1}{p}]}(D\otimes _{K_0}\gS[\frac{1}{p}],\gN)
\otimes _{\gS[\frac{1}{p}]}\O$ et sa limite $\xi$ vérifie  \eqref{eq-fonct}. 
 \end{rem}

Nous allons maintenant définir la structure de Hodge-Pink $V_{D}$ sur $D$ et terminer ainsi la construction de $\mathbb D_{\mathrm{iso}}$. 

\begin{defi}\label{defi-HP-gen} Le foncteur $\mathbb D_{\mathrm{iso}}$ envoie  $(\gN,\vi_{\gN})$ sur $(D,\vi_D,V_{D})$, où 
$D$ et $\vi_{D}$ sont définis par \eqref{debut-def-Diso}, et où la structure de Hodge-Pink $V_{D}$  est le $\wS$-réseau   de $\ta D\otimes _{K_0} \wS[\frac{1}{E}] $ défini par 
\begin{gather} V_{D}=(\vi_{D}\otimes 1)^{-1}(\xi^{-1}(\gN\otimes _{\gS[\frac{1}{p}]}\wS))
  \label{defi-HP-vd}=\ta \xi^{-1}((\vi_\gN\otimes 1)^{-1}(\gN\otimes _{\gS[\frac{1}{p}]}\wS)).\end{gather}
\end{defi}

On va montrer que le foncteur $\mathbb D_{\mathrm{iso}}$ est pleinement fidèle  en construisant un foncteur quasi-inverse sur son image. 
Autrement dit on va montrer que  $V_{D}$ détermine de manière unique 
$\xi^{-1}(\gN)\subset D\otimes_{K_0}\O[\frac{1}{\lambda}]$, et par conséquent détermine le couple $(\gN,\vi_{\gN})$ à un unique isomorphisme près.

Soit  
$\wt \gN$ l'ensemble des $x\in D\otimes _{K_0} \O[\frac{1}{\lambda}]$, tels que 
\begin{itemize}
\item
pour tout $n\in \N$, 
l'image de $x$ dans $D\otimes _{K_0} \wS_{n}[\frac{1}{\vi^{n}(E)}]$ appartienne à $\vi_{D}\dotss  \tu{n}\vi_{D}\tu{n}V_{D}
$, 
\item il existe une constante $C\in \R$ telle que pour tout $n\in \N$, 
\begin{gather}\label{derniere-cond-tildeN}|(\vi_{D}\dotss  \tu{n-1}\vi_{D})^{-1}x|_{|\pi_{K}|^{p^{-n}}}\leq C.\end{gather}
\end{itemize}

On précise que   $\tu{n}V_{D}$ est un $\wS_{n}$-réseau de 
$\tu{n+1} D\otimes _{K_0}\wS_{n}[\frac{1}{\vi^{n}(E)}]$ et donc 
$\vi_{D}\dotss  \tu{n}\vi_{D}\tu{n}V_{D}$ est un $\wS_{n}$-réseau de 
$ D\otimes _{K_0}\wS_{n}[\frac{1}{\vi^{n}(E)}]$. 

Pour donner un sens à \eqref{derniere-cond-tildeN} il faut choisir un $W$-réseau $\Delta$  de $D$ : dans une base 
de $\tu{n}\Delta$  sur $ W$  on écrit 
$(\vi_{D}\dotss  \tu{n-1}\vi_{D})^{-1}x=\begin{pmatrix}
y_{1}\\ \vdots \\y_{r}\end{pmatrix}$ avec $y_{i}\in \O[\frac{1}{\lambda}]$ et on  pose $|(\vi_{D}\dotss  \tu{n-1}\vi_{D})^{-1}x|_{|\pi_{K}|^{p^{-n}}}=\max(|y_{1}|_{|\pi_{K}|^{p^{-n}}},...,|y_{r}|_{|\pi_{K}|^{p^{-n}}})$. 
 Comme $|pb|_{|\pi_{K}|^{p^{-n}}}=|p|.|b|_{|\pi_{K}|^{p^{-n}}}$  pour tout $b\in \O[\frac{1}{\lambda}]$, si on change $\Delta$ la dernière condition reste vraie pour une autre constante $C$, donc la définition de $\wt \gN$ ne dépend pas du choix de $\Delta$. On vérifie facilement  que $ (\vi_{D}\otimes 1)(\ta {\wt \gN}[\frac{1}{E}])= \wt \gN[\frac{1}{E}]$.

\begin{prop}~\label{VdetermineN}
On a  $\xi^{-1}(\gN)= \wt \gN$. 
\end{prop}
\noindent {\bf Démonstration.} On suppose $(\gN,\vi_{\gN})$ d'amplitude  $\subset [s,t]$. On fixe  un $\gS$-module libre $\gM$ muni d'un isomorphisme $\gN=\gM[\frac{1}{p}]$ et une constante $C_{0}$ telle 
que la condition (PIM) de la définition~\ref{defiisochtouca} soit satisfaite (avec $C_{0}$ au lieu de $C$). On fixe aussi un réseau $\Delta$ de $D$. Grâce à $\Delta $ et $\gM$ toutes les normes qui vont suivre seront bien définies. On choisit  des bases de $\gM$ et  de $\Delta$. Ce sont donc aussi des bases de $\gN$ et de $D$. 
On va 
 appliquer les lemmes~\ref{enlever-poles} et~\ref{regul}
aux cordonnées d'un élément de $\xi(\wt \gN)$. 
En effet $\xi(\wt \gN)$ est formé des $y\in  
\gN\otimes _{\gS[\frac{1}{p}]}\O[\frac{1}{\lambda}]$ tels que
\begin{itemize}
\item (C1) pour tout $n\in \N$, 
l'image de $y$ dans $\gN\otimes _{\gS[\frac{1}{p}]}\wS_{n}[\frac{1}{\vi^{n}(E)}]$ appartienne à 
$\xi(\vi_{D}\dotss  \tu{n}\vi_{D}\tu{n}V_{D})$, 
\item (C2) 
il existe $C\in \R$ telle que pour tout $n\in \N$, 
$$|(\vi_{D}\dotss  \tu{n-1}\vi_{D})^{-1}\xi^{-1}y|_{|\pi_{K}|^{p^{-n}}}\leq C.$$
\end{itemize}

On a  
\begin{gather*}\xi(\vi_{D}\dotss  \tu{n}\vi_{D}\tu{n}V_{D}) =
\xi(\vi_{D}\dotss  \tu{n-1}\vi_{D}\tu{n}\xi^{-1}(\tu{n}\gN\otimes _{\gS[\frac{1}{p}]}\wS_{n}))\\ =\vi_{\gN}\dotss  \tu{n-1}\vi_{\gN}(\tu{n}\gN\otimes _{\gS[\frac{1}{p}]}\wS_{n}) =\gN\otimes _{\gS[\frac{1}{p}]}\wS_{n}\end{gather*}
car  $\vi_{\gN},\dotss  ,\tu{n-1}\vi_{\gN}$ sont inversibles  à coefficients dans $\wS_{n}$. D'après le lemme~\ref{enlever-poles} la  condition (C1) équivaut donc à $y\in \gN\otimes _{\gS[\frac{1}{p}]}\O$.

Comme $\tu{n-1}\vi_{D}^{-1}\dotss  \vi_{D}^{-1}\xi^{-1}y=\tu{n}\xi^{-1}\tu{n-1}\vi_{\gN}^{-1}\dotss  \vi_{\gN}^{-1} y$ la  condition (C2) équivaut à l'existence de   $C\in \R$ tel que pour tout $n\in \N$, 
$$|\tu{n}\xi^{-1} \tu{n-1}\vi_{\gN}^{-1}\dotss  \vi_{\gN}^{-1}y|_{|\pi_{K}|^{p^{-n}}}\leq C.$$

Il existe $C_{1}\in \N$ tel que les coefficients de $\xi$ et $\xi^{-1}$ 
appartiennent à 
$p^{-C_{1}}\gS[[\frac{u^{e}}{p}]]$, d'où 
$|\tu{n}\xi|_{|\pi_{K}|^{p^{-n}}}\leq |p|^{-C_{1}}$ et 
$|\tu{n}\xi^{-1}|_{|\pi_{K}|^{p^{-n}}}\leq |p|^{-C_{1}}$
 $$\text{ et donc \ }|p|^{C_{1}}
\leq  \frac{|\tu{n}\xi^{-1} \tu{n-1}\vi_{\gN}^{-1}\dotss  \vi_{\gN}^{-1}y|_{|\pi_{K}|^{p^{-n}}}}
 {| \tu{n-1}\vi_{\gN}^{-1}\dotss  \vi_{\gN}^{-1}y|_{|\pi_{K}|^{p^{-n}}}}\leq |p|^{-C_{1}}.
 $$
 
Enfin, grâce à  \eqref{vi-m-n-r}  la condition (PIM)   donne l'encadrement suivant : 
\begin{gather*}|p|^{-s(p^{-1}+\dotss  +p^{-n}) +C_{0}}.|y|_{|\pi_{K}|^{p^{-n}}}
\leq 
|
 \tu{n-1}\vi_{\gN}^{-1}\dotss  \vi_{\gN}^{-1}y|_{|\pi_{K}|^{p^{-n}}}
\\ \leq   |p|^{-t(p^{-1}+\dotss  +p^{-n}) -C_{0}} |y|_{|\pi_{K}|^{p^{-n}}}.\end{gather*}
On voit alors que la  condition (C2) équivaut à l'existence de  
 $C\in \R$ tel que pour tout $n\in \N$, 
$|y|_{|\pi_{K}|^{p^{-n}}}\leq C$.

La  condition (C1) équivaut à $y\in \gN\otimes _{\gS[\frac{1}{p}]}\O$. Sous cette hypothèse,  grâce au 
  lemme~\ref{regul},  la  condition (C2) équivaut à $y\in \gN$, ce qui termine la démonstration de la proposition~\ref{VdetermineN}.   \cqfd

\begin{prop}\label{parag2-plein-fidel} Le foncteur $\mathbb D_{\mathrm{iso}}$ 
est pleinement fidèle. 
\end{prop}
\noindent{\bf Démonstration}
Cela résulte immédiatement de la  proposition~\ref{VdetermineN}. \cqfd

On note que la proposition~\ref{VdetermineN} construit explicitement un quasi-inverse du foncteur $\mathbb D_{\mathrm{iso}}$ sur son image essentielle. On appellera admissibles les objets de son image essentielle. 

La proposition suivante est l'implication  ``admissible implique faiblement admissible''. 

\begin{prop}\label{prop-HP-fa}
L'image du  foncteur $\mathbb D_{\mathrm{iso}}$ est contenue dans la sous-catégorie pleine des 
$\vi$-modules de Hodge-Pink faiblement admissibles.
\end{prop}

Soit $(\gN,\vi_{\gN})$ un   pseudo-iso-$\vi / \gS$-module et 
 $(D,\vi_{D},V_{D})$ 
 le $\vi$-module de Hodge-Pink qui est son image par $\mathbb D_{\mathrm{iso}}$.

\begin{lem}\label{propR} Pour tout $\gS$-réseau $\gM$ dans  $\gN$ et pour tout $W$-réseau $\Delta$ dans $D$,
 il existe une constante $C$ telle que pour tout entier $n\in
\N$, on a 
\begin{gather}\label{ineg-e-1}|(\vi_{D}\ta \vi_{D} \dotss  \tu{n-1}\vi_{D})^{-1}\xi^{-1}|_{|\pi_{K}|^{p^{-n}}}\leq C \\ \label{ineg-e-2}
\text{et \ \ \ }
|\xi(\vi_{D}\ta \vi_{D} \dotss  \tu{n-1}\vi_{D})|_{|\pi_{K}|^{p^{-n}}}\leq C .\end{gather}
\end{lem}

\noindent {\bf Démonstration.}
On suppose $(\gN,\vi_{\gN})$ d'amplitude $\subset [s,t]$. 
D'après la condition (PIM) il  existe $C\in \N$ tel que pour tout $n\in \N$, 
  \begin{gather*}\vi_{\gN}\ta \vi_{\gN} \dotss  \tu{n-1}\vi_{\gN}{}\in p^{-C}E^{s}\dotss  \vi^{n-1}(E)^{s}\mathrm{Hom}_{\gS}(\tu{n}\gM,\gM)
  \\
  \text{ et \ \ } 
 (\vi_{\gN}\ta \vi_{\gN} \dotss  \tu{n-1}\vi_{\gN}{})^{-1}\in p^{-C}E^{-t}\dotss  \vi^{n-1}(E)^{-t}\mathrm{Hom}_{\gS}(\gM,\tu{n} \gM)
.\end{gather*}
  D'après \eqref{vi-m-n-r},  pour tout $n\in \N$, 
$$|\vi_{\gN}\ta \vi_{\gN}\dotss  \tu{n-1}\vi_{\gN}|_{|\pi_{K}|^{p^{-n}}}\leq 
|p^{-C}E^{s}\dotss  \vi^{n-1}(E)^{s}|_{|\pi_{K}|^{p^{-n}}}=|p|^{
-C+s(p^{-1}+\dotss  +p^{-n})}
$$ et 
\begin{gather*}|(\vi_{\gN}\ta \vi_{\gN} \dotss  \tu{n-1}\vi_{\gN}{})^{-1}|_{|\pi_{K}|^{p^{-n}}}\leq 
|p^{-C}E^{-t}\dotss  \vi^{n-1}(E)^{-t}|_{|\pi_{K}|^{p^{-n}}}
\\
=|p|^{
-C-t(p^{-1}+\dotss  +p^{-n})}.\end{gather*} 
Or  \begin{gather*}(\vi_{D}\ta \vi_{D} \dotss  \tu{n-1}\vi_{D})^{-1}\xi^{-1}=\tu{n}\xi^{-1}(\vi_{\gN}\ta \vi_{\gN}\dotss  \tu{n-1}\vi_{\gN})^{-1},
\\
\xi(\vi_{D}\ta \vi_{D} \dotss  \tu{n-1}\vi_{D})=(\vi_{\gN}\ta \vi_{\gN}\dotss  \tu{n-1}\vi_{\gN})\tu{n}\xi 
\\
\text{et  d'autre part \ \ \ }|\tu{n} \xi|_{|\pi_{K}|^{p^{-n}}}=|\xi|_{|\pi_{K}|}\text{\ \ 
et \ \ }
|\tu{n} \xi^{-1}|_{|\pi_{K}|^{p^{-n}}}=|\xi^{-1}|_{|\pi_{K}|}.\end{gather*} Le lemme~\ref{propR} est démontré. \cqfd

\noindent{\bf Démonstration de 
 la proposition~\ref{prop-HP-fa}. } 
Pour tout entier $l$, $\Lambda^l(\gN)$ est un pseudo-iso-$\vi / \gS$-module et le $\vi$-module de Hodge-Pink associé est $\Lambda^l(D)$. En 
prenant $l$ égal au rang de $\gN$, on voit qu'il existe un entier $t\in \Z$ et $a\in (\gS[\frac{1}{p}])^{\times}$ tels que $\vi_{\Lambda^l(\gN)}=E^{t}a$ dans une base de $\Lambda^l(\gN)$. 
On a bien sûr 
$V_{\Lambda^l(D)}=E^{-t}U_{\Lambda^l(D)}$ donc 
$t_{H}(D)=t$. 
Par la condition (PIM), il existe $C$ tel que pour tout $n\in \N$, $a\vi( a) \dotss  \vi^{n-1}( a)$ et son inverse appartiennent à $p^{-C}\gS$. On en déduit $a\in \gS^{\times}$ et donc $t_{N}(D)=t$. Donc on a bien $t_{H}(D)=t_{N}(D)$.

Supposons par l'absurde que
$D$ n'est pas faiblement admissible. Il existe alors un
sous-$\vi$-module $D'$ de $D$ tel que $t_H(D')>t_N(D')$. 
En rempla\c cant $\gN$ par $\Lambda^l(\gN)$, où $l$ est la dimension de
$D'$, on se ramène à la situation où $D'$ est de
dimension $1$. 
Pour toute base $(e_1,\dotss  ,e_r)$ de $\gN$ sur $\gS[\frac{1}{p}]$, on note
 $(f_1,\dotss  ,f_r)$
la base de $D$ sur $K_0$ obtenue par réduction modulo
$u $. On a  $K_0\subset \gS[\frac{1}{p}]$ et $GL_{r}(K_0)\subset GL_{r}(\gS[\frac{1}{p}])$, donc on 
 peut  choisir la base $(e_1,\dotss  ,e_r)$ telle que
$D'=K_0f_1$. Notons  $j=t_N(D')$ de sorte que 
$\vi_D(f_1)=ap^jf_1$ avec $a\in W^{\times}$. On suppose  par l'absurde que 
$t_H(D')>t_N(D')$, c'est-à-dire 
$$E^{-(j+1)}\ta D'\otimes _{K_0}\wS\subset (\vi_{D}\otimes 1)^{-1}(\xi^{-1}(\gN\otimes _{\gS}\wS))$$
ou encore, de fa\c con équivalente, 
$f_1\in E^{j+1}\xi^{-1}(\gN\otimes_{\gS}\wS)$.  
Comme  $\xi=(\vi_{\gN}\otimes 1) \ta \xi(\vi_{D}\otimes 1)^{-1}$, comme $\vi_{\gN}$ induit un isomorphisme de 
$\ta \gN\otimes_{\gS}\wS_{m}$ dans $\gN\otimes_{\gS}\wS_{m}$ pour tout $m>0$,  et comme  $D'$ est stable par $\vi_{D}$, on voit que 
\begin{gather}\label{R-1f1-appartient-a}\xi(f_1)\text{ \   appartient à \ }
\vi^{m}(E)^{j+1}(\gN\otimes_{\gS[\frac{1}{p}]}\wS_{m})\text{ \   pour tout  \ } m\in 
\N.
\end{gather} On note $\xi(f_1)_1$ le  premier coefficient  de $\xi(f_1)$ dans la base $(e_1,\dotss  ,e_r)$. On a alors 
\begin{itemize}
\item $\lambda^{-(j+1)}\xi(f_1)_1$ est congru à $1$ modulo $u $, 
\item $\lambda^{-(j+1)}\xi(f_1)_1\in \O$ d'après \eqref{R-1f1-appartient-a} et grâce au  lemme~\ref{enlever-poles}, 
 \item pour tout entier $n\in
\N^*$, 
$|\xi(f_1)_1|_{|\pi_{K}|^{p^{-n}}}\leq C|p|^{-jn} $, grâce à \eqref{ineg-e-2} et au fait que $$\vi_{D}\ta \vi_{D}... \tu{n-1}\vi_{D}\tu{n}f_{1} =(a\ta a ... \tu{n-1}a) p^{nj}f_{1}$$
et donc 
$|\lambda^{-(j+1)}\xi(f_1)_1|_{|\pi_{K}|^{p^{-n}}}\leq C'|p|^{n} $ grâce à \eqref{vi-m-n-r}, pour une autre constante $C'$. 
\end{itemize}

Or ces propriétés sont contradictoires 
car $\O$ ne contient aucun élément ayant de telles propriétés, car tout $x\in \O$ congru à $1$ modulo $u $ vérifie $|x|_{r}\geq 1$ pour tout $r\in ]0,1[$. La proposition~\ref{prop-HP-fa} est démontrée. \cqfd

Nous allons énoncer une 
 proposition qui  permettra dans le paragraphe suivant une réduction de la preuve de ``faiblement admissible implique admissible'' aux objets irréductibles. 

Une suite exacte courte  dans la catégorie des $\vi$-modules de Hodge-Pink est une suite exacte $0\to D'\to D\to D''\to 0$ dans la catégorie des $\vi$-modules  telle que $V_{D ''}$ soit l'image de $V_{D}$ dans $D''\otimes_{K_0}\wS[\frac{1}{E}]$ et $V_{D'}$ l'intersection de $V_{D}$ avec $D'\otimes_{K_0}\wS[\frac{1}{E}]$. Le foncteur $\mathbb D_{\mathrm{iso}}$ transforme les suites exactes en suites exactes. 

\begin{prop}\label{prop3generalisee} Une extension de deux $\vi$-modules de Hodge-Pink admissibles est admissible. 
\end{prop}

\noindent {\bf Début de la démonstration de la proposition~\ref{prop3generalisee}.} 
Soient  $(\gN',\vi_{\gN'})$ et $(\gN'',\vi_{\gN''})$ des  pseudo-iso-$\vi / \gS$-modules, et $(D',\vi_D',V_{D'})$ et $ (D '',\vi_{D ''},V_{D ''})$  leurs images par $\mathbb D_{\mathrm{iso}}$. Soit $D=(D,\vi_{D},V_{D})$ une extension de $ (D '',\vi_{D ''},V_{D ''})$ par  $(D',\vi_D',V_{D'})$. On doit montrer que $D$  provient d'une extension $(\gN,\vi_{\gN})$ de 
$(\gN'',\vi_{\gN''})$ par $(\gN',\vi_{\gN'})$.

On se ramène  à montrer la proposition~\ref{prop3generalisee} dans le cas 
où $(D'',\vi_{D''},V_{D''})$ est trivial (on indiquera  dans la remarque~\ref{rem-tens-dual}  comment le lecteur qui le souhaiterait peut traiter directement le cas général, au prix de notations plus compliquées). En effet la catégorie des pseudo-iso-$\vi / \gS$-modules et celle des $\vi$-modules de Hodge-Pink  possèdent des opérations produit tensoriel et dual. Par exemple  le dual de $(\gN'',\vi_{\gN''})$ est $({\gN''}^{*},{}^{t} \vi_{\gN''}^{-1})$
 où ${}^{t} \vi_{\gN''}: {\gN''}^{*} \to \ta {\gN''}^{*}$ est le transposé de $\vi_{\gN''}$. 
Le produit tensoriel de $(\gN',\vi_{\gN'})$ et du dual de $(\gN'',\vi_{\gN''})$ est $(\gN'\otimes {\gN''}^{*},\vi_{\gN'}\otimes {}^{t} \vi_{\gN''}^{-1})$ et il est équivalent de se donner une extension de $(\gS[\frac{1}{p}],\Id)$ par cet objet ou une extension de $(\gN',\vi_{\gN'})$ par $(\gN'',\vi_{\gN''})$. On est donc ramené à montrer la proposition~\ref{prop3generalisee} dans le cas 
où $(D'',\vi_{D''},V_{D''})$ est trivial, c'est-à-dire $$D''=K_0, \ \ \vi_{D''}=\mathrm{Id}\text{ \ \  et \ \ } V_{D''}=\wS. $$ 

On fixe un isomorphisme de $K_0$-espaces vectoriels 
$D=D'\oplus K_0$. On a alors $\vi_{D}=
\begin{pmatrix} \vi_{D'} & Y
\\ 0 & 1\end{pmatrix}$ pour un certain $Y\in D'$ et $V_{D}$  est un $\wS$-réseau de $\ta D\otimes_{K_0}\wS[\frac{1}{E}]$ qui est 
une extension de $\wS$ par $V_{D'}$. 
Il existe  $T\in \ta D'\otimes _{K_0} \wS[\frac{1}{E}]$ tel que $$(T, 1)\in \ta D'\otimes _{K_0} \wS[\frac{1}{E}] \oplus \wS[\frac{1}{E}] =
\ta D\otimes_{K_0}\wS[\frac{1}{E}]$$ appartienne à $V_{D}$ et $T$ est unique modulo $V_{D'}$. En d'autres termes $V_{D}$ est déterminé par un élément $T\in \ta D'\otimes _{K_0} \wS[\frac{1}{E}]/V_{D'}$.

On pose   $\gN=\gN'\oplus \gS[\frac{1}{p}]$ comme  $\gS[\frac{1}{p}]$-module.  Soit $X_{0}\in \gN'[\frac{1}{E}]$ égal à $Y$  modulo $u$. On pose 
$\vi_{\gN}=\begin{pmatrix} \vi_{\gN'} & X_{0}+X
\\ 0 & 1\end{pmatrix}$, avec $X \in u \gN'[\frac{1}{E}]$ arbitraire (nous choisirons ensuite $X$ tel que la structure de Hodge-Pink associée à $\vi_{\gN}$ soit $V_{D}$).

\begin{lem}\label{interm-ext}
Avec les notations précédentes $(\gN,\vi_{\gN})$ est un pseudo-iso-$\vi / \gS$-module. 
\end{lem}
\noindent {\bf Démonstration.}
La  condition (PIM) de la définition~\ref{defiisochtouca} est vérifiée car dans le produit 
$\vi_{\gN} \ta \vi_{\gN} \dotss  \tu{n-1}\vi_{\gN}$ le terme non diagonal apparaît $0$ ou $1$ fois, donc n'ajoute pas de dénominateur  en $p$ d'ordre plus grand qu'une constante indépendante de $n$, et il en va de même pour  $   \tu{n-1}\vi_{\gN}^{-1}\dotss  \ta \vi_{\gN}^{-1}\vi_{\gN}^{-1}$. \cqfd

\noindent {\bf Fin de la démonstration de la proposition~\ref{prop3generalisee}.} 
On va montrer qu'on peut  choisir $X$ de sorte que si $$\xi:
D\otimes _{K_0} \O[\frac{1}{\lambda}] \to \gN\otimes _{\gS[\frac{1}{p}]}\O[\frac{1}{\lambda}]
$$ est l'unique morphisme congru à $1$ modulo $u$  vérifiant  \begin{gather*}\xi(\vi_{D}\otimes 1)=(\vi_{\gN} \otimes 1)\ta \xi \text{\ \ \ \ 
  alors} \\ \ta \xi^{-1} (\vi_{\gN}^{-1}(\gN\otimes _{\gS[\frac{1}{p}]}\wS))= \vi_{D}^{-1}(\xi^{-1} (\gN\otimes _{\gS[\frac{1}{p}]}\wS)) 
\subset \ta D\otimes _{K_0}\wS[\frac{1}{E}]\end{gather*} 
 soit l'extension $V_{D}$  de $\wS$ par $V_{D'}$ qui est donnée au départ (et qui est une extension absolument arbitraire). 
On a $\xi=\begin{pmatrix} {\xi'} & Z
\\ 0 & 1\end{pmatrix}$ pour un certain $$Z\in \gN'\otimes _{\gS[\frac{1}{p}]} \O[\frac{1}{\lambda}]
\subset \gN'\otimes _{\gS[\frac{1}{p}]}\wS[\frac{1}{E}]$$ dépendant de $X$. 
Comme $\xi^{-1}=\begin{pmatrix} {\xi'}^{-1} & -{\xi'}^{-1} (Z)
\\ 0 & 1\end{pmatrix}$, 
 il s'agit de montrer que l'on peut trouver $X$ tel que $-\vi_{D'}^{-1}({\xi'}^{-1} (Z))\in \ta D'\otimes _{K_0}\wS[\frac{1}{E}]$ soit égal à $T$ modulo 
$V_{D'}$.  Cela est équivalent à la condition suivante : \begin{gather*}Z\in \gN'\otimes _{\gS[\frac{1}{p}]}\wS[\frac{1}{E}]\text{ est égal modulo }\gN'\otimes _{\gS[\frac{1}{p}]}\wS\text{ à } 
\\ -\vi_{\gN'}(\ta \xi'(T))=-\xi'(\vi_{D'}(T))\in 
\gN'\otimes _{\gS[\frac{1}{p}]}\big(\wS[\frac{1}{E}]/\wS\big). \end{gather*}
Calculons $Z$ en fonction de $X$. 
Comme 
$\xi=\begin{pmatrix} {\xi'} & Z
\\ 0 & 1\end{pmatrix}$ vérifie 
$$\xi
\begin{pmatrix} \vi_{D'} & Y
\\ 0 & 1\end{pmatrix}=
\begin{pmatrix} \vi_{\gN'} & X_{0}+X
\\ 0 & 1\end{pmatrix}\ta \xi,$$ on obtient  l'équation  $$Z=-\xi'(Y)+X_{0}+X+\vi_{\gN'}\ta Z.$$ 
On a donc 
$$Z=Z_{0}+(X+\vi_{\gN'}\ta X+\vi_{\gN'}\ta \vi_{\gN'}\tu{2} X+\dotss )$$
où $Z_{0}$ correspond à $X=0$ et où la somme dans le membre de droite converge au sens suivant : si $C$ est tel que $\vi_{\gN'}\in E^{-C}\Hom(\ta \gN',\gN')$ et $X\in \frac{u}{E^{C}}\gN'$, la série $\lambda^{C}(X+\vi_{\gN'}\ta X+\vi_{\gN'}\ta \vi_{\gN'}\tu{2} X+\dotss )$ converge dans $\gN'\otimes _{\gS[\frac{1}{p}]}\O$. 
Grâce au lemme~\ref{cont-epi-k} la proposition~\ref{prop3generalisee} résulte alors du lemme suivant (avec $\phi$ égal à la matrice de $\vi_{\gN'}$ dans une base de $\gN'$). \cqfd
\begin{lem}\label{point-fixegeneralise}
Soit $C\in \N^{*}$ et 
$\phi\in  p^{-C}E^{-C}M_{r}(\gS)$. Alors 
l'application
$$\theta: \big(\frac{u}{E^{C}}\gS[\frac{1}{p}]\big)^{r}\to 
\frac{1}{E^{C}}\wS^{r}/\wS^{r}$$
qui à $X$ associe $X+\phi\ta X+\phi\ta \phi\tu{2} X+\dotss  $ 
est surjective. 
\end{lem}
\noindent {\bf Démonstration.}
On rappelle que $E=u^{e}+c_{e-1}u^{e-1}+...+c_{0}$ avec $c_{i}\in pW$ pour $i\in \{1,...,e-1\}$ et $c_{0}=E(0)\in pW^{\times}$. 
On a $$\frac{u^{e}}{p}=a+\frac{E}{p}\text{ \ avec \ }a=-\big(\frac{c_{0}}{p}+...+\frac{c_{e-1}}{p}u^{e-1}\big)\in W^{\times}+u\gS\subset \gS^{\times}.$$ Donc  $\frac{u^{e}}{p}$ est un élément inversible de $\gS/E\gS$. 

On en déduit que $
\frac{u^{Ce}}{p}$ appartient à $\gS/E^{C}\gS\subset \gS[\frac{1}{p}]/E^{C}\gS[\frac{1}{p}]$ : la formule pour cet élément est 
$$\frac{u^{Ce}}{p}=a\Big(u^{(C-1)e}+u^{(C-2)e}E+...+E^{C-1}\Big)\in \gS/E^{C}\gS.$$ 
 On en déduit aussi qu'il existe un élément noté par abus 
$\frac{p^{C}}{u^{e}}\in \gS/E^{C}\gS$ tel que $u^{e}\frac{p^{C}}{u^{e}}=p^{C}$ dans $\gS/E^{C}\gS$  : la formule pour cet élément est $$\frac{p^{C}}{u^{e}}=
a^{-1}p^{C-1}-a^{-2}p^{C-2}E+...+(-1)^{C-1}a^{-C}E^{C-1}\in \gS/E^{C}\gS.$$

L'espace d'arrivée  $\frac{1}{E^{C}}\wS^{r}/\wS^{r}$ est   muni  de la topologie ``$u$-adique''
telle qu'une base de voisinage de $0$ soit formée par les 
$u^{k}\big(\frac{1}{E^{C}}\gS^{r}/\gS^{r}\big)$ 
 (comme dans le lemme~\ref{cont-epi-k}). 
 L'existence des éléments $\frac{u^{Ce}}{p} $ et $\frac{p^{C}}{u^{e}}$ dans $ \gS/E^{C}\gS$ montre qu'une base de voisinage de $0$ est donnée aussi par les 
 $p^{k}\big(\frac{1}{E^{C}}\gS^{r}/\gS^{r}\big)$. On en déduit aussi que 
 \begin{gather}\label{ukz-l}
\text{ pour tout \ }k\in \N, \text{ \ on a \ } 
\bigcup_{l\in \N}\frac{u^{k}}{p^{l}}\big(\frac{1}{E^{C}}\gS^{r}/\gS^{r}\big)
 =\frac{1}{E^{C}}\wS^{r}/\wS^{r}. \end{gather}
 
On écrit 
$\theta=\sum_{i=0}^{\infty}\theta_{i}$ où $$\theta_{i}:  \big(\frac{u}{E^{C}}\gS[\frac{1}{p}]\big)^{r}\to 
\frac{1}{E^{C}}\wS^{r}/\wS^{r}\text{ \ est défini par \ }\theta_{i}(X)=
\phi\ta \phi ... \tu{i-1}\phi   \tu{i} X .$$  En particulier $\theta_{0}(X)=X$. 
 
Soient    $k> C^{2}$ et  $l\in \N$.  On a 
\begin{gather}
\label{incl-theta-0}
\theta_{0}(\big(\frac{u^{ek}}{p^{l}E^{C}}\gS\big)^{r})=\frac{u^{ek}}{p^{l}}\Big(
\frac{1}{E^{C}}\gS^{r}/\gS^{r}\Big) \\  
\label{incl-theta-i}
\text{ et, pour  } i>0, \ \theta_{i}\big(\big(\frac{u^{ek}}{p^{l}E^{C}}\gS\big)^{r}\big)\subset \frac{u^{e(p^{i}k-2iC^{2})}}{p^{l}}\Big(
\frac{1}{E^{C}}\gS^{r}/\gS^{r}\Big).\end{gather}
On doit montrer \eqref{incl-theta-i}. Soit $i>0$. Comme $\phi\in  p^{-C}E^{-C}M_{r}(\gS)$ on a   
\begin{gather*}
\theta_{i}\big(\big(\frac{u^{ek}}{p^{l}E^{C}}\gS\big)^{r}\big)\subset 
\frac{1}{
p^{iC}\vi(E)^{C}... \vi^{i}(E)^{C}} \frac{u^{ep^{i}k}}{p^{l}}\Big(
\frac{1}{E^{C}}\gS^{r}/\gS^{r}\Big).
\end{gather*}
 Comme $\frac{\vi^{m}(E)}{p}\in \wS^{\times}$ pour tout $m\in \N^{*}$ et comme   
$\frac{u^{Ce}}{p}$ appartient à $ \gS/E^{C}\gS$ on en déduit \eqref{incl-theta-i}.  Comme $k> C^{2}$ on a
$$\inf_{i\in \N^{*}}p^{i}k-2iC^{2}=pk-2C^{2}$$ 
et on déduit de \eqref{incl-theta-i} que 
\begin{gather}\label{incl-theta-theta-0}(\theta-\theta_{0})(\big(\frac{u^{ek}}{p^{l}E^{C}}\gS\big)^{r})\subset \frac{u^{e(pk-2C^{2})}}{p^{l}}\Big(
\frac{1}{E^{C}}\gS^{r}/\gS^{r}\Big).\end{gather} 
Comme 
$k> C^{2}$ on a $e(pk-2C^{2})>ek$. On déduit donc de \eqref{incl-theta-0} 
et \eqref{incl-theta-theta-0}  que pour tout $k>C^{2}$ et $l\in \N$ on a 
\begin{gather}\label{incl-fin-theta}
\theta(\big(\frac{u^{ek}}{p^{l}E^{C}}\gS\big)^{r})=\frac{u^{ek}}{p^{l}}\Big(
\frac{1}{E^{C}}\gS^{r}/\gS^{r}\Big). 
\end{gather}
En fixant $k$, en  faisant varier $l$ et en utilisant \eqref{ukz-l} on termine alors la démonstration du lemme~\ref{point-fixegeneralise}. \cqfd

\begin{rem}\label{rem-tens-dual} On peut adapter la démonstration précédente au cas général où l'on ne suppose pas $(D'',\vi_{D''}, V_{D''})$ trivial, en prenant $Y\in D'\otimes {D''}^{*}$, $X\in u  \gN'\otimes {\gN''}^{*}[\frac{1}{E}]$, et en appliquant le lemme~\ref{point-fixegeneralise} à la matrice de $\vi_{\gN'}\otimes {}^{t}\vi_{\gN''}^{-1}$ dans une base de $\gN'\otimes {\gN''}^{*}$. \end{rem}

\section{Faiblement admissible implique admissible}\label{faibl-ad=ad}

Le but de ce paragraphe est d'achever la démonstration du théorème~\ref{fa--a}, en d'autres termes de montrer que le foncteur  $\mathbb D_{\mathrm{iso}}$ de la catégorie  des pseudo-iso-$\vi / \gS$-modules dans la catégorie $MHP(\vi)_{fa}$  des $\vi$-modules de Hodge-Pink faiblement admissibles est
essentiellement surjectif. 

On note $MHP(\vi)_{a}$ l'image essentielle de $\mathbb D_{\mathrm{iso}}$, et on appelle $\vi$-modules de Hodge-Pink admissibles les objets de $MHP(\vi)_{a}$. 
On commence par énoncer les propositions~\ref{prop2} et~\ref{prop1}. Puis on montrera le 
théorème~\ref{fa--a}  en admettant les propositions~\ref{prop2} et~\ref{prop1} (cela consistera en fait à se ramener au cas où $k$ est algébriquement clos et où les $\vi$-modules de Hodge-Pink sont anti-effectifs et irréductibles parmi les faiblement admissibles). 
Le reste du paragraphe sera consacré à la démonstration des 
propositions~\ref{prop2} et~\ref{prop1}.

Le lemme suivant est    une variante    du lemme~\ref{langtonOL}.

\begin{lem}\label{langton-zet}
Soit $\gM$ un $\gS$-module libre de rang $r$ et $V$ un $\wS$-réseau de $\gM\otimes _{\gS}\wS[\frac{1}{E}]$. Alors $\gM'=\big(\gM\otimes _{\gS}\gS[\frac{1}{E}]\big)\cap V$ est un $\gS$-module libre et  on a $$\gM'\otimes _{\gS}\gS[\frac{1}{E}]=\gM\otimes _{\gS}\gS[\frac{1}{E}]\text{ \  et  \ }\gM'\otimes  _{\gS}\wS=V.$$ Si $V\subset \gM\otimes _{\gS}\wS$ on a $\gM'\subset \gM$ et en notant $k$ la longueur du $\wS$-module $\gM\otimes _{\gS}\wS/V$ on a $\det(\gM')=E^{k}\det(\gM)$. L'application qui à $\gM$ et $V$ associe $\gM'$  commute aux opérations tensorielles, notamment  le déterminant. 
\end{lem}
\noindent{\bf Démonstration.} Soit $C\in \N$ tel que $E^{C} \gM\otimes _{\gS}\wS
  \subset V\subset E^{-C}\gM\otimes _{\gS}\wS$. Alors $\gM'$ contient $E^{C} \gM$ et  est déterminé par le fait que $\gM'/E^{C} \gM=V/E^{C} \gM\otimes _{\gS}\wS
$ dans $E^{-C}\gM\otimes _{\gS}\wS /E^{C} \gM\otimes _{\gS}\wS$. On en déduit facilement que $\gM'$ est réflexif, et donc libre car tout module de type fini réflexif sur un anneau local régulier de dimension $2$ est libre (lemme 6 de~\cite{serre-iwasawa}).  Les autres assertions sont faciles. 
 \cqfd

\begin{rem}
On pourrait aussi modifier $\gM$ en d'autres diviseurs  que $E$ (par exemple $\vi^k(E)$ pour $k\in \N$) ou en plusieurs diviseurs  simultanément (l'ordre n'importe pas  par factorialité de $\gS$).  \end{rem}

\begin{defi}\label{defibetagamma}
Soit $D=(D,\vi_{D},V_{D})$ un $\vi$-module de Hodge-Pink tel que 
$V_{D}\subset U_{D}$. 

On définit par récurrence, pour tout entier $n\in \N$, et tout W-réseau $\Delta$ de $D$ le  sous-$\gS$-module $\beta_{n}(\Delta)$ de 
$D\otimes_{K_0}\gS[\frac{1}{p}]$ 
de la fa\c con suivante. 
On pose $$\beta_{0}(\Delta)=\Delta \otimes_{W}\gS\text{ \  et \ } 
\beta_{n+1}(\Delta)=(\vi_{D}\otimes 1)(\ta \beta_{n}(\Delta)\cap V_{D}).$$ 

Puis on note   $\gamma_{n}(\Delta)$ la réduction de $\beta_{n}(\Delta)$ modulo $u $, en d'autres termes, 
$\gamma_{n}(\Delta)=\beta_{n}(\Delta)\otimes_{\gS}W$. 
\end{defi}

Il résulte du lemme~\ref{langton-zet} que pour tout $n\in \N$, 
$\beta_{n}(\Delta)$ est un $\gS$-module libre de rang $r$ et 
$$\beta_{n}(\Delta)[\frac{1}{p},\frac{1}{E}, ... ,\frac{1}{\vi^{n-1}(E)}]=
D\otimes_{K_0}\gS[\frac{1}{p},\frac{1}{E}, ... ,\frac{1}{\vi^{n-1}(E)}]$$ donc $\gamma_{n}(\Delta)$ est un $W$-réseau de $D$. Dans la suite on notera $\vi_{D}$ au lieu de $\vi_{D}\otimes 1$. On remarque que 
$\beta_{n}(\Delta)$ {\it n'est pas} un $\gS$-réseau  de $D\otimes_{K_0}\gS[\frac{1}{p}]$ si $n>0$ et $V_{D}\subsetneq U_{D}$. 
On a 
\begin{gather}\label{formule-synthetique-beta-n} \beta_{n}(\Delta)=\vi_{D}\dotss  \tu{n-1}\vi_{D}(\tu{n}\Delta\otimes _{W}\gS)
\cap \big(\vi_{D}\dotss  \tu{n-1}\vi_{D}\tu{n-1}V_{D}\oplus \dotss  \oplus  \vi_{D} V_{D}\big).\end{gather} 
Dans cette formule   $\tu{i}V_{D}$ est un $\wS_{i}$-réseau de 
$\tu{i+1} D\otimes _{K_0}\wS_{i}[\frac{1}{\vi^i(E)}]$ et  l'intersection a lieu dans 
$$D\otimes_{K_{0}} \wS_{n-1}\oplus ...\oplus  D\otimes_{K_{0}} \wS_{0}$$ où $\vi_{D}\dotss  \tu{n-1}\vi_{D}(\tu{n}\Delta\otimes _{W}\gS)$ s'envoie diagonalement.

\begin{prop}\label{prop2}
Soit $D=(D,\vi_{D},V_{D})$ un $\vi$-module de Hodge-Pink tel que $V_{D}\subset U_{D}$. 
Alors les deux assertions suivantes sont équivalentes : 
\begin{itemize}
\item i) pour tout $W$-réseau $\Delta$ de $D$ (ou pour un $W$-réseau $\Delta$) il existe une constante $C$ telle que pour tout $n\in \N$, on a  $p^{C}\Delta\subset \gamma_{n}(\Delta)\subset p^{-C}\Delta$. 
\item ii) $D$ est admissible. 
\end{itemize}
\end{prop}

On rappelle que si $(D,\vi_{D}),(D',\vi_{D'}),(D'',\vi_{D''})$ sont des $\vi$-modules de Hodge-Pink, une suite exacte courte $0\to D'\to D\to D''\to 0$ dans la catégorie des $\vi$-modules de Hodge-Pink est une suite exacte dans la catégorie des $\vi$-modules telle que $V_{D''}$ soit l'image de $V_{D}$ dans $D''\otimes_{K_0}\wS[\frac{1}{E}]$ et $V_{D'}$ l'intersection de $V_{D}$ avec $D'\otimes_{K_0}\wS[\frac{1}{E}]$. 
Dans cette situation on a clairement $t_{N}(D)=t_{N}(D')+t_{N}(D'')$ et $t_{H}(D)=t_{H}(D')+t_{H}(D'')$.
Le lemme suivant est immédiat. 

\begin{lem}\label{extfa}
Les deux assertions suivantes sont équivalentes : 
\begin{itemize}
\item i) $D$ est faiblement admissible et $t_{H}(D')=t_{N}(D')$,
\item ii) $D'$ et $D''$ sont faiblement admissibles.  \cqfd
\end{itemize} \end{lem}

Un $\vi$-module de Hodge-Pink $D=(D,\vi_{D},V_{D})$  faiblement admissible est irréductible dans la catégorie $MHP(\vi)_{fa}$ des $\vi$-modules de Hodge-Pink faiblement admissibles si et seulement si $t_{N}(D)=t_{H}(D)$ et  pour tout sous-$\vi$-module $D'$ autre que $0$ et $D$, $t_{H}(D')<t_{N}(D')$. 

\begin{prop}\label{prop1} On suppose $k$ algébriquement clos. 
Soit $D=(D,\vi_{D},V_{D})$ un objet irréductible dans la catégorie $MHP(\vi)_{fa}$  et tel que 
$V_{D}\subset U_{D}$. 

Alors pour tout réseau $\Delta$ de $D$ (ou de fa\c con équivalente pour un réseau $\Delta$ de $D$) il existe une constante $C\in \N$ telle que pour tout $n\in \N$, on a  $p^{C}\Delta\subset \gamma_{n}(\Delta)\subset p^{-C}\Delta$. \end{prop}

\noindent {\bf Début de la démonstration du théorème~\ref{fa--a}  en admettant les propositions~\ref{prop2} et~\ref{prop1}.} 
On a vu au paragraphe précédent que le foncteur $\mathbb D_{\mathrm{iso}}$ de la catégorie  des pseudo-iso-$\vi / \gS$-modules dans la catégorie $MHP(\vi)_{fa}$ des $\vi$-modules de Hodge-Pink faiblement admissibles est pleinement fidèle. 
Il reste donc à montrer 

\noindent $(A1)$ : tout objet $(D,\vi_{D},V_{D})$ dans la catégorie $MHP(\vi)_{fa}$ est  admissible. 

On considère les énoncés suivants, dont chacun est plus faible que le précédent.  

\noindent $(A2)$ : tout objet $(D,\vi_{D},V_{D})$ dans la catégorie $MHP(\vi)_{fa}$ qui vérifie $V_{D}\subset U_{D}$ est  admissible. 

\noindent $(A3)$ : si $k$ est algébriquement clos, tout objet $(D,\vi_{D},V_{D})$ dans la catégorie $MHP(\vi)_{fa}$ qui vérifie $V_{D}\subset U_{D}$ est  admissible. 

\noindent $(A4)$ : si $k$ est algébriquement clos, tout objet irréductible $(D,\vi_{D},V_{D})$ dans la catégorie $MHP(\vi)_{fa}$ qui vérifie $V_{D}\subset U_{D}$ est  admissible. 

On va montrer $(A4)$, puis $(A3)$, puis $(A2)$, puis $(A1)$. 

D'abord $(A4)$ 
   résulte de la  proposition~\ref{prop1} et de $i) \Rightarrow ii)$ dans la proposition~\ref{prop2}. 

 On montre maintenant $(A3)$ à l'aide de $(A4)$. Soit $D=(D,\vi_{D},V_{D})$ comme dans $(A3)$. 
 D'après le lemme~\ref{extfa}, $D$ 
admet  en tant que $\vi$-module de Hodge-Pink une filtration 
$0=D_{0}\subset D_{1}\subset \dotss  \subset D_{r}= D$  où  chaque quotient $D_{i}/D_{i-1}$ 
est irréductible dans  la catégorie $MHP(\vi)_{fa}$ et vérifie $V_{D_{i}/D_{i-1}}\subset U_{D_{i}/D_{i-1}}$. Grâce à $(A4)$, chaque quotient $D_{i}/D_{i-1}$  est admissible. Or la proposition~\ref{prop3generalisee}  montre qu'une extension de $\vi$-modules de Hodge-Pink admissibles est admissible et donc 
 $D$ est admissible. On a montré $(A3)$. 

 On montre maintenant $(A2)$ à l'aide de $(A3)$. Soit $D=(D,\vi_{D},V_{D})$ comme dans $(A2)$. 
Le  lemme suivant est bien connu des spécialistes. 
On note $\overline{\gS} $ et $\widehat{\overline \gS} $ les anneaux définis comme $\gS$ et $\wS$ en rempla\c cant $k$ par $\overline k$, c'est-à-dire que 
$\overline{\gS} =W(\overline k)[[u]]$ et $\widehat{\overline \gS} $ est le complété de $\overline{\gS} [\frac{1}{p}]$ pour la topologie $E$-adique. On appelle $\vi $-module de Hodge-Pink sur $\overline k$ un triplet 
$\overline D=(\overline D,\vi_{\overline D},V_{\overline D})$, où 
$\overline D$ est un $W(\overline k)[\frac{1}{p}]$-espace vectoriel de dimension finie, $\vi_{\overline D}:\ta \overline D\to \overline D$ est un isomorphisme de $W(\overline k)[\frac{1}{p}]$-espaces vectoriels  
 et $V_{\overline D}$ est un $\widehat{\overline \gS} $-module libre qui est un réseau 
 dans $\ta {\overline D}\otimes_{W(\overline k)[\frac{1}{p}]}\widehat{\overline \gS} [\frac{1}{E}]$.

\begin{lem}\label{galois} 
Soit  $D$ un $\vi $-module de Hodge-Pink, et $\overline  D=D\otimes _{K_0}W(\overline k)[\frac{1}{p}]$ le $\vi $-module de Hodge-Pink sur $\overline k$ qui s'en déduit (on a par exemple $V_{\overline  D}=V_{D}\otimes _{\widehat{ \gS}}\widehat{\overline \gS}$). Alors si $D$ est faiblement admissible, $\overline  D$ est faiblement admissible. 
\end{lem}
 \noindent {\bf Démonstration.} On raisonne par l'absurde. Soit $\overline D_{0}\subset \overline  D$ un sous-$\vi$-module  sur $\overline  k$ tel que $t_{H}(\overline D_{0})>t_{N}(\overline D_{0})$, et minimal pour cette propriété, ce qui fait que tout quotient $\overline D_{0}''$ de $\overline D_{0}$ vérifie $t_{H}(\overline D_{0}'')>t_{N}(\overline D_{0}'')$. 
Alors $\overline D_{1}=\sum_{\gamma \in \mathrm{Gal}(\overline  k/k)}\gamma(\overline D_{0})$ provient d'un sous-$\vi$-module $D_{1}$ de $D$   (parce que $k=\overline  k^{\mathrm{Gal}(\overline  k/k)}$), et en tant que $\vi$-module de Hodge-Pink sur $\overline  k$, $\overline D_{1}$  est une extension successive de quotients de $\gamma(\overline D_{0})$ pour $\gamma \in \mathrm{Gal}(\overline  k/k)$, donc on a  $t_{H}(D_{1})=t_{H}(\overline  D_{1})>t_{N}(\overline  D_{1})=t_{N}(D_{1})$, ce qui contredit la faible admissibilité de $D$.  \cqfd

\noindent {\bf Fin de la démonstration du théorème~\ref{fa--a}  en admettant les propositions~\ref{prop2} et~\ref{prop1}.} On termine la preuve de $(A2)$ à l'aide de $(A3)$. Grâce au lemme~\ref{galois},  $\overline  D=D\otimes _{K_0}W(\overline k)[\frac{1}{p}]$ est faiblement admissible. Grâce à $(A3)$, $\overline  D$  
est admissible en tant que $\vi$-module de Hodge-Pink sur $\overline  k$ (c'est-à-dire qu'il est associé à un pseudo-iso-$\vi / \gS$-module sur $\overline  \gS$). 
Par  $ii) \Rightarrow i)$ de la proposition~\ref{prop2} (appliquée à $\overline k$ au lieu de $k$) il satisfait  la conclusion $i)$ de la proposition~\ref{prop2}. Il en résulte immédiatement que $D$ satisfait la conclusion $i)$ de la proposition~\ref{prop2} et par $i) \Rightarrow ii)$ de la proposition~\ref{prop2}, $D$ est admissible. On  a montré $(A2)$.

 Enfin   $(A1)$ résulte de $(A2)$. 
 En effet, 
on s'y ramène en  multipliant $\vi_{D}$ par $p^{-k}$ et $V_{D}$ par $E^{k}$ pour un entier $k$ assez grand (cela ne change pas la faible admissibilité, et si $(D,p^{-k}\vi_{D},E^{k}V_{D})$  est associé à un pseudo-iso-$\vi / \gS$-module $(\gN,\vi_{\gN})$, $(D,\vi_{D},V_{D})$ est associé 
à $(\gN,E^{k}\vi_{\gN})$). 
On a 
terminé la démonstration du théorème~\ref{fa--a}.  \cqfd

Nous passons maintenant à la démonstration des propositions~\ref{prop2} et~\ref{prop1}. 
Nous commen\c cons par des propriétés générales de $\beta_{n}$ et $\gamma_{n}$. 
Regardons déjà ce qui se passe quand $r=\dim V$ vaut $1$. Plus précisément prenons   $t_{H}\in \Z_{\leq 0}$, $t_{N}\in \Z$ et posons $D=K_0$, $\vi_{D}=p^{t_{N}}$, $V_{D}=E^{-t_{H}}U_{D}$ et $\Delta=W$. On a alors 
$$\beta_{n}(\Delta)=p^{n(t_{N}-t_{H})}\Big(\prod_{m=0}^{n-1}\vi^{m}\Big(\frac{E}{E(0)}\Big)\Big)^{-t_{H}}\gS\text{ \ 
et  \ }\gamma_{n}(\Delta)=p^{n(t_{N}-t_{H})}W.$$ On voit que la condition i) de la proposition~\ref{prop2} est vérifiée si et seulement si $t_{N}=t_{H}$. On suppose $t_{N}=t_{H}=t\in \Z_{\leq 0}$. Alors $(D,\vi_{D},V_{D})$ est admissible (comme le prévoit la proposition~\ref{prop2}) et il est associé au pseudo-iso-$\vi / \gS$-module $(\gN,\vi_{\gN})=\big(\gS[\frac{1}{p}], \big(\frac{p}{E(0)}E\big)^{t}\big)$. Si on note $\gM=\gS\subset \gN$ on a $\xi^{-1}(\gM)=\lambda^{-t}\gS=\tu{n}\lambda^{-t}\beta_{n}(\Delta)$. On voit que la suite $\beta_{n}(\Delta)$ ``converge'' vers $\xi^{-1}(\gM)$
car $\vi^{n}(\lambda)$ est une unité dans $\gS[[\frac{u^{ep^{n}}}{p}]]$. 
Plus généralement pour tout $\gS$-module $\gM$ libre de rang $1$ muni d'un isomorphisme $\gM[\frac{1}{p}]=\gN$ il existe $C$ tel que pour tout $n$ on ait 
\begin{gather}\label{incl-rang1}p^{C} \beta_{n}(\Delta)[[\frac{u^{ep^{n}}}{p}]]
\subset 
\xi^{-1}(\gM)[[\frac{u^{ep^{n}}}{p}]]
\subset 
p^{-C} \beta_{n}(\Delta) [[\frac{u^{ep^{n}}}{p}]].\end{gather} 
Dans la formule précédente on a utilisé la convention suivante, qui servira dans toute la suite : 
$$\text{ \ 
si \ }Q \text{ \ est un  \ }\gS\text{-module libre on note \ }Q[[\frac{u^{ep^{n}}}{p}]]=Q\otimes_{\gS}\gS[[\frac{u^{ep^{n}}}{p}]].$$
 On verra dans la remarque~\ref{rem-incl-ii-impl-i} que 
 les inclusions (\ref{incl-rang1}) ont lieu  pour tout pseudo-iso-$\vi / \gS$-module $(\gN,\vi_{\gN})$  et elles joueront un rôle heuristique fondamental  dans la preuve de ``i) implique ii)'' dans la proposition~\ref{prop2}.

Le lemme suivant établit  quelques propriétés satisfaites par $\beta_{n}(\Delta)$ et $\gamma_{n}(\Delta)$. 
\begin{lem}\label{propbetagamma}
Soient $(D,\vi_{D},V_{D})$ et  $\Delta$ comme dans la définition~\ref{defibetagamma}. On suppose $t_{N}(D)=t_{H}(D)$. 
Alors il existe une constante $C$ telle que pour tout $n\in \N$, 
\begin{itemize}
\item a)
$\beta_{n+1}(\Delta)\subset (\vi_{D}\otimes 1)(\ta \beta_{n}(\Delta))$ et $\beta_{n+1}(\Delta)\subset p^{-C}\beta_{n}(\Delta)$, 
\item b) les réseaux $\det(\gamma_{n}(\Delta))$ et $\det(\Delta)$ dans $\det(D)$ sont égaux et deux réseaux de $D$ quelconques $U$ et $U'$ parmi 
$\gamma_{n}(\Delta)$, $\gamma_{n+1}(\Delta)$, $\vi_{D}(\ta \gamma_{n}(\Delta))$, $\vi_{D}(\ta \gamma_{n+1}(\Delta))$ vérifient $U\subset p^{-C}U'$, 
\item c) pour tout entier $k\in \N$, 
$$\beta_{n}(\Delta)\text{ mod }u^{p^{k}}\subset p^{-kC}\gamma_{n}(\Delta)\otimes_{W}\gS \text{ mod }u^{p^{k}},$$
\item d) on a \begin{gather*} p^{C}\gamma_{n}(\Delta)\otimes_{W}
\gS[[\frac{u^{e}}{p}]]\subset  \beta_{n}(\Delta)[[\frac{u^{e}}{p}]] \subset  p^{-C}\gamma_{n}(\Delta)\otimes_{W}
\gS[[\frac{u^{e}}{p}]],\end{gather*}
\item e)
 pour tout $n'\in \N$  on a 
$$p^{C}\gamma_{n+n'}(\Delta)\subset \gamma_{n}\circ \gamma_{n'}(\Delta)\subset p^{-C}\gamma_{n+n'}(\Delta).$$
\end{itemize}
\end{lem}
\noindent {\bf Démonstration.}
D'abord l'inclusion évidente 
$\beta_{n+1}(\Delta)\subset (\vi_{D}\otimes 1)(\ta \beta_{n}(\Delta))$ implique  
$\gamma_{n+1}(\Delta)\subset \vi_{D}(\ta \gamma_{n}(\Delta))$. Il existe $C\in \N$ tel que $\vi_{D}(\ta \Delta)\subset p^{-C}\Delta$. Cela implique $\beta_{1}(\Delta)\subset p^{-C}\beta_{0}(\Delta)$. On a alors pour tout $n$, $\beta_{n+1}(\Delta)\subset p^{-C}\beta_{n}(\Delta)$, ce qui termine la preuve de a). On en déduit pour tout $n\in \N$,   
$\gamma_{n+1}(\Delta)\subset p^{-C}\gamma_{n}(\Delta)$. 
La construction de $\beta_{n}(\Delta)$ et de $\gamma_{n}(\Delta)$ commute au déterminant : le fait que l'intersection avec $V_{D}$ commute au déterminant, et en fait à toute opération tensorielle, n'est pas évident, c'est une conséquence du lemme~\ref{langton-zet}. 
 On a donc 
$$\det(\beta_{n}(\Delta))=\Big(\prod_{m=0}^{n-1}\vi^{m}\Big(\frac{E}{E(0)}\Big)\Big)^{-t_{N}(D)}\det(\Delta)\otimes_{W}\gS$$ comme sous-$\gS$-module de $\det(D)\otimes_{K_0}\gS[\frac{1}{p}]$, et donc 
$\det(\gamma_{n}(\Delta))=\det(\Delta)$ comme $W$-réseau de 
$\det(D)$. Comme  $$\gamma_{n+1}(\Delta)\subset \vi_{D}(\ta \gamma_{n}(\Delta)) \ , \ \gamma_{n+1}(\Delta)\subset p^{-C}\gamma_{n}(\Delta)\text{ \  et  \ } \det(\gamma_{n}(\Delta))=\det(\Delta)$$
il existe 
une nouvelle constante $C$  indépendante de $n$ telle  que  la propriété b) du lemme  soit vraie.

On déduit de l'inclusion évidente $\beta_{n+1}(\Delta)\subset (\vi_{D}\otimes 1)(\ta \beta_{n}(\Delta))$  que pour tout $n$ on a 
\begin{gather*}\beta_{n}(\Delta)\text{ mod }u^{p}\subset \vi_{D}(\ta\gamma_{n-1}(\Delta))\otimes_{W}\gS \text{ mod }u^{p},
\\
\beta_{n}(\Delta)\text{ mod }u^{p^{2}}\subset \vi_{D}\ta\vi_{D}(\tu{2}\gamma_{n-2}(\Delta))\otimes_{W}\gS \text{ mod }u^{p^{2}}...\end{gather*}
Grâce à b), la propriété c) en découle. 

Par la propriété c) il existe une constante $C$  indépendante de $n$ telle  que  
$$\beta_{n}(\Delta) \subset p^{-C}\gamma_{n}(\Delta)\otimes_{W}
\gS[[\frac{u^{e}}{p}]].$$ Mais le déterminant de $\beta_{n}(\Delta)$ et celui de $\gamma_{n}(\Delta)$ diffèrent par 
$$\Big(\prod_{m=0}^{n-1}\vi^{m}\Big(\frac{E}{E(0)}\Big)\Big)^{t_{N}(D)}$$ qui est une unité de 
$\gS[[\frac{u^{e}}{p}]]$. 
On en déduit qu'il existe une constante $C$  indépendante de $n$ telle  que  
 d) soit vrai.

La propriété   e) est plus difficile à démontrer. 
 Elle jouera un grand rôle dans la démonstration de la proposition~\ref{prop1}. 
Bien sûr il suffit de montrer l'une des deux inclusions car 
les réseaux 
$ \gamma_{n}\circ \gamma_{n'}(\Delta)$ et $\gamma_{n+n'}(\Delta)$
ont des  déterminants égaux. C'est l'inclusion de droite que nous allons démontrer. 

Par \eqref{formule-synthetique-beta-n} on a 
\begin{gather}\nonumber \beta_{n}(\gamma_{n'}(\Delta))=\vi_{D}\dotss  \tu{n-1}\vi_{D}(\tu{n}\gamma_{n'}(\Delta)\otimes _{W}\gS)
\\ \label{betan-gamman'}
\cap \big(\vi_{D}\dotss  \tu{n-1}\vi_{D}\tu{n-1}V_{D}\oplus \dotss  \oplus  \vi_{D} V_{D}\big).\end{gather} 
%Dans cette formule   $\tu{i}V_{D}$ est un $\wS_{i}$-réseau de 
%$\tu{i+1} D\otimes _{K_0}\wS_{i}[\frac{1}{\vi^i(E)}]$ et  l'intersection a lieu %dans 
%$$D\otimes_{K_{0}} \wS_{n-1}\oplus ...\oplus  D\otimes_{K_{0}} \wS_{0}$$ où $\vi_{D}\dotss  \tu{n-1}\vi_{D}(\tu{n}\gamma_{n'}(\Delta)\otimes _{W}\gS)$ s'envoie diagonalement. 
%$$\subset 
%\vi_{D}\dotss  \tu{n-1}\vi_{D}(\tu{n}\gamma_{n'}(\Delta)\otimes _{W}\gS).$$
D'après  d) on a 
 \begin{gather}\label{incl-consequence-d}\tu{n} \gamma_{n'}(\Delta)\otimes_{W}
\gS\subset p^{-C }\ \tu{n}\beta_{n'}(\Delta)[[\frac{u^{ep^{n}}}{p}]]\end{gather}  et \eqref{betan-gamman'} implique alors 
\begin{gather}\nonumber 
\beta_{n}(\gamma_{n'}(\Delta))\subset p^{-C } \vi_{D}... \tu{n-1}\vi_{D}(\tu{n}\beta_{n'}(\Delta)[[\frac{u^{ep^{n}}}{p}]])
\\
\label{incl-betangamman'Delta}
\cap \big(\vi_{D}... \tu{n-1}\vi_{D}\tu{n-1}V_{D}\oplus ... \oplus  \vi_{D} V_{D}\big).\end{gather}
On rappelle que   $s\in \Z_{\leq 0}$ est  tel que $V_{D}\supset E^{-s}U_{D}$. 
On applique  
le lemme~\ref{perte-modules} ci-dessous  à 
$$\M=\vi_{D}\dotss  \tu{n-1}\vi_{D}(\tu{n}\beta_{n'}(\Delta)), \ \V_{0}=\vi_{D}V_{D},..., \ \V_{n-1}=\vi_{D}\dotss  \tu{n-1}\vi_{D}\tu{n-1}V_{D}.$$ Comme 
$\M\cap (\V_{0}\oplus ... \oplus \V_{n-1})=\beta_{n+n'}(\Delta)$ le lemme~\ref{perte-modules} implique 
\begin{gather}\label{conseq-lemme-perte-modules}\big(\M[[\frac{u^{ep^{n}}}{p}]]\big)\cap (\V_{0}\oplus ... \oplus \V_{n-1})
\subset p^{-C_{1}}\beta_{n+n'}(\Delta)[[\frac{u^{ep^{n}}}{p}]].\end{gather}
L'inclusion \eqref{incl-betangamman'Delta} se réécrit 
\begin{gather}
\label{incl-betangamman'Delta2}
\beta_{n}(\gamma_{n'}(\Delta))\subset p^{-C }\big( \M[[\frac{u^{ep^{n}}}{p}]]\big)\cap (\V_{0}\oplus ... \oplus \V_{n-1}).\end{gather}
Les inclusions  \eqref{conseq-lemme-perte-modules}  et \eqref{incl-betangamman'Delta2} impliquent 
 \begin{gather*}\beta_{n}(\gamma_{n'}(\Delta))\subset p^{-C-C_{1}}
\beta_{n+n'}(\Delta)[[\frac{u^{ep^{n}}}{p}]]
 \\
 \text{d'où\ \ \ }\gamma_{n}(\gamma_{n'}(\Delta))=\beta_{n}(\gamma_{n'}(\Delta))\ \mod\ u \subset p^{-C-C_{1}}\gamma_{n+n'}(\Delta),\end{gather*} ce qui achève la démonstration 
du lemme~\ref{propbetagamma}.   \cqfd
%de la proposition~\ref{prop1}. 
\begin{lem}\label{perte-modules}
Pour tout $s\in \Z_{\leq 0}$  il existe une constante 
$C_{1}$ telle que l'énoncé suivant soit vrai. 
Soit $n\in \N$, $\M$ un $\gS$-module libre de rang $r$, et pour $i\in \{0,...,n-1\}$ soit $\V_{i}$ un $\wS_{i}$-module  libre de rang $r$  vérifiant 
$$\vi^i(E)^{-s}\M\otimes_{\gS}\wS_{i}
\subset \V_{i} \subset
\M\otimes_{\gS}\wS_{i}.$$
Alors $\M\cap (\V_{0}\oplus ... \oplus \V_{n-1})$ est un $\gS$-module libre de rang $r$ et on a 
\begin{gather}
\nonumber 
\big(\M\cap (\V_{0}\oplus ... \oplus \V_{n-1})\big)[[\frac{u^{ep^{n}}}{p}]]\subset 
\big(\M[[\frac{u^{ep^{n}}}{p}]]\big)\cap (\V_{0}\oplus ... \oplus \V_{n-1})
\\ \label{incl-perte-modules}
\subset p^{-C_{1}}\big(\M\cap (\V_{0}\oplus ... \oplus \V_{n-1})\big)[[\frac{u^{ep^{n}}}{p}]].\end{gather}
\end{lem}
 \noindent {\bf Démonstration.}
D'après le lemme~\ref{langton-zet}, 
 $\check \M=\M\cap (\V_{0}\oplus ... \oplus \V_{n-1})$
est un $\gS$-module libre de rang $r$ et    
\begin{gather}
\label{form-Vi}\V_{i}=\check \M\otimes_{\gS}\wS_{i}\text{ \  pour tout \ }i\in \{0,...,n-1\}.\end{gather} La première inclusion de \eqref{incl-perte-modules} est évidente et pour la deuxième on remarque que, grâce à  \eqref{form-Vi}, 
\begin{gather*}\big(\M[[\frac{u^{ep^{n}}}{p}]]\big)\cap (\V_{0}\oplus ... \oplus \V_{n-1})
\\
\subset \check \M \otimes_{\gS}
\Big(\Big(\big(E\dotss  \vi^{n-1}(E)\big)^{s}
 \gS[[\frac{u^{ep^{n}}}{p}]]\Big)\cap \big(\wS_{n-1}\oplus\dotss  \oplus \wS_0\big)\Big)
 \end{gather*}
et on applique le lemme suivant. 
\cqfd

\begin{lem}\label{perte}
Pour tout $s\in \Z_{\leq 0}$  il existe une constante 
$C_{1}$, telle que pour tout $n\in \N$, 
\begin{gather*}\Big(\big(E\dotss  \vi^{n-1}(E)\big)^{s}
 \gS[[\frac{u^{ep^{n}}}{p}]]\Big)\cap \big(\wS_{n-1}\oplus\dotss  \oplus \wS_0\big)
 \subset p^{-C_{1}}
 %\big((E^{p^{n-a+1}})\dotss  (E^{p^{n-1}})\big)^{s}
 \gS[[\frac{u^{ep^{n}}}{p}]].\end{gather*}
 \end{lem}
 \noindent {\bf Démonstration.}
On prend pour $C_{1}$ 
le plus petit entier supérieur ou égal à  $1+(-s)(\frac{1}{p}+\frac{1}{p^{2}}+\dotss  )$. 
On a  
\begin{gather*}\Big(\big(E\dotss  \vi^{n-1}(E)\big)^{s}
 \gS[[\frac{u^{ep^{n}}}{p}]]\Big)\cap \big(\wS_{n-1}\oplus\dotss  \oplus \wS_0\big)
 \subset  \gS[[\frac{u^{ep^{n}}}{p}]][\frac{1}{p}].\end{gather*}
 En effet pour tout $i\in \{0,...,n-1\}$ on a 
 $$\vi^i(E)^{s}
 \gS[[\frac{u^{ep^{n}}}{p}]][\frac{1}{p}]\cap \wS_{i} = \gS[[\frac{u^{ep^{n}}}{p}]][\frac{1}{p}].$$ 
 On possède  sur 
$\gS[[\frac{u^{ep^{n}}}{p}]][\frac{1}{p}]$
 la norme 
$|.|_{|\pi_{K}|^{p^{-n}}}$.   On a alors pour tout entier $r\in \Z$, 
$$p^{r}\gS[[\frac{u^{ep^{n}}}{p}]]\subset \{x\in \gS[[\frac{u^{ep^{n}}}{p}]][\frac{1}{p}],|x|_{|\pi_{K}|^{p^{-n}}}\leq |p|^{ r}\}\subset p^{r-1}\gS[[\frac{u^{ep^{n}}}{p}]].$$
Grâce à \eqref{vi-m-n-r}, le lemme  en résulte.  \cqfd

\noindent{\bf Démonstration de ``ii) implique i)'' dans  la proposition~\ref{prop2}.} 
Soit $(\gM,\vi_{\gM})$ un $\vi / \gS$-module dont l'image par $\mathbb D_{\mathrm{iso}}$  est $D=(D,\vi_{D},V_{D})$.  
Pour vérifier i) choisissons $\Delta=\gM \ \mod \ u$. 
Comme $\O[\frac{1}{\lambda}]\subset \gS[[\frac{u^{e}}{p}]][\frac{1}{p}]$, le 
  lemme~\ref{rigidite2} montre qu'il existe  $C$ tel que  $\xi^{-1}(\gM)\subset p^{-C}\beta_{0}(\Delta)[[\frac{u^{e}}{p}]]$. 
Quitte à augmenter $C$ on suppose de plus  que l'énoncé du  lemme~\ref{propbetagamma} est  vrai pour cette valeur de $C$ (le lemme~\ref{propbetagamma} s'applique car $t_{N}(D)=t_{H}(D)$ puisque $D$ est admissible). 
 
 On en déduit 
$\tu{n}\big(\xi^{-1}(\gM)\big)\subset p^{-C}\ \tu{n}\beta_{0}(\Delta)[[\frac{u^{ep^{n}}}{p}]]$. 
Par factorialité de $\gS$ on a $$\gM=\vi_{\gM}\big(\ta \gM\cap (\vi_{\gM}^{-1}(\gM)\otimes _{\gS}\wS)\big). $$ Autrement dit 
   $\xi^{-1}(\gM)=\vi_{D}(\ta (\xi^{-1}(\gM))\cap V_{D})$ d'où 
   \begin{gather}\nonumber \xi^{-1}(\gM)=
   \vi_{D}\dotss  \tu{n-1}\vi_{D}\Big(
 \tu{n}(\xi^{-1}(\gM))\Big) \cap\big( \vi_{D}\dotss  \tu{n-1}\vi_{D}\tu{n-1}V_{D}\oplus\dotss  \oplus \vi_{D}V_{D}\big) 
\\ \!\!\! \subset p^{-C}\vi_{D}\dotss  \tu{n-1}\vi_{D}\Big(
 \tu{n}\beta_{0}(\Delta)\Big)[[\frac{u^{ep^{n}}}{p}]]
\label{incl-RM-ppppphi-VVVVVV} \cap\big( \vi_{D}\dotss  \tu{n-1}\vi_{D}\tu{n-1}V_{D}\oplus\dotss  \oplus \vi_{D}V_{D}\big).\end{gather}
On applique  le lemme~\ref{perte-modules}  à $$\M=\vi_{D}\ta \vi_{D} \dotss  {}\tu{n-1}\vi_{D}\Big(
 \tu{n}\beta_{0}(\Delta)\Big), \ \V_{0}=\vi_{D}V_{D},..., \ \V_{n-1}=\vi_{D}\dotss  \tu{n-1}\vi_{D}\tu{n-1}V_{D}.$$
L'inclusion \eqref{incl-RM-ppppphi-VVVVVV} se réécrit $$
 \xi^{-1}(\gM)
\subset p^{-C}\M[[\frac{u^{ep^{n}}}{p}]]\cap (\V_{0}\oplus ... \oplus \V_{n-1})$$
 et par définition de  $\beta_{n}(\Delta)$ on a  
$$\beta_{n}(\Delta)=\M\cap (\V_{0}\oplus ... \oplus \V_{n-1}).$$
Donc le lemme~\ref{perte-modules}  implique 
    \begin{gather}\label{inclu-RM-betan-C2} \xi^{-1}(\gM)
\subset p^{-C-C_{1}}
\beta_{n}(\Delta)[[\frac{u^{ep^{n}}}{p}]].\end{gather}
 On en déduit  $\Delta=\xi^{-1}(\gM)\ \mod\  u \subset p^{-C-C_{1}}\gamma_{n}(\Delta)$. Ceci montre l'une des inclusions de i) de la proposition~\ref{prop2}, 
mais l'autre en résulte car $t_{N}(D)=t_{H}(D)$ et donc $\det(\gamma_{n}(\Delta))=\det(\Delta)$. Ceci achève la démonstration de ``ii) implique i)'' dans la proposition~\ref{prop2}. 

\begin{rem} \label{rem-incl-ii-impl-i} L'inclusion (\ref{inclu-RM-betan-C2}) et un argument de déterminants montrent qu'il existe $C$ tel que pour tout $n\in \N$ on ait 
\begin{gather}\label{incl-rem-ii-impl-i}p^{C} \beta_{n}(\Delta)[[\frac{u^{ep^{n}}}{p}]]
\subset 
\xi^{-1}(\gM)[[\frac{u^{ep^{n}}}{p}]]
\subset 
p^{-C} \beta_{n}(\Delta) [[\frac{u^{ep^{n}}}{p}]].\end{gather} 
On pourrait en déduire une autre preuve de la proposition~\ref{prop-HP-fa} ``admissible implique faiblement admissible''. 
Mais le plus important est que (\ref{incl-rem-ii-impl-i})  justifie  a posteriori la construction du pseudo-iso-$\vi/\gS$-module  $(\gN,\vi_{\gN})$ dans la preuve de ``i) implique ii)'' que nous allons donner maintenant. \end{rem}

\noindent{\bf 
Début de la démonstration de ``i) implique ii)'' dans la proposition~\ref{prop2}.}
On a $\det(\gamma_{n}(\Delta))=p^{n(t_{N}(D)-t_{H}(D))}$. La condition i) entraîne que $t_{N}(D)=t_{H}(D)$, ce qui permet d'appliquer le lemme~\ref{propbetagamma}. 
 On fixe un réseau $\Delta$ de $D$. 
 
\begin{lem}\label{comparaison-beta}
Sous la condition  i) de la proposition~\ref{prop2}  il existe $C$ tel que pour tous les entiers $m,n\in \N$ avec $n\geq m$, on ait 
\begin{gather}
\label{form-comparaison-beta1}
p^{C} \beta_{m}(\Delta)[[\frac{u^{ep^{m}}}{p}]]
\subset 
\beta_{n}(\Delta)[[\frac{u^{ep^{m}}}{p}]]
\subset 
p^{-C} \beta_{m}(\Delta) [[\frac{u^{ep^{m}}}{p}]]. 
% \\ \label{form-comparaison-beta2}\text{ et }\ \ p^{C} \beta_{m}(\Delta)\ \text{mod}\ u^{ep^{m}}
%\subset 
%\beta_{n}(\Delta)\ \text{mod}\ u^{ep^{m}}\subset 
%p^{-C} \beta_{m}(\Delta)\ \text{mod}\ u^{ep^{m}}.
\end{gather}
\end{lem}
\begin{rem}\label{rem-cauchy-converge}
Si on considère \eqref{incl-rem-ii-impl-i} comme une propriété de ``convergence'' de la suite $\beta_{n}(\Delta)[[\frac{u^{ep^{n}}}{p}]]
$ vers $\xi^{-1}(\gM)[\frac{1}{p}]$, le lemme précédent affirme que (sous l'hypothèse i))
 la suite $\beta_{n}(\Delta)[[\frac{u^{ep^{n}}}{p}]]
$ est ``de Cauchy''. Le lemme~\ref{lem-Mn-Sigma[[]]} ci-dessous affirme que si une suite est ``de Cauchy'', elle ``converge'' et la combinaison des lemmes~\ref{comparaison-beta}  et~\ref{lem-Mn-Sigma[[]]} permettra de construire le pseudo-iso-$\vi/\gS$-module  $(\gN,\vi_{\gN})$ comme ``limite'' de la suite $\beta_{n}(\Delta)[[\frac{u^{ep^{n}}}{p}]]
$. 

\end{rem}
 \noindent {\bf Démonstration du lemme~\ref{comparaison-beta}.}
%Les trois termes de (\ref{form-comparaison-beta2}) sont des $(\gS/u^{ep^{m}}\gS)$-modules libres de rang $r$ dans 
%$D\otimes _{K_0}(\gS/u^{ep^{m}}\gS)[\frac{1}{p}]$. 
%La  suite d'inclusions (\ref{form-comparaison-beta2}) résulte de (\ref{form-comparaison-beta1}) par réduction modulo l'idéal de $\gS[[\frac{u^{ep^{m}}}{p}]]$ engendré par $\frac{u^{ep^{m}}}{p}$. 
Pour montrer l'inégalité de droite de (\ref{form-comparaison-beta1}) 
on applique  le lemme~\ref{perte-modules}  avec $m$ à la place de $n$ et  $$\M=\vi_{D}\dotss  \tu{m-1}\vi_{D}(\tu{m}\beta_{0}(\Delta)), \ \V_{0}=\vi_{D}V_{D},..., \ \V_{m-1}=\vi_{D}\dotss  \tu{m-1}\vi_{D}\tu{m-1}V_{D}.$$
On a \begin{gather}\label{betamDelta}
\beta_{m}(\Delta)=\M \cap (\V_{0}\oplus ... \oplus \V_{m-1})
\\ 
\label{betanDelta}
\text{et \ \  \ }
\beta_{n}(\Delta)=\vi_{D}\dotss  \tu{m-1}\vi_{D}(\tu{m}\beta_{n-m}(\Delta))\cap (\V_{0}\oplus ... \oplus \V_{m-1}).\end{gather}
Or \begin{gather*}\beta_{n-m}(\Delta)\subset p^{-C}
\gamma_{n-m}(\Delta)\otimes _{W}\gS[[\frac{u^{e}}{p}]] \\
 \subset p^{-C-C'} \gamma_{0}(\Delta)\otimes_{W} \gS[[\frac{u^{e}}{p}]]
=p^{-C-C'}\beta_{0}(\Delta) [[\frac{u^{e}}{p}]]\end{gather*} grâce à  la propriété d) du lemme~\ref{propbetagamma} et à la condition  i) de la proposition~\ref{prop2} (en notant $C'$ la constante qui y apparaît). 
D'où $$\tu{m}\beta_{n-m}(\Delta) \subset p^{-C-C'}\ \tu{m}\beta_{0}(\Delta)[[\frac{u^{ep^{m}}}{p}]]$$
et par \eqref{betanDelta} on a donc 
\begin{gather}\label{incl-betanDelta2}
\beta_{n}(\Delta)\subset p^{-C-C'}\M[[\frac{u^{ep^{m}}}{p}]]
\cap (\V_{0}\oplus ... \oplus \V_{m-1}).
\end{gather}
Grâce à \eqref{betamDelta} et \eqref{incl-betanDelta2}
  le lemme~\ref{perte-modules}  implique 
$$\beta_{n}(\Delta)\subset p^{-C-C'-C_{1}} \beta_{m}(\Delta)[[\frac{u^{ep^{m}}}{p}]].$$
 Donc on a établi   l'inégalité de droite de (\ref{form-comparaison-beta1}), et un argument  de déterminants permet d'en  déduire
l'inégalité de gauche. \cqfd

\noindent{\bf Suite de la démonstration de 
  ``i) implique ii)'' dans la proposition~\ref{prop2}.}  On va appliquer le lemme suivant à $\gM_{n}=\beta_{n}(\Delta)[[\frac{u^{ep^{n}}}{p}]]$. On rappelle que  $\A=K_{0}[[u]]$  est le complété de 
$\gS[\frac{1}{p}]$ pour la topologie $u$-adique.  

\begin{lem}\label{lem-Mn-Sigma[[]]}
Soit $P$ un $\A$-module libre de rang  $r$ et pour tout $n\in \N$ soit $\gM_{n}$ un  $\gS[[\frac{u^{ep^{n}}}{p}]]$-module libre de rang $r$ muni d'un isomorphisme $\gM_{n}\otimes_{\gS[[\frac{u^{ep^{n}}}{p}]]}\A=P$. On suppose qu'il existe $C\in \N$ tel que pour tout $n\geq m$ on ait 
\begin{gather}\label{cond-Mn-Mm}
p^{C}\gM_{m}\subset \gM_{n}\otimes_{\gS[[\frac{u^{ep^{n}}}{p}]]}\gS[[\frac{u^{ep^{m}}}{p}]] \subset p^{-C}\gM_{m}\text{ \ dans \ } P. 
\end{gather}
Alors il existe un unique $\gS[\frac{1}{p}]$-module $\gN$ libre de rang $r$  muni d'un isomorphisme $\gN\otimes_{\gS[\frac{1}{p}]}\A=P$ tel que pour tout  $\gS$-réseau $\gM$ dans $\gN$
 \begin{gather}\label{cond-M-Mm}
\text{    il existe  }C'\in \N \text{   tel que }
p^{C'}\gM_{m}\subset \gM[[\frac{u^{ep^{m}}}{p}]] \subset p^{-C'}\gM_{m} 
\text{  pour tout } m\in \N. 
\end{gather}
On a de plus 
\begin{gather}\label{formule-gN}
\gN=\{x\in P, \exists k\in \Z, \forall m\in \N, 
x\in p^{k}\gM_{m}\}. 
\end{gather}
\end{lem}
\noindent {\bf Début de la démonstration du  lemme~\ref{lem-Mn-Sigma[[]]}.} Cela résultera  du lemme suivant, appliqué à $Q_{m}=\gM_{m}\otimes_{\gS[[\frac{u^{ep^{m}}}{p}]]}\gS/u^{ep^{m}}\gS$. On remarque que 
 $\A/u^{k}\A=(\gS/u^{k}\gS)[\frac{1}{p}]$ pour tout $k\in \N$.

 \begin{lem}\label{N-module-libre} 
 Soit $P$ un $\A$-module libre de rang $r$, et pour tout $m\in \N$ soit $Q_{m}$ un sous-$(\gS/u^{ep^{m}}\gS)$-module libre de rang $r$ de $P/u^{ep^{m}}P$ tel que 
 $Q_{m}\otimes _{\gS/u^{ep^{m}}\gS}\A/u^{ep^{m}}\A=P/u^{ep^{m}}P$.
 On suppose qu'il existe $C\in \N$, tel que pour $m,n\in \N$, avec $n\geq m$,  en notant $\big(Q_{n}\ \text{mod}\ u^{ep^{m}}\big)$ l'image de $Q_{n}$ dans $P/u^{ep^{m}}P$, on ait 
\begin{gather}\label{cond-N-module-libre}p^{C} Q_{m}
\subset \big(Q_{n}\ \text{mod}\ u^{ep^{m}}\big)
\subset 
p^{-C} Q_{m}.\end{gather}
 Alors \begin{gather}\label{def-N-2eme}\gN=\{x\in P, \exists k\in \Z, \forall m\in \N, 
x\ \text{mod}\ u^{ep^{m}}\in p^{k}Q_{m}\}\end{gather}
 est un $\gS[\frac{1}{p}]$-module libre de rang $r$ muni naturellement d'un isomorphisme  
 $\gN\otimes _{\gS[\frac{1}{p}]}\A=P$
 et si $\gM$ est un $\gS$-réseau de $\gN$,  \begin{gather}\label{condition-C'}\text{il  existe \ }C'\text{\   tel que pour tout }m\in \N,\ \  
 p^{C'} Q_{m}
\subset \big(\gM\ \text{mod}\ u^{ep^{m}}\big)
\subset 
p^{-C'} Q_{m}.\end{gather} 
 \end{lem}
\begin{rem} L'hypothèse que $Q_{m}$ est libre est bien sûr nécessaire car si on prenait $P=\A$ et $Q_{m}$ engendré par $1,\frac{u^{e}}{p}, \big(\frac{u^{e}}{p}\big)^{2},...,\big(\frac{u^{e}}{p}\big)^{p^{m}-1}$ la conclusion du lemme  serait  fausse. \end{rem}
 
\noindent {\bf Début de la démonstration du lemme~\ref{N-module-libre}.} 
Pour tout $n \geq m \geq 0$ on note $$ Q_{m}^{n}=
\sum _{i\geq n}p^{C}Q_{i}\ \text{mod}\ u^{ep^{m}}. $$
Alors $(Q_{m}^{n})_{n\geq m}$ est une suite décroissante de sous-$(\gS/u^{ep^{m}}\gS)$-modules de $Q_{m}$ contenant  $p^{2C}Q_{m}$. 
 Elle est stationnaire, car $(Q_{m}^{n}/p^{2C}Q_{m})_{n\geq m}$ est une suite décroissante de sous-modules du $(\gS/u^{ep^{m}}\gS+p^{2C}\gS)$-module de longueur finie $Q_{m}/p^{2C}Q_{m}$. On note $\wt Q_{m}$ la limite de la suite stationnaire 
$(Q_{m}^{n})_{n\geq m}$. 
Donc $\wt Q_{m}$ est un sous-$(\gS/u^{ep^{m}}\gS)$-module de $Q_{m}$, on a $p^{2C}Q_{m}\subset \wt Q_{m} \subset Q_{m}$, et pour $n\geq m$ on a une surjection $\wt Q_{n}\to \wt Q_{m}$. 
 Soit   $(x_{1}^{0},\dotss  ,x_{r}^{0})$ une base du $(\gS/u^{e}\gS)$-module libre $p^{2C}Q_{0}$ (qui est inclus dans $\wt Q_{0}$). Par récurrence sur $n$ on choisit 
    $x_{1}^{n}\in \wt Q_{n},\dotss  ,x_{r}^{n}\in \wt Q_{n}$ des relèvements de $x_{1}^{n-1},\dotss  ,x_{r}^{n-1}$ de $\wt Q_{n-1}$ à $\wt Q_{n}$. 
 Pour tout $i=1,\dotss  ,r$ la limite projective des $x_{i}^{n}$ définit un élément $x_{i}\in P$. On a même $x_{i}\in \gN$.

 On note $P^{*}=\mathrm{Hom}_{\A}(P,\A)$ et pour  tout $m\in \N$ on note  
 $$Q_{m}^{*}=\mathrm{Hom}_{\gS/u^{ep^{m}}\gS}(Q_{m},\gS/u^{ep^{m}}\gS).$$ Comme  
 $P^{*}$ et le système des $Q_{m}^{*}$ vérifient les conditions du lemme    on peut  leur associer par la construction précédente $r$ éléments 
 $\eta_{1},\dotss  ,\eta_{r}\in P^{*}$ tels que pour $i\in \{1,...,r\}$ et  $n\in \N$ la réduction  $\eta_{i}^{n}$ de  $\eta_{i}$ modulo $u^{ep^{n}}$ appartienne à $Q^{*}_{n}$ et que $\eta_{1}^{0},\dotss  ,\eta_{r}^{0}$ forment une base de $p^{2C}Q_{0}^{*}$. Pour $n\in \N$ et $i,j\in \{1,...,r\}$, on a $\eta_{i}^{n}\in 
 Q_{n}^{*}$ et $x_{j}^{n}\in Q_{n}$, d'où $\s{\eta_{i},x_{j}}\in \gS/u^{ep^{m}}\gS$. En passant à la limite on obtient $\s{\eta_{i},x_{j}}\in \gS$. Donc 
 $a=\det(\s{\eta_{i},x_{j}})$ est un élément  de $\gS$ dont la réduction modulo $u$ est non nulle (et appartient même à $p^{4Cr}W^{\times }$). Pour tout $x\in \gN$ on a $\s{\eta_{i},x}\in \gS[\frac{1}{p}]$. Donc 
 on a un diagramme commutatif  
$$\begin{matrix}(\gS[\frac{1}{p}])^{r}&\hookrightarrow& \gN&\to& (\gS[\frac{1}{p}])^{r}\\
\downarrow&&\downarrow&&\downarrow \\
\A^{r} &\xrightarrow{\sim}&P&\xrightarrow{\sim}&\A^{r}\end{matrix}
$$ 
où les applications horizontales sont 
données par $(\alpha_{1},\dotss  ,\alpha_{r})\mapsto \alpha_{1}x_{1}+\dotss  +\alpha_{r}x_{r}$ et 
$x\mapsto (\s{\eta_{1},x},\dotss  ,\s{\eta_{r},x})$ 
et où les flèches verticales sont les inclusions évidentes. Il résulte du diagramme que l'application $\gN\to (\gS[\frac{1}{p}])^{r}$ est aussi une injection. 
Donc 
 $\gN$ est un $\gS[\frac{1}{p}]$-module 
 de type fini sans torsion et il est 
 libre de rang $r$ grâce au  lemme suivant et au théorème de structure des modules de type fini sur les anneaux principaux.  

\begin{lem}\label{fracz-reflexif}
L'anneau $\gS[\frac{1}{p}]$ est principal. 
\end{lem}
 \noindent {\bf Démonstration.}
 L'anneau $\gS$ est intègre, noethérien et intégralement clos dans son corps des fractions ; il en donc de même de $\gS[\frac{1}{p}]$. Les idéaux premiers non nuls de $\gS$ sont 
 \begin{itemize}
 \item les idéaux principaux $(f)$, où $f$ est un élément irréductible de $\gS$, 
 \item l'idéal maximal $(p,u)$. 
 \end{itemize}
 Par conséquent les idéaux premiers non nuls de 
 $\gS[\frac{1}{p}]$ sont maximaux. Il en résulte déjà que  $\gS[\frac{1}{p}]$ est un anneau de Dedekind (voir par exemple~\cite{serre}, chapitre 1, proposition 4). De plus les idéaux premiers de $\gS[\frac{1}{p}]$ sont principaux, comme on vient de le voir. Donc $\gS[\frac{1}{p}]$ est principal (ibidem, chapitre 1, proposition 7). \cqfd 
 
  \noindent {\bf Fin de la démonstration du lemme~\ref{N-module-libre}.}
 Il reste à montrer (\ref{condition-C'}). 
 Notons $\gM_{0}$ le $\gS$-module libre engendré par $x_{1},...,x_{r}$. 
 Soit $\wt Q$ la limite projective des $\wt Q_{n}$, c'est-à-dire 
 $$\wt Q=\{x\in P, \forall m\in \N, 
x\ \text{mod}\ u^{ep^{m}}\in \wt Q_{m}\}. $$
Alors $\wt Q$  est un sous-$\gS$-module de $P$ qui contient 
$x_{1},...,x_{r}$. 
Comme $\wt Q_{m}\subset Q_{m}$ et $\eta_{i}^{m}\in Q_{m}^{*}$ pour tout $m\in \N$, on a $\s{\eta_{i},x}\in \gS$ pour $i\in \{1,...,r\}$ et $x\in \wt Q$, d'où
 $\wt Q\subset a^{-1}\gM_{0}$.  Comme $\gS$ est noethérien, $\wt Q$ est donc un $\gS$-module de type fini. 
 On a bien sûr 
\begin{gather*}\{x\in P, \forall m\in \N, 
x\ \text{mod}\ u^{ep^{m}}\in p^{2C} Q_{m}\}\subset \wt Q
\\
\subset \{x\in P, \forall m\in \N, 
x\ \text{mod}\ u^{ep^{m}}\in  Q_{m}\},\end{gather*}
donc $\wt Q$ engendre le $\gS[\frac{1}{p}]$-module $\gN$. 
Soit $\gM$ un $\gS$-module   libre de rang $r$ tel que 
 $\gM\otimes_{\gS}\gS[\frac{1}{p}]=\gN$. Comme $\gM$ et $\wt Q$ sont des $\gS$-modules de type fini qui engendrent $\gN$, il existe 
  $C'$ tel que  $p^{C'}\gM\subset \wt Q\subset p^{-C'}\gM$. 
  Pour tout $m\in \N$ on a  
  $$p^{2C} Q_{m}\subset \big(\wt Q \ \text{mod}\ u^{ep^{m}}\big)=\wt Q_{m}\subset  Q_{m}$$ et donc   
  $p^{2C+C'} Q_{m} \subset \big( \gM\ \text{mod}\ u^{ep^{m}}\big)\subset p^{-C'} Q_{m}$. 
 \cqfd
 
 \noindent {\bf Fin de la démonstration du  lemme~\ref{lem-Mn-Sigma[[]]}.} 
 Comme cela a été annoncé, on applique le lemme~\ref{N-module-libre} à 
  $$Q_{m}=\gM_{m}\otimes_{\gS[[\frac{u^{ep^{m}}}{p}]]}\gS/u^{ep^{m}}\gS.$$ 
 On a le droit d'appliquer le lemme~\ref{N-module-libre} car en tensorisant \eqref{cond-Mn-Mm} par $\gS/u^{ep^{m}}\gS$ au-dessus de $\gS[[\frac{u^{ep^{m}}}{p}]]$ on obtient \eqref{cond-N-module-libre}. Soient $\gN$ et $\gM$ comme dans le lemme~\ref{N-module-libre}. Pour $n\geq m$ on note, jusqu'à la fin de la démonstration du  lemme~\ref{lem-Mn-Sigma[[]]},  $\gS_{m,n}$ l'image de 
 $\gS[[\frac{u^{ep^{m}}}{p}]]$ dans 
 $\big(\gS/u^{ep^{n}} \gS\big)[\frac{1}{p}]$. 
  Grâce à l'hypothèse \eqref{cond-Mn-Mm}, 
  il existe $C$ tel que pour  $n\geq m$ on a 
   \begin{gather}\label{cond-Mn-Mm-comp}
p^{C}\gM_{m}\otimes_{\gS[[\frac{u^{ep^{m}}}{p}]]}\gS_{m,n}\subset \gM_{n}\otimes_{\gS[[\frac{u^{ep^{n}}}{p}]]}\gS_{m,n} \subset p^{-C}\gM_{m}\otimes_{\gS[[\frac{u^{ep^{m}}}{p}]]}\gS_{m,n}.
\end{gather}
Comme $\gM_{n}\otimes_{\gS[[\frac{u^{ep^{n}}}{p}]]}\gS_{m,n}= Q_{n}\otimes_{  \gS/u^{ep^{n}} \gS}\gS_{m,n}$ 
 il résulte alors de \eqref{condition-C'} (avec $n$ au lieu de $m$) qu'il existe $C$ tel que pour  $n\geq m$ on a  
  \begin{gather}\label{condition-C'-comp} 
  p^{C}\gM_{m}\otimes_{\gS[[\frac{u^{ep^{m}}}{p}]]}\gS_{m,n}\subset \gM \otimes_{  \gS}\gS_{m,n} \subset p^{-C}\gM_{m}\otimes_{\gS[[\frac{u^{ep^{m}}}{p}]]}\gS_{m,n}.
  \end{gather} 
En fixant $m$ et en faisant tendre $n$ vers l'infini on en déduit \eqref{cond-M-Mm} car $\gS[[\frac{u^{ep^{m}}}{p}]]=\varprojlim \gS_{m,n}$. 
Enfin \eqref{formule-gN}  résulte de \eqref{cond-M-Mm}. 
  \cqfd 
 
\noindent{\bf Fin de la démonstration de ``i) implique ii)'' dans la proposition~\ref{prop2}.}
Grâce à \eqref{form-comparaison-beta1}  on peut appliquer le 
 lemme~\ref{lem-Mn-Sigma[[]]}    à $\gM_{n}=\beta_{n}(\Delta)[[\frac{u^{ep^{n}}}{p}]]$ et donc 
\begin{gather} \label{construc-gN} 
\gN=\{x\in D\otimes_{K_0}\A, \exists k\in \Z, \forall m\in \N, 
x\in p^{k}\beta_{m}(\Delta)[[\frac{u^{ep^{m}}}{p}]]\}
\end{gather}
est un $\gS[\frac{1}{p}]$-module libre de rang $r$. 
On note que d'après la démonstration du lemme~\ref{lem-Mn-Sigma[[]]} on a aussi 
\begin{gather*}
\gN=\{x\in D\otimes_{K_0}\A, \exists k\in \Z, \forall m\in \N, 
x\ \text{mod}\ u^{ep^{m}}\in p^{k}\beta_{m}(\Delta)\ \text{mod}\ u^{ep^{m}}\}
\end{gather*} mais on n'utilisera pas cette formule. 

 On définit $\vi_{\gN}^{-1}: \gN\to \ta \gN$ comme la restriction  de $$\vi_{D}^{-1}:  D\otimes_{K_0}\A \to \ta D\otimes_{K_0}\A$$
 (a posteriori comme $V_{D}\subset U_{D}$,  $\vi_{\gN}^{-1}$ n'a pas de dénominateur en $E$, la matrice $\xi^{-1}$ associée à $\gN$ est à coefficients dans $\O$ au lieu de $\O[\frac{1}{\lambda}]$ et  
 $\gN\subset D\otimes_{K_0}\O\subset D\otimes_{K_0}\A$). 
 
 Il nous reste à montrer que  $\vi_{\gN}^{-1}:\gN\to \ta \gN$ est bien défini  et détermine un  isomorphisme de $\gS[\frac{1}{p},\frac{1}{E}]$-modules $\vi_\gN:  \ta
\gN[\frac{1}{E}]\to \gN[\frac{1}{E}]$, que $(\gN,\vi_{\gN})$ est un pseudo-iso-$\vi / \gS$-module et  que $D=(D,\vi_{D},V_{D})$ est le $\vi$-module de Hodge-Pink associé.

Montrons que $\vi_{D}^{-1}:  D\otimes_{K_0}\A \to \ta D\otimes_{K_0}\A$ envoie $ \gN$ dans $\ta \gN$ et que son image contient $E^{-s}\ \ta \gN$ où $s\in \Z_{\leq 0}$ est tel que $V_{D}\supset E^{-s}U_{D}$. 
Pour tout $n\in \N$ on a 
$$E^{-s} \vi_{D}(\ta \beta_{n}(\Delta))
\subset
\beta_{n+1}(\Delta)\subset \vi_{D}(\ta \beta_{n}(\Delta)).$$ Grâce à \eqref{cond-M-Mm}  on en déduit  
$$ E^{-s} \vi_{D}(\ta \gN)\subset \gN\subset 
 \vi_{D}(\ta \gN).$$
 Donc $\vi_{\gN}^{-1}:\gN\to \ta \gN$ est bien défini comme la restriction de $\vi_{D}^{-1}$  à $ \gN$, et $\vi_{\gN}^{-1}(\gN)\supset E^{-s}\ \ta \gN$. On montre maintenant que $(\gN,\vi_{\gN})$ satisfait la condition (PIM).  
On a 
\begin{gather}\label{est-incl-m-n1}\vi_{D}\ta\vi_{D}\dotss  \tu{n-1}\vi_{D}(\tu{n}\beta_{m}(\Delta))\subset 
\big(E\dotss  \vi^{n-1}(E)\big)^{s}\beta_{m+n}(\Delta) \\ \label{est-incl-m-n2}\text{et\ \ }\beta_{m+n}(\Delta)\subset 
\vi_{D}\ta\vi_{D}\dotss  \tu{n-1}\vi_{D}(\tu{n}\beta_{m}(\Delta)). \end{gather}
Soit $\gM$ un $\gS$-module libre de rang $r$ muni d'un isomorphisme $\gM[\frac{1}{p}]=\gN$. D'après \eqref{cond-M-Mm} il  existe $C'\in \N$ tel que pour tout $l,m\in \N$, 
\begin{gather}\label{est-incl-m-n-l}p^{C'}\gM[[\frac{u^{ep^{m}}}{p}]]\subset
\beta_{m+l}(\Delta)[[\frac{u^{ep^{m}}}{p}]]\subset
p^{-C'}\gM[[\frac{u^{ep^{m}}}{p}]]. \end{gather} En appliquant (\ref{est-incl-m-n-l})  à $l=0$ et à  $l=n$ on   déduit de (\ref{est-incl-m-n1}) que, pour  $m,n\in \N$,
$$p^{2C'}\big(E\dotss  \vi^{n-1}(E)\big)^{-s}\vi_{\gN}\ta\vi_{\gN}\dotss  \tu{n-1}\vi_{\gN}(\tu{n}\gM)[[\frac{u^{ep^{m}}}{p}]]\subset \gM[[\frac{u^{ep^{m}}}{p}]],$$
et donc pour tout $n\in \N$, 
\begin{gather}\label{est-phi-phi-M1}p^{2C'}\big(E\dotss  \vi^{n-1}(E)\big)^{-s}\vi_{\gN}\ta\vi_{\gN}\dotss  \tu{n-1}\vi_{\gN}(\tu{n}\gM)\subset \gM.\end{gather}
De la même fa\c con on déduit de  (\ref{est-incl-m-n2}) que, pour tout $m,n\in \N$,
$$p^{2C'}\gM[[\frac{u^{ep^{m}}}{p}]]\subset \vi_{\gN}\ta\vi_{\gN}\dotss  \tu{n-1}\vi_{\gN}(\tu{n}\gM)[[\frac{u^{ep^{m}}}{p}]]$$ 
et donc pour tout $n\in \N$, 
\begin{gather}\label{est-phi-phi-M2}p^{2C'} \gM\subset \vi_{\gN}\ta\vi_{\gN}\dotss  \tu{n-1}\vi_{\gN}(\tu{n}\gM).\end{gather}
Grâce à (\ref{est-phi-phi-M1}) et (\ref{est-phi-phi-M2}), la condition (PIM) est satisfaite et donc 
  $(\gN,\vi_{\gN})$ est un pseudo-iso-$\vi / \gS$-module. 
 
 Montrons maintenant que $(D,\vi_{D},V_{D})$ est associé à $\gN$. D'abord $(D,\vi_{D})$ est la réduction de $(\gN,\vi_{\gN})$ modulo $u$. 
Le raisonnement établissant l'unicité de la matrice $\xi$ 
dans le lemme~\ref{rigidite2}
montre que l'inclusion de $\gN$ dans $D\otimes _{K_0}\A$ est $\xi^{-1}$. 
 
 Il reste à démontrer que $V_{D}$ est bien associé à $\gN$,  c'est-à-dire que $V_D=V$ 
  où $V=\vi_{D}^{-1}(\gN\otimes_{\gS[\frac{1}{p}]} \wS)$ est la structure de Hodge-Pink associée à $(\gN,\vi_{\gN})$ (puisque $\gN$ désigne ici $\xi^{-1}(\gN)\subset D\otimes _{K_0}\A$).
On choisit   $m>0 $. 
D'après \eqref{cond-M-Mm}  on a 
 \begin{gather} \label{encadrement-N-verif-VD-HP}
  \gN\otimes _{\gS[\frac{1}{p}]}\gS[[\frac{u^{ep^{m}}}{p}]][\frac{1}{p}]=
\beta_{m}(\Delta)  \otimes_{\gS} \gS[[\frac{u^{ep^{m}}}{p}]][\frac{1}{p}].
  \end{gather}
  Or on a  une inclusion 
  $\gS[[\frac{u^{ep^{m}}}{p}]][\frac{1}{p}]\subset \wS$.
  D'après le lemme~\ref{langton-zet},   
  $\vi_{D}^{-1}(\beta_{m}(\Delta))\otimes_{\gS}\wS= V_{D}$ et on déduit alors de  (\ref{encadrement-N-verif-VD-HP}) que $V=V_{D}$.     Ceci termine la démonstration de la proposition~\ref{prop2}. 
\cqfd

 \noindent {\bf Début de la démonstration de la proposition~\ref{prop1}.}
 Dans la suite de ce paragraphe, on a besoin de la convention suivante. Si $k'$ est un corps parfait contenant $\Fp$ (qui sera en fait $k$ ou une extension finie de $\Fp$), on appelle $\vi $-module  sur $k'$ un couple  
$D'=(D',\vi_{D'})$, où 
$ D'$ est un $W(k')[\frac{1}{p}]$-espace vectoriel de dimension finie et  $\vi_{D'}:\ta D'\to D'$ est un isomorphisme de $W(k')[\frac{1}{p}]$-espaces vectoriels, et bien sûr on continue à appeler $\vi$-module un  
 $\vi$-module sur $k$. On note que $W(\Fp)=\Z_p$. 
 
Grâce à l'hypothèse que $k$ est  algébriquement clos,   
on  suppose que $(D,\vi_{D})$ provient d'un $\vi$-module $(\underline  D,\vi_{\underline  D})$ sur $\Fp$ qui est    décent au sens de~\cite{RZ}. 

Si $\Delta$ et $\Delta '$ sont deux $W$-réseaux de $D$, on note $l_{1}(\Delta,\Delta'),\dotss  ,l_{r}(\Delta,\Delta')$ les diviseurs élémentaires de $\Delta'$ par rapport à $\Delta$ : ce sont les uniques entiers relatifs tels que $l_{1}(\Delta,\Delta')\geq \dotss  \geq 
l_{r}(\Delta,\Delta')$ et que pour une certaine base $e_{1},\dotss  ,e_{r}$ de $\Delta$ sur $W$, on ait $$\Delta'=
p^{-l_{1}(\Delta,\Delta')}We_{1}+\dotss  +p^{-l_{r}(\Delta,\Delta')}We_{r}. $$

Dans la suite on soulignera d'un trait les $\vi$-modules sur $\Fp$ et de deux traits les  $\vi$-modules sur une extension finie  $\mathbb{F}_{p^{t}}$ de $\Fp$. On notera $\Z_{p^{t}}= W(\mathbb{F}_{p^{t}}) $ et $\QQ_{p^{t}}=W(\mathbb{F}_{p^{t}})[\frac{1}{p}]  $

\begin{lem}\label{etape1}
Soit  $(D,\vi_{D},V_{D})$ un $\vi$-module de Hodge-Pink   tel que $D=\underline   D\otimes _{\QQ_p}K_{0}$  et $\vi_{D}=\vi_{\underline   D}\otimes 1$ pour un certain $\vi$-module décent $(\underline  D,\vi_{\underline  D})$ sur $\Fp$. Supposons 
$V_{D}\subset U_{D}$, $t_{H}(D)=t_{N}(D)$ et pour tout sous-$\vi$-module $D'$ de dimension $1$, $t_{H}(D')<t_{N}(D')$. 
Soit $\underline  \Delta$ un $\Z_p$-réseau de $\underline  D$ et $\Delta=
\underline  \Delta \otimes _{\Z_p}W$. Soit $a\in \N$. 
Il existe une constante $C_{0}$ telle que la propriété suivante soit vraie. 
Pour tout  $n\in \N^{*}$, il existe une constante $C$ telle que 
 pour tout réseau  $\Delta '$ de $D$ vérifiant   $$\vi_{D}(\ta \Delta')\subset p^{-a}\Delta'\text{ \  et \ }l_{1}(\Delta,\Delta')\geq l_{2}(\Delta,\Delta')+C,$$ 
on ait 
$l_{1}(\Delta, \gamma_{n}(\Delta'))\leq l_{1}(\Delta,\Delta')-n+C_{0}$. 
\end{lem}
 \noindent {\bf Démonstration.}
%Quitte à augmenter $a$ on peut supposer 
%$\vi_{D}(\ta \Delta)\subset p^{-a}\Delta$.  
D'après  le lemme 2.18 de~\cite{RZ} et la proposition 1.6 de~\cite{RZ-landvogt},  il existe des entiers $t, \tilde a, C_{1}$ (ne dépendant que de $D,\vi_{D}$ et $a$) tels que, pour tout  $W$-réseau $\Delta'$ de $D$ vérifiant   $\vi_{D}(\ta \Delta')\subset p^{-a}\Delta'$, il existe un $W$-réseau $\Delta''$ de $D$ défini sur $\mathbb{F}_{p^{t}}$ (c'est-à-dire provenant d'un $\Z_{p^{t}}$-réseau de $\underline  D\otimes_{\QQ_p}\QQ_{p^{t}}$) et vérifiant   $$\vi_{D}(\ta \Delta'')\subset p^{-\tilde a}\Delta''\text{ \  et \  
}\Delta'\subset \Delta'' \subset p^{-C_{1}}\Delta'.  $$
En fait, en examinant les démonstrations de~\cite{RZ} et~\cite{RZ-landvogt}, on voit qu'on pourrait prendre $\tilde a=a$, mais on n'en a pas besoin ici.  
On a alors \begin{gather*}l_{i}(\Delta,\Delta')\leq l_{i}(\Delta,\Delta'')\leq l_{i}(\Delta,\Delta')+C_{1} \\   \text{et \ \ }l_{i}(\Delta,\gamma_{n}(\Delta'))\leq l_{i}(\Delta,\gamma_{n}(\Delta''))\leq l_{i}(\Delta,\gamma_{n}(\Delta'))+C_{1}\end{gather*} 
pour tout $i$. Quitte à remplacer $C_{0}$ par $C_{0}+2C_{1}$, $C$ par $C+2C_{1}$ et $a$ par $\tilde a$, il 
suffit donc de montrer le lemme~\ref{etape1}  en supposant de plus $\Delta '$ défini sur $\mathbb{F}_{p^{t}}$. Le lemme suivant montre qu'avec cette condition supplémentaire l'énoncé du  lemme~\ref{etape1} est vrai 
avec  $C_{0}=0$. On est donc ramené à montrer le lemme suivant. \cqfd 

\begin{lem}\label{etape1-interm}
Soient $(D,\vi_{D},V_{D})$, $(\underline  D,\vi_{\underline  D})$, $\underline  \Delta$ et $\Delta$ comme dans le lemme précédent. 
Soit $a\in \N$ et $t,n\in \N^{*}$. On note  $B(t,a)$ l'ensemble des réseaux 
$\Delta '$ de $D$ qui sont définis sur $\mathbb F_{p^{t}}$ (c'est-à-dire  associés à des 
$\Z_{p^{t}}$-réseaux $\underline{\underline \Delta} '$ de $\underline  D\otimes_{\QQ_p}\QQ_{p^{t}}$) et vérifient 
$\vi_{D}(\ta \Delta')\subset p^{-a}\Delta'$. 
Il existe une constante $C$ telle que 
 \begin{gather*}\text{pour tout   \ }\Delta '\in B(t,a)\text{ \  
 vérifiant \ }l_{1}(\Delta,\Delta')\geq l_{2}(\Delta,\Delta')+C, \\ 
\text{ \ on ait 
\ }l_{1}(\Delta, \gamma_{n}(\Delta'))\leq l_{1}(\Delta,\Delta')-n.\end{gather*}  
\end{lem}
 \noindent {\bf Démonstration.}
 Si l'énoncé du lemme~\ref{etape1-interm} est faux 
  il existe une suite de réseaux $\Delta_{m}$ dans l'ensemble $B(t,a)$ vérifiant 
$$l_{1}(\Delta,\gamma_{n}(\Delta_{m}))>l_{1}(\Delta,\Delta_{m})-n$$ et telle que 
 $l_{1}(\Delta,\Delta_{m})-l_{2}(\Delta,\Delta_{m})$ tende vers l'infini.
 On note $\underline{\underline  \Delta}_{m}$ le $\Z_{p^{t}}$-réseau de $\underline  D\otimes_{\QQ_p}\QQ_{p^{t}}$ tel que $ \Delta_{m}= \underline{\underline  \Delta}_{m}\otimes_{\Z_{p^{t}}}W$. En notant $a_{m}=l_{1}(\Delta,\Delta_{m})-l_{2}(\Delta,\Delta_{m})$ on voit que $$\big(p^{l_{1}(\Delta,\Delta_{m})}{\underline{\underline  \Delta}}_{m} +(p^{a_{m}}\underline  \Delta\otimes_{\Z_p}\Z_{p^{t}})\big)/(p^{a_{m}}\underline  \Delta\otimes_{\Z_p}\Z_{p^{t}})$$ est un sous-$(\Z_{p^{t}}/p^{a_{m}})$-module libre de rang 1 de
 $\underline  \Delta\otimes_{\Z_p}(\Z_{p^{t}}/p^{a_{m}})$ et définit donc un  point de $\mathbb P^{r-1}(\Z_{p^{t}}/p^{a_{m}})$. Par compacité de $\mathbb P^{r-1}(\Z_{p^{t}})$ muni de la topologie $p$-adique, on peut extraire une sous-suite ayant une limite et on note $\underline{\underline  D}'$ le sous-$\QQ_{p^{t}}$-espace vectoriel de dimension $1$ de $\underline  D\otimes_{\QQ_p}\QQ_{p^{t}}$ tel que cette limite soit
 $\underline{\underline  D}' \cap (\underline  \Delta\otimes_{\Z_p}\Z_{p^{t}})$, 
  sous-$\Z_{p^{t}}$-module libre de rang $1$ de $\underline  \Delta\otimes_{\Z_p}\Z_{p^{t}}$. On pose alors  $D'= \underline{\underline  D}'\otimes _{\QQ_{p^{t}}}K_{0}$ qui est un sous-$K_0$-espace vectoriel de dimension $1$ de $D$. Autrement dit, quitte à extraire une sous-suite il existe une  suite d'entiers $b_{m}$ tendant vers $+\infty$ quand $m$ tend vers l'infini telle que $$p^{l_{1}(\Delta,\Delta_{m})}\Delta_{m}+p^{b_{m}}\Delta=(D' \cap \Delta )+p^{b_{m}}\Delta\text{ \  pour tout }m.$$ 
Comme les réseaux $\Delta_{m}$ appartiennent à $B(t,a)$, en passant à la limite on voit que 
$D'$ est un sous-$\vi$-module de $D$. Par construction $D'$ est  défini sur 
 $\mathbb{F}_{p^{t}}$ mais nous n'avons plus besoin de cette information. Le lemme~\ref{grand-dim1} ci-dessous montre que  si $m$ est assez grand, 
$$l_{1}(\Delta,\gamma_{n}(\Delta_{m}))\leq l_{1}(\Delta,\Delta_{m})+n(t_{H}(D')- t_{N}(D')). $$
 On obtient ainsi un sous-$\vi$-module  $D'$ de dimension $1$ de $D$
tel que $t_{H}(D')\geq t_{N}(D')$, ce qui contredit l'hypothèse. Le lemme~\ref{etape1-interm} est donc ramené au lemme suivant. \cqfd

\begin{lem}\label{grand-dim1} 
Soit $(D,\vi_{D},V_{D})$
 un $\vi$-module de Hodge-Pink avec $V_{D}\subset U_{D}$ et $D'=(D',\vi_{D'},V_{D'})$ un sous-$\vi$-module de Hodge-Pink de dimension $1$. Soit $\Delta\subset D$ un $W$-réseau.  
Soit $x$ un générateur de $D'$. Pour tout $n\in \N$, il existe $C\in \N$ tel que pour $r\in \N$ assez grand, 
$$\gamma_{n}(Wx+p^{r}\Delta)\subset p^{n(t_{N}(D')- t_{H}(D'))}Wx+p^{r-C}\Delta.$$
\end{lem}
 \noindent {\bf Démonstration.}
On rappelle que  $V_{D'}= D'\otimes_{K_0}\wS[\frac{1}{E}]\cap V_{D}$. On note  $D''=D/D'$ avec $\vi_{D''}$ induit par $\vi_{D}$ et $V_{D''}$  l'image de $V_{D}$ dans $D''\otimes_{K_0}\wS[\frac{1}{E}]$, de sorte que  $0\to D'\to D\to D''\to 0$ est une suite 
exacte courte dans la catégorie des $\vi$-modules de Hodge-Pink. On note   $\beta_{n}'$, $\beta_{n}''$, $\gamma_{n}'$ et $\gamma_{n}''$ les applications correspondant à $\beta_{n}$ et $\gamma_{n}$ pour $D'$ et $D''$. On note $\Delta ''$ l'image de $\Delta$ dans $D''$. Soit $n\in \N$. 
Il existe $r_{0}$ tel que 
$\beta_{n}(p^{-r_{0}}Wx+\Delta)\to \beta_{n}''(\Delta'')$ soit surjectif : on relève les vecteurs d'une base du $\gS$-module libre de rang $r$ 
\begin{gather*}\beta_{n}''(\Delta'')=((\vi_{D''}\ta \vi_{D''} \dotss  {}\tu{n-1}\vi_{D''}\tu{n}\Delta'') \otimes_{W}\gS)
\\
\cap\big( \vi_{D''}\dotss  \tu{n-1}\vi_{D''}\tu{n-1}V_{D''}
\oplus \dotss  
\oplus\vi_{D''}V_{D''}\big)\end{gather*} en des vecteurs de $(\vi_{D}\ta \vi_{D} \dotss  {}\tu{n-1}\vi_{D}\tu{n}\Delta) \otimes_{W}\gS$ que l'on corrige par des vecteurs de $(\vi_{D'}\ta \vi_{D'} \dotss  {}\tu{n-1}\vi_{D'}\tu{n}D') \otimes_{K_0}\gS[\frac{1}{p}]$ pour qu'ils appartiennent à $\vi_{D}\dotss  \tu{n-1}\vi_{D}\tu{n-1}V_{D} \oplus \dotss  \oplus 
\vi_{D}V_{D}$, ce qui est possible, puisque 
\begin{gather*}\gS[\frac{1}{p}]\to \wS_{0}/E^{-s}\wS_{0}\oplus \dotss   \oplus 
\wS_{n-1}/\vi^{n-1}(E)^{-s}\wS_{n-1}\end{gather*} est surjectif (où $s\in \Z_{\leq 0}$ est tel que $E^{-s}U_{D}\subset V_{D}\subset U_{D}$). 
Soit  $r\geq r_{0}$. 
Comme $\beta_{n}(p^{-r_{0}}Wx+\Delta)\to \beta_{n}''(\Delta'')$ est surjectif de noyau $\beta_{n}'(p^{-r_{0}}Wx) $ et 
 $\beta_{n}(p^{-r}Wx+\Delta)\to \beta_{n}''(\Delta'')$ est surjectif de noyau $\beta_{n}'(p^{-r}Wx) $ 
 on a 
\begin{gather*}\beta_{n}(Wx+p^{r}\Delta)=\beta_{n}'(Wx)
+p^{r-r_{0}}\beta_{n}(Wx+p^{r_{0}}\Delta),  \\
\text{ donc \  \ }
\gamma_{n}(Wx+p^{r}\Delta)=\gamma_{n}'(Wx)
+p^{r-r_{0}}\gamma_{n}(Wx+p^{r_{0}}\Delta).\end{gather*} Comme 
$\gamma_{n}'(Wx)=p^{n(t_{N}(D')- t_{H}(D'))}Wx
$, le lemme \ref{grand-dim1}  est démontré, ce qui achève aussi la démonstration des lemmes~\ref{etape1} et~\ref{etape1-interm}. 
 \cqfd

 \noindent {\bf Fin de la démonstration de   la proposition~\ref{prop1}. } On rappelle que $k$ est supposé algébriquement clos. Soit $D=(D,\vi_{D},V_{D})$ un $\vi$-module de Hodge-Pink  qui est un objet irréductible de $MHP(\vi)_{fa}$ et tel que 
$V_{D}\subset U_{D}$.

 Soit $a$ la constante du b) du lemme~\ref{propbetagamma}. On a $\vi_{D}(\ta \gamma_{m}(\Delta))\subset p^{-a}\gamma_{m}(\Delta)$ pour tout $m\in  \N$. Soit  $i\in \{1,\dotss  ,r-1\}$. 
 Le lemme~\ref{etape1}, appliqué à $\Lambda^{i}D$   montre l'existence de $C_{0}\in \N$ tel que pour tout 
 $n\in \N$ il existe $C_{1}$ tel que 
 si $\Delta'$ est un réseau de $D$ vérifiant  $$\vi_{D}(\ta \Delta')\subset p^{-a}\Delta'\text{\ \   et \ }l_{i}(\Delta, \Delta')\geq l_{i+1}(\Delta, \Delta')+C_{1}, $$alors 
 $(l_{1}+\dotss  +l_{i})(\Delta, \gamma_{n}(\Delta'))\leq (l_{1}+\dotss  +l_{i})(\Delta,\Delta')-n+C_{0}$. 

En notant $C_{2}$ 
 la constante qui apparaît dans 
 e) du lemme~\ref{propbetagamma} on a donc, pour tout  $m\in  \N$ tel que $l_{i}(\Delta, \gamma_{m}(\Delta))\geq l_{i+1}(\Delta, \gamma_{m}(\Delta))+C_{1}$, 
 \begin{gather}\label{llli}(l_{1}+\dotss  +l_{i})(\Delta, \gamma_{n+m}(\Delta))\leq (l_{1}+\dotss  +l_{i})(\Delta,\gamma_{m}(\Delta))-n+C_{0}+iC_{2}.\end{gather}
 
 Comme $\det(\gamma_{m}(\Delta))=\det (\Delta)$ pour tout $m\in \N$, on a $$(l_{1}+\dotss  +l_{r})(\Delta,\gamma_{m}(\Delta))=0.$$ 
On a  $p^{a}\gamma_{m}(\Delta)\subset \gamma_{m+1}(\Delta)\subset p^{-a}\gamma_{m}(\Delta)$ pour tout $m\in \N$. 
 En appliquant l'inégalité~(\ref{llli}) avec $n=C_{0}+(r-1)C_{2}$ et $i\in \{1,...,r-1\}$, et en appliquant le lemme~\ref{suite-polygones} ci-dessous avec $$C=\max(C_{1},na)\text{ \  et \ } l_{i}^{m}=l_{i}(\Delta,\gamma_{mn}(\Delta))\text{ \ pour \ }i\in \{1,...,r\}\text{ \ et \ }m\in \N,$$ on voit qu'il existe $C_{3}$ tel que  \begin{gather*}p^{C_{3}}\Delta\subset  \gamma_{mn}(\Delta)\subset p^{-C_{3}}\Delta \text{ pour tout }m\in \N. \end{gather*} On a donc 
 $p^{C_{3}+[\frac{m}{2}]a}\Delta\subset  \gamma_{m}(\Delta)\subset p^{-C_{3}-[\frac{m}{2}]a}\Delta$ pour tout $m\in \N$. On est donc ramené à montrer le  lemme~\ref{suite-polygones}. \cqfd 
 
 Pour  donner un sens  géométrique au lemme suivant, on signale que si $\Delta'$ est un réseau, le polygône concave ayant pour sommets
 $(i,(l_{1}+...+l_{i})(\Delta,\Delta'))$ pour $i=0,...,r$ peut être appelé ``polygône des invariants'' de $\Delta'$ relativement à $\Delta$ et que les $l_{i}(\Delta,\Delta')$  sont les pentes de ce polygône. Dans le cas que l'on considère on a toujours $(l_{1}+...+l_{r})(\Delta,\Delta')=0$, c'est-à-dire que 
 le premier sommet est $(0,0)$ et le dernier est $(r,0)$.

  \begin{lem}\label{suite-polygones}
  Soit $r\in \N^{*}$ et $C\in \R_{+}$. Soit $l_{i}^{n}\in \R$ pour $i\in \{1,...,r\}$ et $n\in \N$ vérifiant
  \begin{itemize}
  \item i) pour tout $n\in \N$,  $l_{1}^{n}\geq ...\geq l_{r}^{n}$, 
  \item ii) pour tout $n\in \N$, $\sum_{i=1}^{r}l_{i}^{n}=0$, 
  \item iii) on a $l_{1}^{0}=...=l_{r}^{0}=0$, 
  \item iv) pour $i\in \{1,...,r\}$ et $n\in \N$,  $l_{i}^{n}-C\leq l_{i}^{n+1}\leq l_{i}^{n}+C$, 
  \item v) pour $i\in \{1,...,r-1\}$ et $n\in \N$, si $l_{i}^{n}-l_{i+1}^{n}\geq C$ on a 
  $$l_{1}^{n+1}+...+l_{i}^{n+1}\leq l_{1}^{n}+...+l_{i}^{n}.$$ 
  \end{itemize}
  Alors on a 
  $l_{i}^{n}\in [-C',C']$ pour tout $i\in \{1,...,r\}$ et pour tout $n\in \N$, avec $C'$ ne dépendant que de $r$ et de $C$. 
      \end{lem}
    \noindent{\bf Démonstration.} Pour tout $n\in \N$ on note $P_{n}$ le polygône concave dont les sommets sont les $(i,P_{n}(i))$ pour $i\in \{0,...,r\}$, avec 
  $P_{n}(i)=l_{1}^{n}+...
+l_{i}^{n}$ (en particulier le premier sommet est $(0,0)$ et le dernier est $(r,0)$). La différence des pentes au sommet $(i,P_{n}(i))$, que nous appellerons brisure,  est $l_{i}^{n}-l_{i+1}^{n}$. 
La condition v) assure que   $P_{n+1}(i)\leq P_{n}(i)$  si la brisure de $P_{n}$  en  $(i,P_{n}(i))$ est $\geq C$.  L'idée naïve de la démonstration est la suivante : si $\max_{i}P_{n}(i)$ est atteint en un sommet $(j,P_{n}(j))$ de brisure $\geq C$ la condition v) impose $P_{n+1}(j)\leq P_{n}(j)$ et 
si $\max_{i}P_{n+1}(i)$ est atteint pour la même valeur de $i$ (c'est-à-dire $i=j$) on en déduit $\max_{i}P_{n+1}(i)\leq \max_{i}P_{n}(i)$, ce qui impose à  $P_{n}$ de rester borné de fa\c con uniforme en $n$. Nous allons voir qu'en rempla\c cant $P_{n}$ par $P_{n}-Q$ pour un certain $Q$ de brisure constante $3C$ cet argument 
devient correct. 
On pose donc $Q(i)=\frac{3C}{2}i(r-i)$ pour $i\in\{0,...,r\}$. 
   On pose ensuite $$m_{n}=\max_{i\in \{0,...,r\}}(P_{n}(i)-Q(i))\in \R_{+}. $$ Nous allons montrer que $m_{n+1}\leq m_{n}$ pour tout $n\in \N$. Comme $m_{n}\geq0$,  il n'y a rien à démontrer si $m_{n+1}=0$. Supposons donc $m_{n+1}>0$ et soit $j\in \{1,...,r-1\}$ tel que le maximum soit atteint en $j$. 
  La brisure de $P_{n+1}-Q$ en $j$ est donc $\geq 0$. Or la brisure de $Q$ est égale à $3C$, donc $l_{j}^{n+1}-l_{j+1}^{n+1}\geq 3C$. Grâce à la condition iv) on en déduit $  l_{j}^{n}-l_{j+1}^{n}\geq C$, d'où par la condition v), $P_{n}(j)\geq P_{n+1}(j)$. On a donc 
  $$m_{n}\geq P_{n}(j)-Q(j)\geq P_{n+1}(j)-Q(j)=m_{n+1}.$$ Par la condition iii) on a $m_{0}=0$. 
  Donc $m_{n}=0$ pour tout $n\in \N$. Donc \begin{gather*}l_{n}^{1}=P_{n}(1)\leq Q(1)=\frac{3}{2}C(r-1) \\ \text{ et \ }l_{n}^{r}=-P_{n}(r-1)\geq -Q(r-1)=-\frac{3}{2}C(r-1)\text{ \   pour tout }n\in \N.\end{gather*} 
  Ceci termine la démonstration du lemme~\ref{suite-polygones} (avec $C'=\frac{3}{2}C(r-1)$) et donc celle  de la proposition~\ref{prop1}.    \cqfd

 \section{Théorie entière à la Breuil}   \label{fontainelaffaille1}
 
  On note $\mathcal S$ le complété $p$-adique de l'enveloppe à puissance divisées 
 $$\gS[\frac{u^{ei}}{i!}]_{i\in \N^{*}}=\gS[\frac{u^{ep}}{p}, \frac{u^{ep^{2}}}{p^{p+1}}, \frac{u^{ep^{3}}}{p^{p^{2}+p+1}},\dotss  ].$$ Cet anneau a été introduit par Faltings~\cite{faltings} et Breuil~\cite{breuil-invent99,breuil}. 
 
 Nous allons donner une autre démonstration du théorème 2.2.1 de~\cite{carliu} (voir le théorème~\ref{breuil-2261} ci-dessous). 
 %Ce résultat est indépendant de la condition \eqref{cond-gr} de transversalité de Griffiths. Cependant, lorsque cette condition est vérifiée, 
 %  ce résultat est relié aux conjectures de Breuil en théorie entière (voir 
 % la conjecture 2.2.6 (1) dans~\cite{breuil}). En fait cette conjecture de Breuil  a été démontrée par Tong Liu (théorème 2.3.5 de~\cite{liu}), en utilisant les résultats de Kisin~\cite{kisin}, après des résultas partiels de Breuil et Caruso~\cite{breuil-compositio,breuil,caruso}. Le théorème 2.2.1 de~\cite{carliu} est postérieur à la preuve de Liu,  et sert dans~\cite{carliu} comme résultat préliminaire pour étudier des objets de torsion.  

On fixe  $m\in \{1,..., p-2\}$. 
 On rappelle que $\mathrm{Mod}^{\vi}_{/\gS,[0,m]} $ désigne la catégorie des $\vi / \gS$-modules $(\gM,\vi_{\gM})$  d'amplitude $\subset [0,m]$, c'est-à-dire tels que $\vi_{\gM}$ et $E^{m}\vi_{\gM}^{-1}$ n'aient pas de dénominateurs en $E$. 
Soit $\mathrm{Mod}^{\vi}_{/\mathcal S,[0,m]}$ la catégorie des $(D,\vi_{D},V_{D},\M)$ où 
\begin{itemize}
\item 
$(D,\vi_{D})$ est un $\vi$-module 
\item $V_{D}$ est une structure de Hodge-Pink  telle que 
 $U_{D}\subset V_{D}\subset E^{-m}U_{D}$, 
\item  $\M$ est un $\mathcal S$-module libre, muni d'un isomorphisme 
\begin{gather}\label{isom-M-D} \M \otimes _{\mathcal S}\mathcal S[\frac{1}{p}]\simeq \ta D\otimes _{K_0}\mathcal S[\frac{1}{p}],\end{gather} 
\item on suppose que  $\M$ est ``fortement divisible'', c'est-à-dire 
  \begin{gather}\label{cond-fort-div}\mathcal S.p^{-m}\ \ta( (\vi_{D}\otimes 1)(\M \cap E^{m}V_{D}))=\M.\end{gather}
 \end{itemize}
 
Comme $\M\cap E^{m}V_{D}\supset E^{m}\M$ et $\frac{\vi(E)}{p}$ est une unité dans $\mathcal S$, l'hypothèse \eqref{cond-fort-div} implique que $ \vi_{D}(\ta \M)\subset \M$ (condition demandée par Breuil dans~\cite{breuil}).

 On note $\mathbb D:\mathrm{Mod}^{\vi}_{/\gS,[0,m]}\to \mathrm{Mod}^{\vi}_{/\mathcal S,[0,m]}$ le foncteur qui à $(\gM,\vi_{\gM})$ associe \break $(D,\vi_{D},V_{D},\mathcal M)$ 
  avec 
$(D,\vi_{D},V_{D})=\mathbb D_{\mathrm{iso}}(\gM[\frac{1}{p}],\vi_{M}\otimes 1)$ et   $\M=\ta \gM\otimes_{\gS}\mathcal S$ qui est un  $\mathcal S$-module libre  muni de l'isomorphisme \eqref{isom-M-D} associé à $\ta \xi^{-1}$.

  Pour montrer que le foncteur $\mathbb D$ est bien défini on doit vérifier que $\M$ est fortement divisible. On  rappelle cette démonstration, due à Breuil.

 L'inclusion $$\mathcal S.p^{-m}\ \ta( (\vi_{D}\otimes 1)(\M \cap E^{m}V_{D}))\supset \M$$ est évidente. En effet 
on a $$\M\cap  E^{m}V_{D}\supset 
\ta \xi^{-1}(\ta \gM)\cap  E^{m}V_{D}\supset 
\ta \xi^{-1}( E^{m}\vi_{\gM}^{-1}(\gM))$$  puisque par définition de la catégorie $\mathrm{Mod}^{\vi}_{/\gS,[0,m]}$, $ E^{m}\vi_{\gM}^{-1}\in \mathrm{Hom}_{\gS}(\gM,{}\ta \gM)$ (en fait 
$\ta \xi^{-1}(\ta \gM)\cap  E^{m}V_{D}=\ta \xi^{-1}( E^{m}\vi_{\gM}^{-1}(\gM))$). Donc 
\begin{gather*}p^{-m}\ \ta ((\vi_{D}\otimes 1)(\M\cap  E^{m}V_{D}))\supset p^{-m}\ \ta((\vi_{D}\otimes 1)(\ta \xi^{-1}( E^{m}\vi_{\gM}^{-1}(\gM))))
\\
=\Big(\frac{\vi(E)}{p}\Big)^{m}\ta \xi^{-1}(\ta \gM)
\text{ puisque } (\vi_{D}\otimes 1)\ta \xi^{-1}\vi_{\gM}^{-1}=\xi^{-1}
.\end{gather*}
 Enfin $\frac{\vi(E)}{p}$ est inversible dans $\mathcal S$. 
Il  reste donc à montrer l'inclusion inverse 
$$p^{-m}\ \ta( (\vi_{D}\otimes 1)(\M \cap E^{m}V_{D}))\subset \M. $$
Comme $\vi_{\gM}$ appartient à $\mathrm{Hom}_{\gS}(\ta \gM,\gM)$, 
et par la définition de $V_{D}$, on a 
\begin{gather*} \xi((\vi_{D}\otimes 1)(\M \cap E^{m}V_{D})) \subset 
(\gM\otimes_{\gS}\mathcal S)\cap 
(E^{m}\gM\otimes_{\gS}\wS)\\ 
=\gM\otimes_{\gS}
(\mathcal S\cap E^{m}\wS). \end{gather*}
Mais  $E^{-m}\mathcal S\cap \wS$  est inclus dans le complété 
    $p$-adique de $$\mathcal S+\frac{u^{e(p-m)}}{p}\mathcal S+\frac{u^{e(p^{2}-m)}}{p^{p+1}}\mathcal S+\dotss ,$$ donc 
$\ta(\mathcal S\cap E^{m}\wS)\subset \vi(E)^{m}\mathcal S$. 
Comme $\frac{\vi(E)}{p}\in \mathcal S$, on a  $$p^{-m}\ \ta(\mathcal S\cap E^{m}\wS)\subset \mathcal S,$$ donc $\M$ est fortement divisible et le foncteur $\mathbb D$ est bien défini.

 Le théorème suivant est le théorème 2.2.1 de~\cite{carliu}. 
 
 \begin{thm}\label{breuil-2261}
 Le foncteur  $\mathbb D:\mathrm{Mod}^{\vi}_{/\gS,[0,m]}\to \mathrm{Mod}^{\vi}_{/\mathcal S,[0,m]}$ 
  est  une équivalence de catégories. 
    \end{thm}
  \begin{rem}  Il résulte du théorème que pour tout objet $(D,\vi_{D},V_{D},\M)$ de $ \mathrm{Mod}^{\vi}_{/\mathcal S,[0,m]}$, $(D,\vi_{D},V_{D})$ est faiblement admissible. \end{rem}

  \noindent {\bf Début de la démonstration du théorème~\ref{breuil-2261}.}
Montrons d'abord  que $\mathbb D$  est pleinement fidèle.
Soient $(\gM,\vi_{\gM})$ et $(\gM',\vi_{\gM'})$ deux objets de $\mathrm{Mod}^{\vi}_{/\gS,[0,m]}$ et $(D,\vi_{D},V_{D},\M)$ et $(D',\vi_{D'},V_{D'},\M ')$ leurs images dans $\mathrm{Mod}^{\vi}_{/\mathcal S,[0,m]}$. 
L'application \begin{gather*}\mathrm{Hom}_{\mathrm{Mod}^{\vi}_{/\gS,[0,m]}}((\gM,\vi_{\gM}),(\gM',\vi_{\gM'})) \\ 
\to \mathrm{Hom}_{\mathrm{Mod}^{\vi}_{/\mathcal S,[0,m]}}((D,\vi_{D},V_{D},\M),(D',\vi_{D'},V_{D'},\M '))\end{gather*} est injective par la fidélité de $\mathbb D_{\mathrm{iso}}$ et on veut montrer qu'elle est surjective. 
Soit $$h\in \mathrm{Hom}_{\mathrm{Mod}^{\vi}_{/\mathcal S,[0,m]}}((D,\vi_{D},V_{D},\M),(D',\vi_{D'},V_{D'},\M ')).$$ 
Par la pleine fidélité de $\mathbb D_{\mathrm{iso}}$, on sait qu'il existe un unique morphisme $f$ de pseudo-iso-$\vi / \gS$-modules de  $\gM[\frac{1}{p}]$ dans $\gM'[\frac{1}{p}]$
 qui vérifie 
$f\vi_\gM=\vi_\gM'\ta f$
et dont la réduction modulo $u$ est $h$. L'hypothèse $(h\otimes 1)(\M)\subset \M'$ implique que $f\otimes_{\gS[\frac{1}{p}]}1_{\mathcal S[\frac{1}{p}]}$ envoie $\ta \gM\otimes_{\gS}\mathcal S$ dans $\ta \gM'\otimes_{\gS}\mathcal S$. 
Il s'agit de montrer que $f$ appartient  à 
$\mathrm{Hom}_{\gS}(\gM,\gM')$. 

Choisissons des bases de $\gM$ et $\gM'$ sur $\gS$ et notons $r$ et $r'$ les rangs de $\gM$ et $\gM'$. 
On a donc 
 $f\in M_{r'r}(\gS[\frac{1}{p}])$ et $\ta f\in  M_{r'r}(\mathcal S)$ et il s'agit de montrer $f\in M_{r'r}(\gS)$. Comme $\gS/u^{ep}\gS$ est un quotient de $\mathcal S$, l'hypothèse $\ta f\in  M_{r'r}(\mathcal S)$ implique que  $ f$ modulo $u^e$ appartient à $M_{r'r}(\gS/u^{e}\gS)$. 

On a $f=\vi_{\gM'}{}\ta f\vi_\gM^{-1}$. Supposons par l'absurde que $f$ n'appartient pas à $M_{r'r}(\gS)$. Soit $a\in  \N$ le petit entier tel que $p^{a}f$ appartienne à 
$M_{r'r}(\gS)$. Par hypothèse $a>0$. Soit $b\in \N$ le plus grand entier tel que 
$p^{a}f\text{ mod }p \in M_{r'r}(k[[u]])$ appartienne à 
$u^{b} M_{r'r}(k[[u]])$. Par hypothèse, $b\geq  e$. On a l'égalité suivante, où toutes les matrices entre parenthèses ont leurs coefficients dans $\gS$ :
$$E^{m}(p^{a}f)=(\vi_{\gM'})(p^{a}\ {}\ta f)(E^{m}\vi_\gM^{-1}).$$
En réduisant cette égalité modulo $p$, on trouve l'égalité suivante, où toutes les matrices entre parenthèses sont à coefficients dans $k[[u]]$ : 
\begin{gather*}u^{em}\Big((p^{a}f) \text{ mod }p\Big) \\ =
\Big(\vi_{\gM'} \text{ mod }p\Big)
\Big((p^{a}\ {}\ta f) \text{ mod }p\Big)
\Big(E^{m}\vi_\gM^{-1}\text{ mod }p\Big).\end{gather*}
Or $(p^{a}f) \text{ mod }p$ est divisible exactement par $u^{b}$ et $(p^{a}\ {}\ta f) \text{ mod }p$ est divisible (exactement) par $u^{pb}$. Le membre de gauche est divisible exactement par $u^{b+em}$, alors que le membre de droite est divisible par 
$u^{pb}$. Comme $b\geq e$ et $p\geq m+2$ on a 
$b(p-m)>e$ ce qui amène une contradiction.  

   Il nous reste à montrer que  le foncteur de $\mathrm{Mod}^{\vi}_{/\gS,[0,m]}$ dans $\mathrm{Mod}^{\vi}_{/\mathcal S,[0,m]}$ est essentiellement surjectif. 
  On doit construire un $\gS$-module libre $\gM$ de rang $r$, à partir de la donnée d'un $\mathcal S$-module libre $\M$ de rang $r$. L'idée est de construire une suite de modules libres de rang $r$ sur des anneaux qui se rapprochent de plus en plus de $\gS$ et de définir $\gM$ comme la limite de cette suite. Il y a plusieurs choix possibles pour cette suite d'anneaux. La difficulté de cette démonstration n'est pas de construire un $\gS$ module $\gM$, mais de montrer qu'il est libre de rang $r$. Il est important de noter la parenté entre la preuve qui suit et la démonstration du théorème~\ref{fa--a} donnée au paragraphe~\ref{faibl-ad=ad}.

      Le lemme suivant généralise  le lemme~\ref{langton-zet}  (qui correspond au cas où $\CC=\CC'=\gS$).    
  
  \begin{lem}\label{Langtonapproche-S} 
  Soient $m\in \N^{*}$ et $\CC$ et $\CC'$ des anneaux tels que 
  $\gS\subset \CC\subset \CC'\subset \wS$ et 
  \begin{gather}\label{cond-CC-CC'}
  E^{-m}\CC\cap \wS\subset \CC'. 
  \end{gather}
    a) 
  Soit $\NN$ un $\CC$-module libre de rang fini et $\W$ un sous-$\wS$-module de $\NN\otimes _{\CC} 
  \wS$ contenant 
  $E^{m}\NN\otimes _{\CC} 
  \wS$. Alors 
  $\NN'=\CC'.(\NN\cap \W)$ est un 
  sous-$\CC'$-module libre de  
  $\NN\otimes _{\CC}\CC'$, 
   \begin{gather*}\NN'[\frac{1}{E}]=\NN\otimes _{\CC}\CC'[\frac{1}{E}]
  \text{ \ \ et  \ \ } \NN'\otimes _{\CC'} \wS=\W.\end{gather*}
 De plus $\det(\NN')=E^{k}\det(\NN)\otimes_{\CC} \CC'$ 
  où  $k$ désigne la longueur du $\wS$-module $(\NN\otimes _{\CC} 
  \wS)/\W$.
  
  b) Soit de plus $\gN$ un $\gS$-module libre muni d'un isomorphisme $$\gN\otimes _{\gS}
 \CC=\mathcal N. $$ Alors l'inclusion 
 $\CC'. (\gN\cap \W)\subset \CC'. (\mathcal N\cap \W)$ est une égalité. 

 c) Soit de plus 
 $\wt \CC$ un  anneau tel que 
  $\gS\subset \wt \CC\subset \CC$ et 
  $\wt{\mathcal N}$ un $\wt \CC$-module libre muni d'un isomorphisme  $\wt{\mathcal N}\otimes _{\wt \CC}
 \CC=\mathcal N$. Alors l'inclusion 
 $$\CC'. (\wt{\mathcal N}\cap \W)\subset \CC'. (\mathcal N\cap \W)$$ est une égalité. 
  \end{lem}
   \noindent {\bf Démonstration.} 
   On montre simultanément a) et b). Soit $\gN$ un $\gS$-module libre muni d'un isomorphisme  $\gN\otimes _{\gS}
 \CC=\mathcal N$.
     D'après le lemme~\ref{langton-zet}, $\gN'=\gN\cap \W$ est un $\gS$-module libre, et on a $$E^{m}\gN\subset \gN'\subset \gN, \ \  \det(\gN')=E^{k}\det(\gN)\text{ \ et \ }\W=\gN'\otimes _{\gS}\wS.$$ 
  Comme  $\NN=\gN\otimes _{\gS}\CC$ il résulte de ce qui précède que  $$\NN\cap \W\subset \gN'\otimes _{\gS}(E^{-m}\CC\cap \wS).$$  
   Donc, grâce à \eqref{cond-CC-CC'},  $$\NN'\subset \gN'\otimes _{\gS}\CC'$$ et on a l'égalité puisque $\NN\cap \W$ contient évidemment $\gN'$. 
 De plus $$\det(\NN')=\det(\gN')\otimes _{\gS}\CC'.  $$
 Enfin c) découle de b), en choisissant $\gN$ tel que $\wt\NN=\gN \otimes_{\gS}\wt\CC$. 
 \cqfd
 
 Dans la suite, pour $a\in \N^{*}$ on note 
 $$\CC_{a}=\gS[[\frac{u^{ea}}{p}]].$$
 On appliquera le lemme~\ref{Langtonapproche-S} dans les deux situations suivantes
 \begin{itemize}
 \item
 $\CC=\mathcal S$ et $\CC'=\CC_{p-m}$,  
 \item
 pour $a>m$, $\CC=\CC_{a}$ et $\CC'=\CC_{a-m}$. 
 \end{itemize}
 Dans ces deux situations, la condition \eqref{cond-CC-CC'} est vérifiée. Voici la justification dans le premier cas : $(E^{-m}\mathcal S\cap \wS)$ est inclus dans le complété $p$-adique de  $\gS[\frac{u^{e(p-m)}}{p}, \frac{u^{e(p^{2}-m)}}{p^{p+1}}, \frac{u^{e(p^{3}-m)}}{p^{p^{2}+p+1}},\dotss  ]$. Cet anneau est inclus dans $\CC_{p-m}$ car $\frac{p^{i+1}-m}{p^{i}+...+p+1}\geq p-m$ pour tout $i\geq0$, puisque $m\geq 1$.

 \noindent{\bf Fin de la démonstration du théorème~\ref{breuil-2261}.} 
 Montrons maintenant que le foncteur  $\mathbb D:\mathrm{Mod}^{\vi}_{/\gS,[0,m]}\to \mathrm{Mod}^{\vi}_{/\mathcal S,[0,m]}$ est essentiellement surjectif. 
 Soit $(D,\vi_{D},V_{D},\M)$ un objet de $\mathrm{Mod}^{\vi}_{/\mathcal S,[0,m]}$.
 
  D'après le a) du lemme~\ref{Langtonapproche-S} 
 (appliqué à $\CC=\mathcal S$, $\CC'=\CC_{p-m}$ et $\mathcal N=\M$), $$\CC_{p-m}.(\M\cap E^{m}V_{D})$$ est un $\CC_{p-m} $-module libre de rang $r$. On note $\gM_{0}$ son image par $E^{-m}\vi_{D}\otimes 1$, si bien que 
 $$\gM_{0}=\CC_{p-m}.E^{-m}(\vi_{D}\otimes 1)(\M\cap E^{m}V_{D})
 $$ est un $\CC_{p-m}$-module libre de rang $r$ naturellement inclus dans $$E^{-m}D\otimes _{K_0}\CC_{p-m}[\frac{1}{p}].$$ Alors $\ta \gM_{0}$ est un 
 $\CC_{p(p-m)}$-module libre de rang $r$ qui est  naturellement inclus dans $$\vi(E)^{-m}\ \ta D\otimes _{K_0}\CC_{p(p-m)}[\frac{1}{p}]. $$
 Par l'hypothèse de forte divisibilité de $\M$, et comme 
 $\frac{\vi(E)}{p}$ est une unité dans $\mathcal S$, on a 
  \begin{gather}\label{relation-gM0-calM}\ta \gM_{0}\otimes_{\CC_{p(p-m)}}\mathcal S=\M.\end{gather}
 
  On définit par récurrence $\gM_{n}$ pour $n>0$ en posant 
\begin{gather}\label{def-rec-Mn}\gM_{n}=\CC_{m_{n}}.E^{-m}(\vi_{D}\otimes 1)(\ta \gM_{n-1}\cap E^{m}V_{D}),\end{gather} où la suite $m_{n}$ est définie par $$m_{0}=p-m \text{ \ \ et \ \ } m_{n}=pm_{n-1}-m$$ (comme $m\leq p-2$, cette suite
est strictement croissante et tend vers l'infini). 
 
 Par le a) du  lemme~\ref{Langtonapproche-S} (appliqué à $\CC=\CC_{pm_{n-1}}$, $\CC'=\CC_{m_{n}}$ et $\mathcal N=\ta \gM_{n-1}$), $\gM_{n}$ est un $\CC_{m_{n}}$-module libre de rang $r$,  inclus dans $E^{-m}\dotss  \vi^n(E)^{-m} D\otimes _{K_0}\CC_{m_{n}}[\frac{1}{p}]$. 
 %L'hypothèse de forte divisibilité de $\M$, et le c) du lemme~\ref{Langtonapproche-S} 
 %(appliqué à $a=p, b=pm_{n-1}, \mathcal N=\mathcal N$ et $\wt{\mathcal N}=\ta \gM_{n-1}$) 
 %montrent, par récurrence sur $n$,  
 %$\ta \gM_{n}\otimes_{\CC_{pm_{n}}}\mathcal S=\M$. 

D'après \eqref{relation-gM0-calM} on a  $\ta \gM_{0}\subset \M$. On en déduit 
  $\gM_{1}\subset \gM_{0}$. Cette assertion a un sens car les deux sont plongés dans $E^{-m}D\otimes _{K_0}\CC_{p-m}[\frac{1}{p}]$. On en déduit, pour tout 
 $n\in \N$, $\gM_{n+1}\subset \gM_{n}$. 
 En fait \eqref{relation-gM0-calM} implique, grâce au c) du lemme~\ref{Langtonapproche-S} (appliqué à $\CC=\mathcal S$, $\CC'=\CC_{m_{0}}$, $\wt\CC=\CC_{pm_{0}}$,  $\mathcal N=\mathcal M$,   $\wt\NN=\ta \gM_{0}$), que l'inclusion $\gM_{1}\subset \gM_{0}$ induit une égalité 
 $$\gM_{1}\otimes _{\CC_{m_{1}}}\CC_{m_{0}}=\gM_{0}.$$ 
 On montre alors, par récurrence sur $n$, grâce au c)  du lemme~\ref{Langtonapproche-S} (appliqué à $\CC=\CC_{pm_{n-1}}$, $\CC'=\CC_{m_{n}}$, $\wt\CC=\CC_{pm_{n}}$,  $\mathcal N=\ta \gM_{n-1}$,   $\wt\NN=\ta \gM_{n}$), que l'inclusion $\gM_{n+1}\subset \gM_{n}$ induit une égalité $$\gM_{n+1}\otimes _{\CC_{m_{n+1}}}\CC_{m_{n}}=\gM_{n}.$$
 En particulier comme $m_{n+1}>m_{n}>p^{n}$, l'inclusion 
   $\gM_{n+1}\subset \gM_{n}$ induit  une égalité modulo $u^{ep^{n}}$, 
c'est-à-dire 
$$\gM_{n+1}\otimes_{\CC_{m_{n+1}}}\gS/u^{ep^{n}}= \gM_{n}\otimes_{\CC_{m_{n}}}\gS/u^{ep^{n}}.$$
 
 On prend alors $\gM$ égal à l'intersection des $\gM_{n}$. On remarque que 
 $$Q_{n}=\gM_{n}\otimes_{\CC_{m_{n}}}\gS/u^{ep^{n}}$$ est un $\gS/u^{ep^{n}}$-module libre de rang $r$ et que $$Q_{n}=Q_{n+1}\otimes_{\gS/u^{ep^{n+1}}}\gS/u^{ep^{n}}.$$ On a aussi $\gM=\varprojlim Q_{n}$ donc $\gM$ est un $\gS$-module libre de rang $r$. De plus 
   pour tout $n\in \N$ l'inclusion $\gM\subset \gM_{n}$ induit une égalité $$\gM\otimes_{\gS}\CC_{m_{n}}=\gM_{n}. $$

 D'après \eqref{def-rec-Mn} et grâce au b)  du lemme~\ref{Langtonapproche-S},
  on a 
 $$\gM=E^{-m}(\vi_{D}\otimes 1)(\ta \gM\cap E^{m}V_{D}). $$

 On définit  $E^{m}\vi_{\gM}^{-1}: \gM\to \ta \gM$ comme la restriction de $E^{m}(\vi_{D}\otimes 1)^{-1}$ à $ \gM$. 
On a  alors $\vi_{\gM}\in \mathrm{Hom}_{\gS}(\ta \gM,\gM)$.

Comme $\gM\otimes_{\gS}\CC_{m_{0}}=\gM_{0}$ et 
grâce à  \eqref{relation-gM0-calM} 
on a 
 $\ta \gM\otimes _{\gS}\mathcal S= \M$. 
   \cqfd
   
   \begin{rem}
  On notera la parenté entre la preuve du théorème~\ref{breuil-2261} et  la preuve du théorème~\ref{fa--a} donnée dans le paragraphe~\ref{faibl-ad=ad}.  \end{rem}
 
\section{Un cadre plus général}
% englobant le cas d'inégales caractéristiques} 
  
 On présente ici  un cadre un peu plus général où $\Z_{p}$ est remplacé par l'anneau des entiers d'un corps local non archimédien $\O_{F}$. 
 Ce cadre serait adapté à l'étude des modules $\pi_{F}$-divisibles munis d'une action stricte de $\O_{F}$ au sens de~\cite{faltings-strict}. Il est donc plus restrictif que celui considéré par Kisin dans~\cite{kisin-pss}, où l'action n'est pas nécessairement stricte.

 Les résultats des paragraphes~\ref{pseudo-iso-phi-gS},~\ref{construc-Dieudonne},~\ref{faibl-ad=ad} et~\ref{fontainelaffaille1} s'étendent à ce cadre et y généralisent  ceux de~\cite{fontaine29} (qui correspond au cas où $\O_{F}$ est d'égales caractéristiques) au lieu d'être simplement  analogues. 
 Pour la théorie rationnelle cela n'implique rien de nouveau car le  théorème ``faiblement admissible implique admissible'' pour $\Z_{p}$ implique facilement le théorème à coefficients dans $\O_{F}$ si $F$ est une extension finie de $\QQ_{p}$. En revanche la généralisation à $\O_{F}$ de la théorie entière  du paragraphe~\ref{fontainelaffaille1} n'est pas une conséquence du cas où $\O_{F}=\Z_{p}$. 
 
Les résultats du paragraphe~\ref{sec-griffiths}  ne s'étendent pas de manière agréable à ce nouveau cadre à cause du rôle spécial joué par l'opérateur $N_{D}$. Lorsque  $N_{D}=0$, la condition \eqref{cond-gr} de transversalité de Griffiths est analogue à la condition de tranquillité dans~\cite{fontaine29}. 

Soit $\O_{F}$ l'anneau des entiers d'un corps local non archimédien $F$  de corps résiduel $\F$. On note $\pi_{F}$ une uniformisante de $\O_{F}$. Soit $k$ un corps parfait contenant $\F$. On note $W=W_{\O_{F}}(k)$ où $W_{\O_{F}}(k)$ est le ``$\O_{F}$-anneau de Witt de $k$'', défini par 
\begin{gather*} W_{\O_{F}}(k)=
W(k)\otimes _{W(\F)}\O_{F}\text{ si }\O_{F}\text{ est d'inégales caractéristiques} \\
\text{ et }W_{\O_{F}}(k)=k[[\pi_{F}]]\text{ si }\O_{F}=\F[[\pi_{F}]]. 
\end{gather*} 
 On possède un morphisme de Frobenius $\O_{F}$-linéaire $\vi:W\to W$ qui relève l'endomorphisme $x\mapsto x^{q}$ de $k$.  On note $K_{0}=\mathrm{Frac} W=W[\frac{1}{\pi_{F}}]$. 
On note 
$\gS=W[[u]]$. On note $z\in \gS$ l'image de $\pi_{F}\in W$. 
On note $\vi:\gS\to \gS$ le morphisme égal à $\vi$ sur $W$ et envoyant $u$ sur $u^{q}$. 
 
Soit $E\in \gS$ tel que $E(0)\in zW^{\times}$ et que $E$ ne soit pas multiple de $z$. Soit $\O_{K}=\gS/E\gS$. 
 Le morphisme $W\to \O_{K}$ fait de $\O_{K}$ l'anneau des entiers d'une    extension finie totalement ramifiée $K$ de $K_{0}$. On note 
 $\pi_{K}$ l'image de $u$, qui est alors une uniformisante de $\O_{K}$. Si on écrit $E=E(u)=\sum _{n\in \N}c_{n}u^{n}$ avec $c_{n}\in W$ on a $c_{0}=E(0)\in zW^{\times}$, $c_{1},...,c_{e-1}\in zW$ et $c_{e}\in W^{\times}$ pour un certain  entier $e$ qui est aussi l'indice de ramification de $K$ sur $K_{0}$. Quitte à multiplier $E$ par un élément inversible de $\gS$ (ce qui ne change pas les résultats de cet article) on peut supposer si on veut que $c_{e}=1$ et $c_{e+1}=c_{e+2}=...=0$, c'est-à-dire  que $E$ est 
un polynôme d'Eisenstein, et le  polynôme minimal de $\pi_{K}$  sur $K_{0}$.

Ce formalisme contient les cas particuliers suivants
\begin{itemize}
\item {\bf situation de~\cite{kisin} et de cet article :} on a  $$q=p, \O_{F}=\Z_{p}, \pi_{F}=p, W=W(k), z=p, 
\gS=W[[u]],$$ $K$ est une extension totalement ramifiée de $K_{0}=W[\frac{1}{p}]$, $\pi_{K}$ est une uniformisante de $K$,  
et $E$ est  le polynôme minimal de $\pi_{K}$ sur $K_{0}$, qui est un polynôme d'Eisenstein, 
\item {\bf situation  de~\cite{fontaine29} :} on a \begin{gather*}\O_{F}=\F[[\pi_{F}]], \O_{K}=k[[\pi_{K}]], \gS=\O\tc_{\F}\O_{K}, u=1\otimes \pi_{K},\\  z=\pi_{F}\otimes 1 \text{ \  et \ }E=z-(1\otimes \pi_{F}).\end{gather*} Dans~\cite{fontaine29}, $\O_{F}, \pi_{F}, \O_{K}, \pi_{K}, \lambda$ étaient notés $\O,\pi,\O_{L},\pi_{L}, \alpha$ respectivement. De plus $\O$ était noté $\CC$ et muni d'une topologie légèrement plus fine. 
\end{itemize}

La définition~\ref{def-ModphiS} s'étend de fa\c con évidente à ce nouveau cadre. 
Plus précisément la catégorie $\mathrm{Mod}^{\vi}_{/\gS,[s,t]} $  
 des $\vi / \gS$-modules d'amplitude $\subset [s,t]$  est formé des 
 couples $(\gM, \vi_{\gM})$ où $\gM$ est un $\gS$-module libre de type fini et $\vi_{\gM}: \ta \gM[\frac{1}{E}]\to \gM[\frac{1}{E}]$ est un isomorphisme vérifiant  $E^{t}\ \gM \subset \vi_{\gM}(\ta \gM)\subset E^{s}\ \gM$.
 Dans la situation  de~\cite{fontaine29} où $\O_{F}$ est d'égales caractéristiques, un $\vi / \gS$-module est exactement  un chtouca local sur $\O_{K}$ au sens de~\cite{fontaine29} et l'isomorphisme $\xi$ du
  lemme~\ref{rigidite2} est l'inverse de celui noté  $R$  dans le lemme 7.4    
  de~\cite{fontaine29}.  Dans le cas où $\O_{F}$ est d'inégales caractéristiques, on notera que $E$ est le polynôme d'Eisenstein de $\pi_{K}$ sur $\O_{F}$ et non sur $\Z_{p}$. Cela correspond à la condition d'action stricte mentionnée au début du paragraphe. 

On définit   la  catégorie $MHP(\vi)$ des $\vi$-modules de Hodge-Pink de la même fa\c con que dans l'introduction. Les résultats des paragraphes~\ref{pseudo-iso-phi-gS},~\ref{construc-Dieudonne} et~\ref{faibl-ad=ad} s'étendent de fa\c con évidente, en rempla\c cant $p$ respectivement par 
$\pi_{F}, z$  ou $q$ selon que $p$ est vu comme un élément de 
$\Z_{p}$, de $\gS$ ou comme le cardinal du corps résiduel. En particulier on construit un 
 foncteur de Dieudonné $\mathbb D_{\mathrm{iso}}: \ModphiS \otimes_{\O_{F}}F\to MHP(\vi) $ et on montre qu'il  est pleinement fidèle et que son image essentielle est constituée des $\vi$-modules de Hodge-Pink faiblement admissibles. 
 Dans situation  de~\cite{fontaine29} où $\O_{F}$ est d'égales caractéristiques, ce  théorème  est exactement le théorème 7.3 de~\cite{fontaine29}.  
 
 Enfin on énonce la généralisation du  théorème~\ref{breuil-2261}, car ce n'en est pas une conséquence (mais la preuve est exactement la même).   On note $\mathcal S$ le complété $z$-adique de  
 $\gS[\frac{u^{eq}}{z}, \frac{u^{eq^{2}}}{z^{q+1}}, \frac{u^{eq^{3}}}{z^{q^{2}+q+1}},\dotss  ]$. On fixe  $m\in \{1,..., q-2\}$. 
 Soit $\mathrm{Mod}^{\vi}_{/\mathcal S,[0,m]}$ la catégorie des $(D,\vi_{D},V_{D},\M)$ avec 
$(D,\vi_{D},V_{D})\in  MHP(\vi)$ d'amplitude $\subset [0,m]$, et 
 $\M$  un $\mathcal S$-module libre muni d'un isomorphisme 
$ \M[\frac{1}{z}]\simeq \ta D\otimes _{K_0}\mathcal S[\frac{1}{z}]$ et  ``fortement divisible'' au sens où 
$$\mathcal S.z^{-m}\ \ta( (\vi_{D}\otimes 1)(\M \cap E^{m}V_{D}))=\M.$$
    On note $\mathbb D:\mathrm{Mod}^{\vi}_{/\gS,[0,m]}\to \mathrm{Mod}^{\vi}_{/\mathcal S,[0,m]}$ le foncteur qui à $(\gM,\vi_{\gM})$ associe \break $(D,\vi_{D},V_{D},\mathcal M)$ 
  avec 
$(D,\vi_{D},V_{D})=\mathbb D_{\mathrm{iso}}(\gM[\frac{1}{p}],\vi_{M}\otimes 1)$ et   $\M=\ta \gM\otimes_{\gS}\mathcal S$.

   \begin{thm}\label{breuil-2261-gen}
 Le foncteur  $\mathbb D:\mathrm{Mod}^{\vi}_{/\gS,[0,m]}\to \mathrm{Mod}^{\vi}_{/\mathcal S,[0,m]}$ est bien défini et 
  est  une équivalence de catégories. \cqfd
    \end{thm}

\end{document}